\documentclass{article} 
\usepackage{graphicx}

\usepackage[utf8]{inputenc} 
\usepackage{hyperref}       
\usepackage{url}            
\usepackage{booktabs}       
\usepackage{amsfonts}       
\usepackage{nicefrac}       
\usepackage{microtype}      

\usepackage{amsfonts,amssymb,amsmath,amsbsy}
\usepackage{tabularx}
\usepackage{graphicx}
\usepackage{float}
\usepackage{verbatim}
\usepackage{dsfont}
\usepackage{placeins}
\usepackage{color}
\usepackage{tikz}
\usetikzlibrary{calc,shapes,arrows}
\usepackage{multimedia}
\usepackage{url}
\usepackage{bbm}
\usepackage{latexsym,mathrsfs}
\usepackage{mathtools}
\usepackage{multirow}
\usepackage{algorithm}
\usepackage{algorithmic}
\usepackage{textcomp}
\usepackage{xspace} 

\newcommand{\Popt}{{\overline P}_{\text{opt}}}
\newcommand{\doe}{DoE\xspace}
\newcommand{\initDoEsize}{\texttt{initDoEsize}\xspace}
\newcommand{\budget}{\texttt{budget}\xspace}
\newcommand{\fmin}{\ensuremath{f^\text{min}}\xspace}
\newcommand{\COCO}{\texttt{COCO}\xspace}
\begin{document}

\title{Revisiting Bayesian Optimization in the light of the COCO benchmark
}


\author{Rodolphe Le Riche\thanks{CNRS LIMOS (Mines St-Etienne and UCA), St-Etienne, France, leriche@emse.fr}         \and
        Victor Picheny\thanks{Secondmind, Cambridge, UK, victor@secondmind.ai }
}


\date{
\begin{minipage}{9cm}
{\small(the order of the authors is alphabetical as they have contributed equally to the content of the current study) }
\vskip\baselineskip
July 2021
\vskip\baselineskip
{\small Please cite as: Rodolphe Le Riche and Victor Picheny, \textit{Revisiting Bayesian Optimization in the light of the COCO benchmark}, Journal of Structural and Multidisciplinary Optimization, DOI:10.1007/s00158-021-02977-1, July 2021.}
\end{minipage}
}
\maketitle

\begin{abstract}
It is commonly believed that Bayesian optimization (BO) algorithms are highly efficient for optimizing numerically costly functions. 
However, BO is not often compared to widely different alternatives, and is mostly tested on narrow sets of problems (multimodal, low-dimensional functions), 
which makes it difficult to assess where (or if) they actually achieve state-of-the-art performance.
Moreover, several aspects in the design of these algorithms vary across implementations without a clear recommendation emerging from current practices, and
many of these design choices are not substantiated by authoritative test campaigns.
This article reports a large investigation about the effects on the performance of (Gaussian process based) BO of common and less common design choices. 
The following features are considered: the size of the initial design of experiments, the functional form of the trend, the choice of the kernel, the internal optimization strategy, input or output warping, and the use of the Gaussian process (GP) mean in conjunction with the classical Expected Improvement. 
The experiments are carried out with the established COCO (COmparing Continuous Optimizers) software.
It is found that a small initial budget, a quadratic trend, high-quality optimization of the acquisition criterion bring consistent progress. 
Using the GP mean as an occasional acquisition contributes to a negligible additional improvement.
Warping degrades performance.
The Mat{\'e}rn 5/2 kernel is a good default but it may be surpassed by the exponential kernel on irregular functions.
Overall, the best EGO variants are competitive or improve over state-of-the-art algorithms in dimensions less or equal to 5 for multimodal functions.
The code developed for this study makes the new version (v2.1.1) of the R package DiceOptim available on CRAN.
The structure of the experiments by function groups allows to define priorities for future research on Bayesian optimization. 
\end{abstract}


\section{Introduction}

Bayesian Optimization algorithms (BO) are global optimization methods that iterate by constructing and using conditional Gaussian processes (GP). 
They belong to the larger family of surrogate-assisted optimization algorithms \cite{gramacy2020surrogates} with the distinguishing feature that the GP provides a mathematically funded model for the uncertainty that accompanies unexplored regions of the search space.
The original idea of BO can be traced back to the sixties and seventies with the articles \cite{kushner1964new,mockus1978application} and the reference BO algorithm called EGO (for Efficient Global Optimization) was proposed 20 years later \cite{jones1998efficient} with a scope of application already targeting numerically costly objective functions. 
Since then, new BO algorithms have appeared which differ mainly in the way the GP is used (i.e., in the so-called acquisition criterion) and incidentally in the GP model and in the numerical resolution of internal optimizations. 

BO being a very active field for the past 20 years, there are now a large number of reviews and tutorials, see   \cite{jones2001taxonomy,forrester2008engineering,picheny2013benchmark,shahriari2016taking,frazier2018tutorial,gramacy2020surrogates} for instance. 
There are also a large number of available open-source implementations, such as \texttt{Spearmint} \cite{snoek2012practical}, 
\texttt{DiceOptim} \cite{roustant2012dicekriging}, \texttt{BayesOpt} \cite{martinez2014bayesopt}, 
\texttt{SMAC} \cite{hutter2015manual},
\texttt{GPyOpt} \cite{gonzalez2016gpyopt},
\texttt{GPflowOpt} \cite{knudde2017gpflowopt},
\texttt{RoBO} \cite{klein2017robo},
\texttt{STK} \cite{bect2019stk},
\texttt{Botorch} \cite{balandat2019botorch},
\texttt{SMT} \cite{bouhlel2019python}, among others. 
However, while some design choices and default tunings are shared between those implementations, several fundamental algorithm features like the computing budget devoted to initializing the GP, and several others that we will discuss below, have still not been agreed upon in the community, being generally cast as case-dependent, and their choice being left to the user.

It is also a common claim that BO algorithms are state-of-the-art for costly functions. 
However, most studies substantiate this assertion by comparing BO algorithms with themselves, e.g., by comparing acquisition criteria \cite{de2021greed,rehbach2020expected,picheny2013benchmark} or other BO variants \cite{eggensperger2013towards,salem2019gaussian,spagnol2019global}. 
This can be understood as it allows to make incremental comparisons and draw clear conclusions.
It is also common for studies to be based on test functions that can be well learned by GPs such as the Branin, the Hump Camel or the Hartmann functions \cite{lindauer2019towards}, or even GP predictions (surrogate test cases in \cite{eggensperger2015efficient}, in MOPTA \cite{jones2008large} and in metaNACA \cite{gaudrie2019phd}). However, those do not allow assessing the performance of BO on a wider range of problems, nor how well it compares to other families of algorithms.

The \COCO (COmparing Continuous Optimizers, \cite{hansen2016coco}) software is a recent effort 
to build a testbed that allows the rigorous comparison of optimizers, in general but also when the objective function is expensive \cite{hansen2010comparing}. The noiseless BBOB test suite is available as part of \COCO \cite{bbob_online_functions}. It contains 24 functions that are randomly translated and rotated to make a larger set of tests. The functions are divided into 5 classes (separable, unimodal with moderate conditioning, unimodal with high conditioning, multi-modal with adequate global structure, multi-modal with weak global structure), which eases the analysis of the results.

This study attempts to answer practical and ever-present questions about Bayesian optimization in the light of the \COCO testing program. 
Three goals are pursued.
First, we want to understand what makes BO efficient and answer general questions regarding the choice of the GP kernel, the trend, the initial GP budget, and the suboptimization of the acquisition function. 
Second, we want to find on which function class and dimension BO is relevant when compared to state-of-the-art optimizers for expensive functions.
Through this study, a new version of the BO code \texttt{DiceOptim} was developped and validated. 
The third goal is to expose on the way some tricks of the trade such as how hyperparameters are optimized, local convergence avoided and numerical instabilities handled.


The article starts with a brief introduction to Bayesian optimization,
followed by an explanation about the algorithm's degrees of freedom that will be explored. Then, the implementation of the algorithm is detailed and the testing procedure presented. 
The results, which include comparisons with two state-of-the-art optimizers (\texttt{SMAC} \cite{hutter2015manual} and \texttt{DTA-CMA-ES} \cite{bajer2019gaussian})
conclude the article and are complemented in Appendix by experiments to measure the regression ability of the different GPs.

\section{Basics of Bayesian optimization}

For completeness and to introduce our notations, we now provide basics of BO.
More extensive coverage can be found in e.g., \cite{forrester2008engineering,shahriari2016taking,frazier2018tutorial,gramacy2020surrogates}.
The optimization problem is to best minimize an expensive function $f: \mathcal S \rightarrow \mathbb R$ within a limited budget of $n$ calls. 
The search space $\mathcal S$ is a closed and bounded part of the $d$-dimensional space of real numbers, $\mathbb R^d$. 
In the rest of this paper, $\mathcal S$ is a hyper-rectangle defined by lower and upper bounds, $L, U \in \mathbb R^d$.
BO proceeds by generating a sequence of points $\{x_t\}_{t=1}^T$ at which $f$ is observed, so that $\min_t f(x_t) \rightarrow \min_{x \in \mathcal S} f(x)$ when $T$ is large. The index $t$ stands for time and carries the expensive function assumption that, whatever is done in the BO algorithm, the computing time is dominated by the $t$ calls to $f()$.

In order to save calls to the expensive function $f$, a statistical model of $f$ is built in the form of a Gaussian process (GP), $Y(x)$. 
$Y(x)$ is fully characterized by its trend $\mu(x)$ and its covariance function, also known as kernel, $\text{Cov}(Y(x),Y(x')) = k(x,x')$ and, at any point $x \in \mathcal S$, $Y(x)$ follows a normal law $\mathcal N(\mu(x),k(x,x))$. 
The underlying assumption of the GP model is that there is a trajectory of $Y(x)$ which is equal to $f(x)$. 

At iteration $t$, the observations of the true function $\{x^1,f^1\},\{x^2,f^2\},$\linebreak$\ldots,\{x^t,f^t\}$ (with $f^i = f(x^i)$) are available.
The trend $\mu(.)$ and the kernel $k(.,.)$ depend on parameters that can be learned by maximizing the likelihood or a cross-validation error of the $t$ observations. Learning these parameters is a first nonlinear internal optimization problem, the result of which is a known prior law for $Y(x)$.

The GP is then forced to approximate the observations by conditioning, which yields another GP written 
\begin{equation*}
Y^t(x) \equiv Y(x) \mid \left( Y(x^1) = f(x^1), \ldots , Y(x^t) = f(x^t) \right).
\end{equation*}
The conditional GP follows a posterior law $Y^t(x) \sim \mathcal N(m^t(x),c^t(x,x))$ whose 
trend $m^t()$ and covariance $c^t(,)$ are analytically known from $\mu(.)$ and $k(.,.)$ and the $t$ current observations. 
The formulas for the posterior trend and covariance can be found in \cite{roustant2012dicekriging} or any other textbook on GPs.
The conditioning with the function observations has spawned the adjective ``Bayesian'' in BO. 
At time $t$, $Y^t()$ synthesizes both what is known about the true function, that is the observations, and the uncertainties 
that remain throughout the search space.

The mean and kernel $\mu(.)$ and $k(.,.)$ are generally chosen beforehand from a family of parametric functions. 
Typically, anisotropic stationary kernels are used for $k(.,.)$, with parameters encoding amplitude (kernel variance) and smoothness (lengthscales, also referred to as ranges). $\mu(.)$ is often set to simply zero, or are chosen as polynomials (following the universal kriging framework, \cite{matheron1963principles}).
The parameters of both $\mu(.)$ and $k(.,.)$ are inferred, generally jointly by maximising the model likelihood.

\begin{figure}
\begin{center}
\begin{tikzpicture}
\node[anchor=south west,inner sep=0, outer sep=0] (contour) at (0,0) {\includegraphics[width=0.7\textwidth,trim=0 1.6cm 0 2.cm,clip]{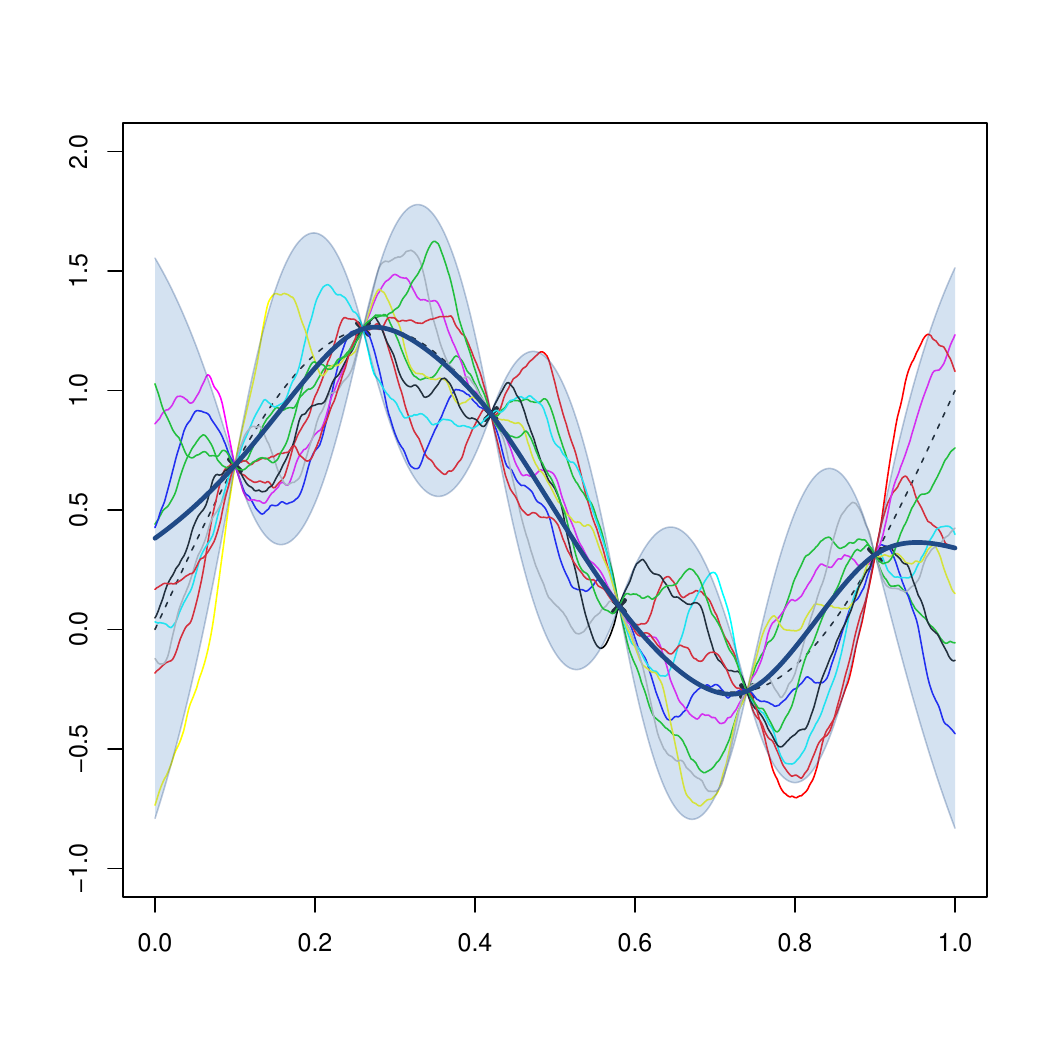}};
\begin{scope}[x={(contour.south east)},y={(contour.north west)}]
\node[] (xlab) at (0.55,-0.05) {$x$};
\node[] (ylab) at (0.,0.6) {$Y^t$ or $f$};
\end{scope}
\end{tikzpicture}
\end{center}
\vspace{-0.5cm}
\caption{
Illustration of a conditional GP model in 1 dimension. Crosses, true function evaluations, i.e., here $t=6$; dashed line, true function $f$ (unknown excepted at the crosses); thin coloured lines, trajectories of the conditional GP; thick blue line, conditional GP mean, $m^t(x)$; blue background, 95\% confidence interval.
\label{fig-explKrig1D}
}
\end{figure}

\vskip\baselineskip
A generic BO method is described in Algorithm~\ref{alg:BO}.
It starts with the creation of a Design of Experiments (\doe), which most often is space-filling. 
The value of the objective function is evaluated at the \doe points. 
Next, an iterative procedure is entered where a conditional GP is learnt from the current \doe, the GP is used to define the next iterate, 
which in turns is evaluated and added to the \doe. 

The current GP represents the beliefs about the true function posterior to the observed data.
The definition of the next iterate from the current GP involves the maximization of an acquisition criterion $a(x;Y^t)$, cf. Equation~(\ref{eq:maxAcquisition}). 
A wide variety of acquisition criteria exist (see \cite{shahriari2015taking,picheny2013benchmark}). 
They all make a tradeoff between the intensification of the search in already known good regions of the search space and the exploration of new regions. 
In this article, only the celebrated Expected Improvement criterion (EI) \cite{jones1998efficient} is considered, 
\begin{equation}
\begin{split}
a(x&;Y^t)  = EI^t(x) = \\ 
& (\fmin - m^t(x)) \Phi\bigg(\frac{\fmin - m^t(x)}{\sqrt{c^t(x,x)}}\bigg) 
+ \sqrt{c^t(x,x)} \phi\bigg( \frac{\fmin - m^t(x)}{\sqrt{c^t(x,x)}} \bigg) ,
\end{split}
\label{eq:EI}
\end{equation}
where \fmin is the best observed true function, and $\phi$ and $\Phi$ denote respectively the probability and cumulative density functions of the standard Gaussian distribution.  
When the acquisition criterion is the Expected Improvement, the BO algorithm is called EGO for Efficient Global Optimization. 
Once the GP mean and covariance is learnt, 
calculating $EI^t$ is relatively inexpensive as it does not imply any call to the true (expensive) function. 
However, $EI^t$ is a multimodal criterion, so solving Problem~(\ref{eq:maxAcquisition}) requires a global optimization 
algorithm. This is further discussed in Section~\ref{sec:acqfuncoptim}.

A BO algorithm classically stops when the budget, i.e., the maximum number of true function evaluations, is exhausted. 
Other stopping criteria (not considered here) are possible, for example a lower bound on the Expected Improvement or a maximum number of function evaluations without progress.
In the implementation discussed in this article, stopping also occurs after repeated iteration failures (either at the GP building or the acquisition maximization steps), cf. Section~\ref{sec:implementation}.
\begin{algorithm}
\caption{Bayesian Optimization}
\label{alg:BO}
\begin{algorithmic}
\REQUIRE a function $f$, a \budget, an \initDoEsize $<$ \budget, a search space $\mathcal S$ (lower and upper bounds when $x \in \mathbb R^d$)
\STATE \textbf{Create an initial \doe} of \initDoEsize points in $\mathcal S$: $x^1,\ldots,x^{\initDoEsize}$. 
\STATE Calculate $f(x^1),\ldots,f(x^{\initDoEsize})$, $t \leftarrow $\initDoEsize
\WHILE {$t \le$ \budget or complementary stopping criterion}
\STATE Build a conditional GP from the current \doe : $Y^t$
\STATE Maximize the acquisition criterion to define the next iterate: 
\begin{equation}
x^{t+1} = \arg \max_{x \in \mathcal S} a(x;Y^t)
\label{eq:maxAcquisition}
\end{equation}
\STATE Calculate $f(x^{t+1})$
\STATE $t \leftarrow t+1$
\ENDWHILE
\RETURN $x^\text{best} = \arg \min_{i=1,\ldots,t} f(x^i)$, $f(x^\text{best})$, the final DoE, the final $Y^t$
\end{algorithmic}
\end{algorithm}

\section{Studying BO's degrees of freedom}

BO offers many degrees of freedom that can significantly affect its performance. 
Since new acquisition functions are regularly proposed, it is often seen as the main degree of freedom. 
However, in practical situations one also has to define the GP model (through its mean and kernel functions),
the proportion of the budget allocated to the initial design and the way this design is chosen,
or the way the acquisition function is optimized.

\subsection{Common practices}\label{sec:commonpractices}

To assess common practices, we studied (without seeking exhaustivity) the default settings of an ensemble of recent BO libraries: 
\texttt{DiceOptim} \cite{roustant2012dicekriging}, 
\texttt{STK} \cite{bect2019stk},
MATLAB$^\copyright$'s Bayesian optimization algorithm \cite{matlabBO}  
\texttt{GPyOpt} \cite{gonzalez2016gpyopt},
\texttt{BayesOpt} \cite{martinez2014bayesopt}, 
\texttt{BayesianOptimization} \cite{BayesianOptimization} and 
\texttt{BoTorch} \cite{balandat2019botorch}, 
\texttt{GPflowOpt} \cite{knudde2017gpflowopt},
as well as practical recommendations from recent reviews, tutorials and textbooks \cite{forrester2008engineering,shahriari2016taking,frazier2018tutorial,gramacy2020surrogates}.

\paragraph{Initial DoE size.} All recent practitioners agree on using a small initial set of observations (between 4 and 10), somehow in disruption with earlier practice (for instance, \cite{forrester2008engineering} recommends using 1/3 of the total budget for the initial DoE). Surprisingly, most toolboxes use a fixed number, regardless of the dimension, which may hinder the likelihood optimization process if the number of inferred parameters is greater than the number of observations.

\paragraph{Kernel.} While most toolboxes support various kernels, there is a general consensus on using the Mat{\'e}rn kernel with smoothing parameter $\nu=5/2$ (this is the default choice for all the toolboxes cited above).

\paragraph{Trend function.} All toolboxes either rescale the data for a zero mean and ignore this term, or use a constant trend as default.

\paragraph{EI optimization.} There is also a general consensus on using derivative-based optimization after a ``warm-up'' phase, sometimes paired with a multi-start scheme, although some toolboxes also offer the use of global optimizers (\texttt{GPyOpt}, \texttt{DiceOptim}). More involved procedures 
(e.g. sequential Monte-Carlo, see \cite{feliot2017bayesian}) are generally not available by default.

\paragraph{Modifying the exploration / exploitation trade-off.} Although not generalized, 
some authors recommend introducing some randomness in the acquisitions, 
for instance using an $\varepsilon$-greedy strategy (\texttt{BayesOpt}), to favor exploration and increase robustness. 
Conversely, it is possible to favor exploitation by optimising the GP mean rather than the EI (as will be done later).

\paragraph{Non-linear transformations.} Besides linear rescaling, transformations of either the input or output spaces have been proposed, in particular to handle non-stationarity (see e.g. \cite{snelson2004warped,snoek2014input}). However, none of the toolboxes cited above propose this as a default.

\subsection{Factors studied in this work}
\label{sec-factors}

The factors varied in this experiment lie in two categories. The first category corresponds to the parameters every user must choose and for which guidance is not obvious: 
\begin{itemize}
 \item the size of the initial design of experiment
 \item the GP kernel, and 
 \item the trend function.
\end{itemize}
 The second category corresponds to possible refinements of BO, in order to perform an ablation study: 
 \begin{itemize}
  \item using a warping function over the outputs \cite{snelson2004warped},
  \item using a warping function over the inputs \cite{snoek2014input},
  \item modifying the exploration / exploitation trade-off, and
  \item changing the way EI is optimized.
 \end{itemize}

The factors are likely to have a joint effect on BO performance, so we need to study them jointly rather than independently. However, as each factor has 2 or more possible values, a full factorial design is out of reach, due to its computational cost. We choose instead to evaluate only some combinations of factors, based on likely interactions. In addition, based on our first findings, we added a few configurations that correspond to possible optimal choices for all factors. The full set of configurations is reported in Table \ref{tab:allruns}.

\begin{table}
\begin{tabular}{|c|p{0.15\textwidth}|p{0.7\textwidth}|}
\hline
ID & name & comments \\
\hline
1 & M & default setting\\
2 & S & small initial DoE\\
3 & L & large initial DoE\\
4 & LinM & linear trend\\
5 & QuadM & quadratic trend\\
6 & ScalM & scaling (input warping)\\
7 & ScalS & scaling (input warping), small initial DoE size\\
8 & ScalL & scaling (input warping), large initial DoE size\\
9 & WarpM & output warping\\
10 & WarpS & output warping, small initial DoE size\\
11 & WarpL & output warping, large initial DoE size\\
12 & ExpM & exponential kernel\\
13 & ExpS & exponential kernel, small initial DoE\\
14 & ExpScalM & exponential kernel and input warping\\ 
15 & MeanM & mean GP as a proxy criterion\\
16 & EirandM & global search only for EI\\ 
17 & EilocM  & local search only for EI\\
18 & MeanS & mean GP as an occasional acquisition function, small initial DoE\\
19 & ExpMeanS  & exponential kernel, proxy criterion, small initial DoE\\
20 & QuadMean & proxy criterion, quadratic trend, initial DoE of size $2d+1$ (to allow estimating the extra model parameters)\\
\hline 
\end{tabular}
\caption{Different BO configurations. If not specified in the name, the parameters are the default values: medium DoE size, Matérn 5/2 kernel, constant trend and multi-start BFGS for EI optimization.
\label{tab:allruns}
}
\end{table}

\paragraph{Initial DoE size.}
All the DoEs are generated by hypercube sampling, optimized with respect to the \textit{maximin} criterion \cite{fang2005design}, which is the standard choice in BO.
Three sizes of the initial DoE are investigated; all depend linearly on the dimension. 
The small (S) budget is $d+4$ (number of GP hyperparameters + 2), the medium (M) is $7.5 d$ and the large (L) is $20 d$. 
Assuming a total budget of $30d$, M and L correspond respectively to 25\% and 66\% of the total budget.
The medium budget is the default of the algorithm.

\paragraph{Kernel.} Classical choices for the GP kernel include the exponential, squared exponential and Mat{\'e}rn families \cite{shahriari2016taking}. Our default choice, aligned with common practice, is the anisotropic\footnote{also referred to as ARD, for Automatic Relevance Determination} Mat{\'e}rn kernel with smoothing parameter $\nu=5/2$, which provides twice differentiable GP models. For a contrasted alternative, we chose the exponential kernel, which provides continuous but non-differentiable models (for both the GP mean and the trajectories).

\paragraph{Trend function.} We stick here to the \textit{universal kriging} framework \cite{matheron1963principles}, and consider only polynomials as basis functions for the trend. We consider 1) constant (default), 2) linear without interaction, 
3) quadratic without interaction trends\footnote{e.g., for two parameters $x_1$ and $x_2$, 
$\beta_0 + \beta_1 x_1 + \beta_2 x_2 + \beta_{11} x_1^2 + \beta_{22} x_2^2$. 
The interaction terms are omitted to limit the number of parameters of the model.}.

\paragraph{Output warping.} It consists in applying a transformation to the observations. This step is expected to be beneficial when the data is strongly non-Gaussian. We follow here the approach of \cite{snelson2004warped}, which suggests as a transformation a sum of tanh functions, $f^i_\text{warp} = f^i + \sum_{j=1}^J a_j tanh \left(b_j( c_j + f^i) \right)$, with $a_j$, $b_j$ and $c_j$ a set of warping parameters that can be inferred by maximum likelihood. We use here $J=2$. Note that as this implies a substantial computational overhead, we only apply it once at the initial stage, prior to the BO loop.
As this method may benefit from more data, we study its influence jointly with the DoE size.

\paragraph{Input warping.} Similarly, one can apply a transformation over the input space. It is possible to apply independently transformations for each input \cite{sampson1992nonparametric,snoek2014input}. 
We follow here the approach of \cite{xiong2007non}, available in \texttt{DiceKriging}, that allows non-linear deformations to be applied along each canonical axe. In a nutshell, it sets a small number of ``knots'' in the input space, which are used to define piecewise linear deformations. This approach has a relatively moderate number of parameters that can be inferred by maximum likelihood. To avoid confusion with output warping, we refer to this method as \textit{(input) scaling} in the rest of the paper.

\paragraph{GP mean acquisition.} Recent works found that EGO sometimes struggles to find the optimum with precision \cite{Trike_2016,mcleod2018optimization}, either because EI intrinsically favours exploration or because its shape becomes quickly very flat with a couple of narrow peaks, which makes its maximization overly challenging. 
It has also been suggested in \cite{de2021greed,rehbach2020expected} that BO lacks exploitation for functions that are well approximated by the GP mean and in dimensions higher than 5. A strategy mixing GP mean minimization and a random search in the mean-variance Pareto set was found effective in \cite{de2021greed}.
In the same spirit, we propose here to enhance exploitation by replacing for some steps EI by minus the GP mean, in order to favour exploitation and obtain an easier-to-evaluate acquisition function. With this option active, this alternative acquisition is used whenever the EI maximization fails, but also at a regular frequency (once every 5 iterations).

\paragraph{EI optimization.} An important question for practitioners is how much effort they should dedicate to the EI optimization. One might argue that only a local optimum is necessary, as the global optimum might be found at the next BO loop, and a local optimum at one iteration can become the global one during the next. Conversely, one might consider that fine convergence is unnecessary, considering the uncertainty associated with the GP hyperparameters. We compare here a commonly used base strategy (multi-start BFGS) with a purely local search (single BFGS run from a random initial point) and a purely global (random) search (see next section for details).

\section{Implementation}
\label{sec:implementation}
All the experiments presented in the work are based on \texttt{DiceOptim}, an \texttt{R} package dedicated to BO, 
see \cite{roustant2012dicekriging,picheny2014noisy} for a detailed description.
As \texttt{DiceOptim}'s target use is expensive-to-evaluate computer experiments, the previous subroutines for training the model and optimizing the acquisition function were designed for accuracy, not speed, which precludes a direct use of past versions of \texttt{DiceOptim} on \COCO. In addition, 
whenever a BO iteration is unsuccessful, the optimization stopped immediately (to allow the user to fix the problem manually and avoid computing useless expensive experiments). This would lead to many early stopped runs with \COCO.

In this section, we detail our new implementation (the v2.1.1, \cite{diceoptimV211}) and the changes made for the two main BO steps (model update and acquisition) in favour of speed and reliability, while limiting the impact on accuracy.

\subsection{GP training and updating}
Every time a new point is added to the training set, the model hyperparameters are adjusted by re-optimizing the model likelihood, which constitutes the first bottleneck of BO. This step is particularly challenging, as the likelihood is a multimodal function, hence requiring a global optimizer. In addition, some hyperparameter values can lead to covariance matrix ill-conditioning and numerical failures. 

\paragraph{Data rescaling.} In accordance with common practices, we center and normalize the observations prior to model fitting. This rescaling is fixed and based on the initial DoE mean and standard deviation.

\paragraph{Likelihood optimizer.} We use a multi-start L-BFGS strategy. The starting points are chosen using a latin-hypercube sample of size $2d$.
This allows us to leverage the closed-form expressions of the likelihood derivatives, while limiting the risk of getting trapped in a local optimum region.

\paragraph{Clipping lengthscales.} Bounding the search space of the kernel lengthscale parameters reduces the risk of failure. Indeed, too large lengthscales lead to 
covariance matrix ill-conditioning, while very small lengthscales often correspond to a local optimum of the likelihood, but yield a non-informative model (white noise model). Here, we choose as bounds $(U - L) \sqrt d /100$ and $(U - L)\sqrt d$, which offers a good trade-off between risk limitation and model expressivity. The $\sqrt d$ term allows lengthscales to have larger values with more dimensions.
This rule is empirical, but motivated by practical insight. The diagonal of the unit cube (largest possible distance between two points of the input space) is equal to $\sqrt{d}$. Hence, assuming all dimensions are treated similarly and the kernel depends on a normalized distance (e.g., the squared exponential, the Matérn kernels etc), our upper bound for the lengthscale (proportional to $\sqrt{d}$) delivers a constant maximal correlation between the two farthest points in the space, regardless of the dimension.

\paragraph{Regularization.} A well-known method to limit covariance matrix ill-conditioning is to add a small positive constant to its diagonal (often referred to as \textit{nugget effect} or \textit{jitter}). However, this diagonal term makes the GPs not interpolating, which can affect in particular the late exploitation stages. Here, as default we do not use any regularization; if the hyperparameter optimization step fails during one iteration (due to a Cholesky decomposition error for instance), we introduce a regularization term to the covariance matrix, equal to a small fraction  ($10^{-12}$) of the previous model's variance and restart training. The regularization is then used for all further iterations.

\subsection{Acquisition function optimization}
\label{sec:acqfuncoptim}

EI is usually a highly multimodal function with large plateaus with 0 value, and maximizing it is a challenging task that is the second computational bottleneck of BO. 

\paragraph{EI optimizer.} We use the following heuristic, in line with current practice: EI is first evaluated over a large space-filling design of size $\min(2000, 500d)$; then the L-BFGS algorithm is used for local convergence, starting from the best point. This procedure is repeated $\min(10, d)$ times, leveraging parallel processing. When studying the effect of the EI optimizer, the first alternative is to use a single L-BFGS run starting from a single random point. The second alternative is to remove the L-BFGS step, the EI being effectively optimized by random search.

\paragraph{Acquisition failures.} Sometimes, one of the above optimizations gets trapped in one of the large plateaus and returns an acquisition point with $EI=0$ (up to numerical noise). In that case, we either re-run the optimization with minus the GP mean as an acquisition function (if the ``Mean'' option is activated, see \textit{GP mean acquisition}), or pick the acquisition point at random over the search space (hence, introducing this random step amounts to using a form of epsilon-greedy strategy).


\paragraph{Avoid conditioning issues.} EI maximization sometimes result in an acquisition point very close to an existing observation. Adding this point results in covariance matrix ill-conditioning and model failure. Following \cite{binois2019gpareto}, we avoid this issue by checking prior to observation if the proposed point satisfies a proximity criterion\footnote{$\min_{i=1 \ldots t} \left( \left[c^t(x^i,x^i) + c^t(x^*,x^*) - 2 c^t(x^i,x^*)\right]  / k(x^*,x^*) \right) \geq 10^{-6}$}, and replace it by a random sample over the search space if it is found critically close to already sampled points.

\paragraph{Enhancing exploration by decreasing range values.} Repeated EI maximization failures are often due to deceiving functions, that is, for which a small training set indicates a very ``flat'' function, while it also has narrow peaks \cite{forrester2008global}. 
While there exist fully Bayesian approaches to handle this problem \cite{benassi2011robust}, their computational overhead prevents us from using it with \COCO.
Here we use the following very simple trick: if several (3) consecutive failures are recorded, we do not choose the range parameters by maximum likelihood but we set them as a fraction (2/3) of the previous values. This generally increases the prediction uncertainty, hence exploration, which allows recovering from deceiving samples.

\section{Description of the \COCO benchmark}\label{sec:cocodescription}
The comparison of the variants of the optimizers is based on the \COCO \cite{hansen2016coco} software.
The set of objective functions which is used in the comparisons is the BBOB noiseless suite \cite{hansen2010comparing,bbob_online_functions}.
It is made of 24 functions divided into 5 groups: the separable functions (which are written as a sum of univariate functions), the functions with low or moderate conditioning, the unimodal functions with high conditioning, the multimodal functions with adequate global structure and the multimodal functions with weak global structures. 
Table \ref{tab-bbob_functions} gives the function numbers, their names, the \emph{group} they belong to, and additional comments. The conditioning is a quantification of how stretched a function is and is defined as the maximum difference in rate of change of a function at the optimum between two directions. 
The function groups in \COCO are useful to understand the adequation between a specific optimizer and features of the problem.
\begin{table}
\begin{tabular}{|c|p{0.3\textwidth}|p{0.65\textwidth}|}
\hline
ID & name & comments \\
\cline{1-3}
\multicolumn{3}{c}{separable functions}\\
\cline{1-3}
f1 & Sphere & unimodal, allows to checks numerical accuracy at convengence \\
f2 & Ellipsoidal & unimodal, conditioning $\approx 10^6$ \\
f3 & Rastrigin & $10^d$ local minima, spherical global structure \\
f4 & B\"uche-Rastrigin & $10^d$ local minima, asymmetric global structure \\
f5 & Linear Slope & linear, solution on the domain boundary \\
\cline{1-3}
\multicolumn{3}{c}{functions with low or moderate conditioning}\\
\cline{1-3}
f6 & Attractive Sector & unimodal, highly asymmetric \\
f7 & Step Ellipsoidal & unimodal, conditioning $\approx 100$, made of many plateaus \\
f8 & Original Rosenbrock & good points form a curved $d-1$ dimensional valley \\
f9 & Rotated Rosenbrock & rotated f8 \\
\cline{1-3}
\multicolumn{3}{c}{unimodal functions with high conditioning $\approx 10^6$ }\\
\cline{1-3}
f10 & Ellipsoidal & rotated f2 \\
f11 & Discus & a direction is 1000 times more sensitive than the others \\
f12 & Bent Cigar & non-quadratic optimal valley \\
f13 & Sharp Ridge & resembles f12 with a non-differentiable bottom of valley \\
f14 & Different Powers & different sensitivities w.r.t. the $x_i$'s near the optimum \\
\cline{1-3}
\multicolumn{3}{c}{multimodal functions with adequate global structure}\\
\cline{1-3}
f15 & Rastrigin & rotated and asymmetric f3 \\
f16 & Weierstrass & highly rugged and moderately repetitive landscape, non unique optimum \\
f17 & Schaffers F7 & highly multimodal with spatial variation of frequency and amplitude, smoother and more repetitive than f16 \\
f18 & moderately ill-conditioned Schaffers F7 & f17 with conditioning $\approx 1000$ \\
f19 & Composite Griewank-Rosenbrock & highly multimodal version of Rosenbrock \\
\cline{1-3}
\multicolumn{3}{c}{multimodal functions with weak global structure}\\
\cline{1-3}
f20  & Schwefel & $2^d$ most prominent optima close to the corners of a shrinked and rotated rectangle \\
f21 & Gallagher's Gaussian 101-me peaks & 101 optima with random positions and heights, conditioning $\approx 30$ \\
f22 & Gallagher's Gaussian 21-hi peaks & 21 optima with random positions and heights, conditioning $\approx 1000$ \\
f23 & Katsuura & highly rugged and repetitive function with more than $10^d$ global optima \\
f24 & Lunacek bi-Rastrigin & highly multimodal function with 2 funnels, one leading to a local optimum and covering about 70\% of the search space \\
\hline
\end{tabular}
\caption{
Functions of the BBOB noiseless testbed with comments. The functions are divided in groups.
\label{tab-bbob_functions}
}
\end{table}

In \COCO, an optimization algorithm is tested by repeating runs over different instances of each function, where an \emph{instance} is the function with a random rotation of the coordinate system and a random translation of the optimum. 
\emph{A problem} is a pair \{function, minimum value to reach\}. Therefore, for each instance of a function, there are several problems to solve of difficulty varying with the target value.
To set the target values and more generally define a reference performance, \COCO relies on a composite fake algorithm called best09. best09 is made at each optimization iteration of the best performing algorithm of the Black-Box Optimization Benchmarking (BBOB) 2009 \cite{hansen2010comparing}. 
Relevant target values can then be estimated as values reached by best09 at a given number of function evaluations: low budgets (e.g., $2 d$) give easy targets, larger budgets more difficult targets. In our experiments, the targets were set at the values reached by best09 after $[0.5,1,3,5,7,10,15,20]\times d$ function evaluations.
In compliance with the expensive function assumption, all optimization runs are $30 \times d$ function evaluations long, including the initial DoE.
The stopping criterion of the optimization which is a maximal number of calls to $f$ should not be mistaken with the above ``minimum value(s) to reach'' which are only used in the performance metric calculated after the search.

\COCO results are based on several metrics. In this study, we focus on the \emph{Empirical Run Time Distributions} (ERTD) that gives, at a given number of evaluations, the proportion of problems which are solved by an algorithm. ERTD is displayed as a function of the logarithm of the number of evaluations, rescaled by the problem dimension\footnote{Note that due to this logarithm, some ERTD curves extend beyond the lower limit of the x-axis.}. ERTD curves increasing quickly to one indicate good performance. ERTD curves reaching or exceeding the best09 ones indicate excellent performances. In the following, we also report the ERTD of random search, which can be used as a lower bound on performance.


\section{Results}
In this section, we present selected results from the full \COCO experiments to analyze the effects of the EGO factors and compare EGO to the state-of-the-art algorithms. The full results are available as supplementary material (see \cite{supplement_BO_COCO}).


\subsection{Choosing EGO fundamental parameters}

\subsubsection{Global analysis}

\begin{figure}
\begin{tabular}{cc@{}c@{}c@{}}
& \textbf{3D} & \textbf{5D} & \textbf{10D} \\
\rotatebox[origin=l]{90}{\textbf{ DoE size}} & 
\includegraphics[width=0.3\textwidth]{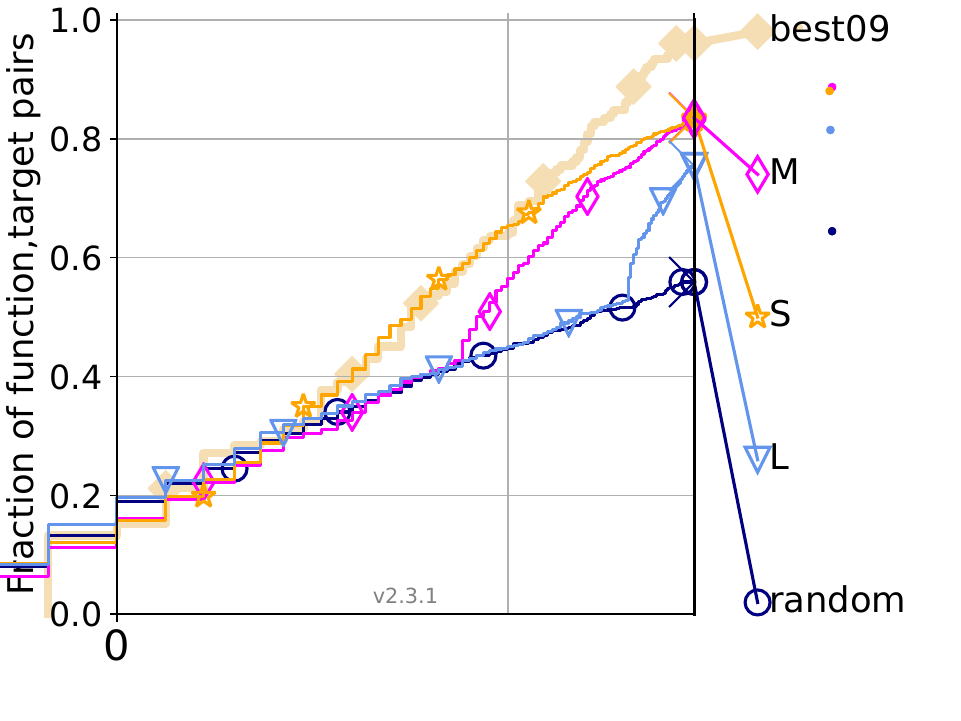}&
\includegraphics[width=0.3\textwidth]{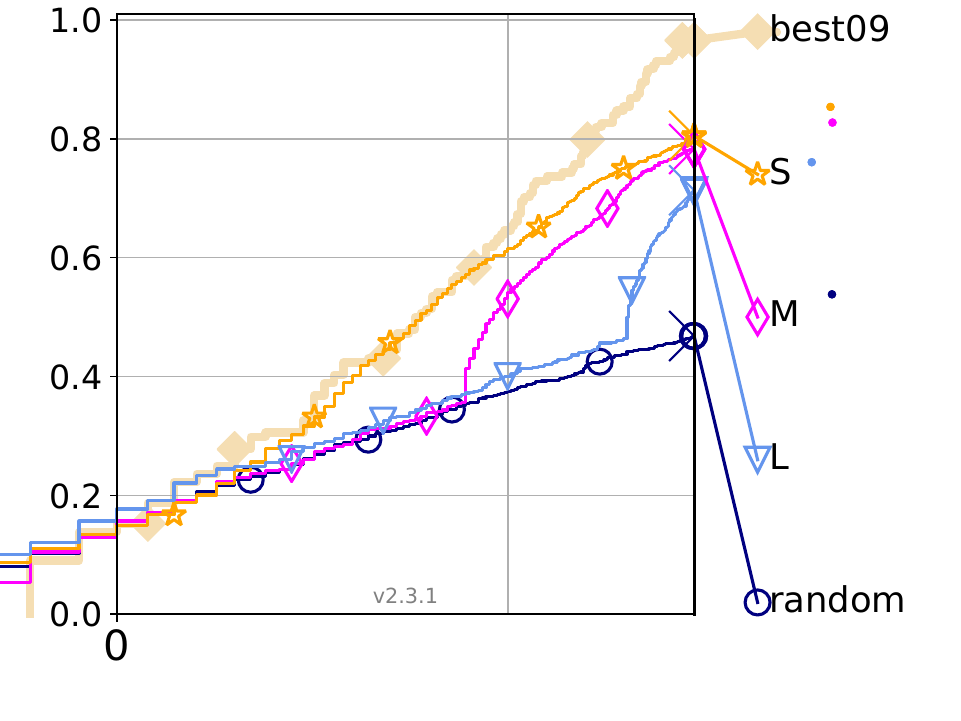}&
\includegraphics[width=0.3\textwidth]{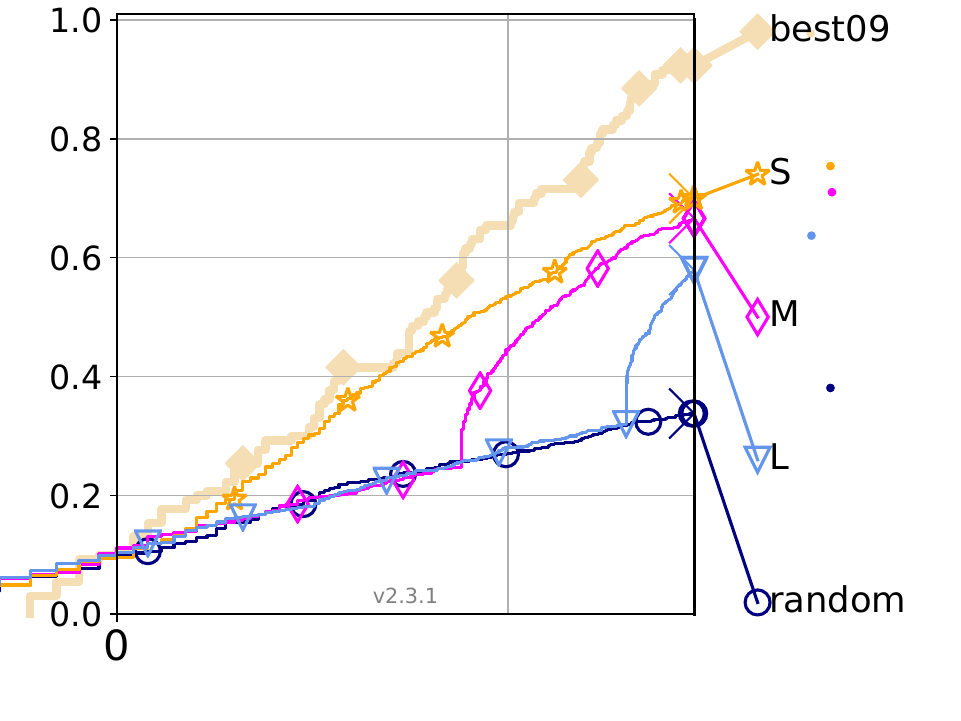}\\
\rotatebox[origin=l]{90}{\textbf{ Exp kernel}} & 
\includegraphics[width=0.3\textwidth]{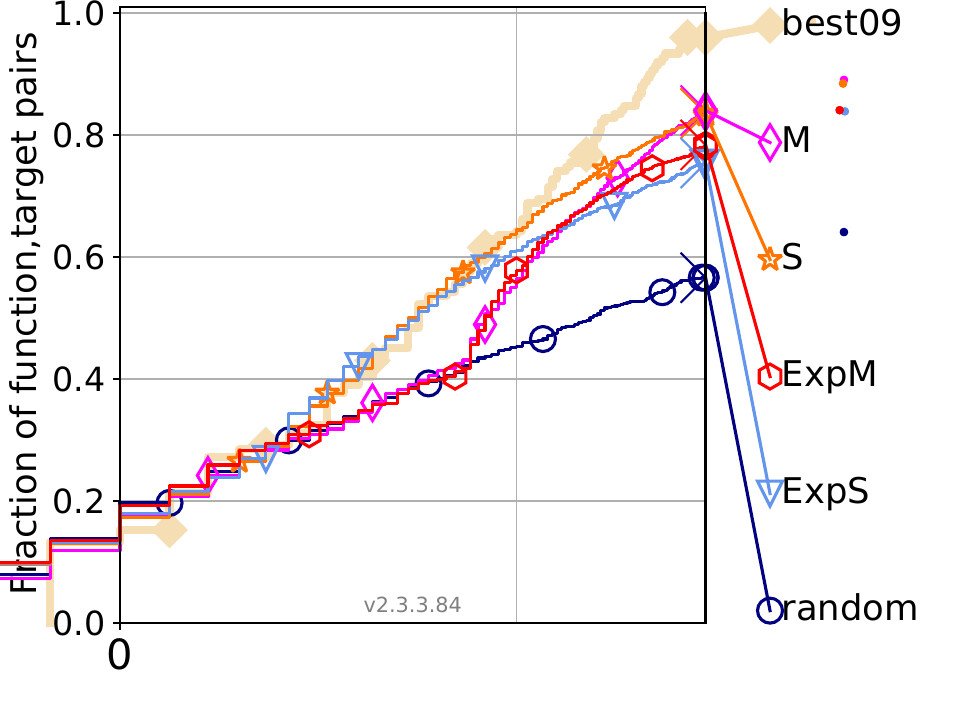}&
\includegraphics[width=0.3\textwidth]{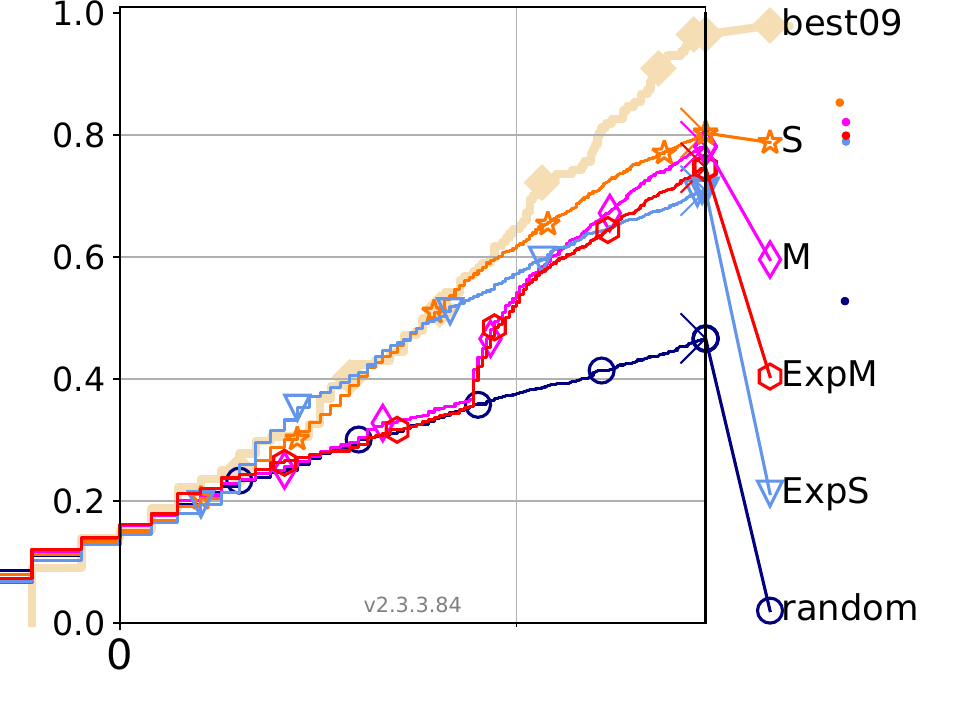}&
\includegraphics[width=0.3\textwidth]{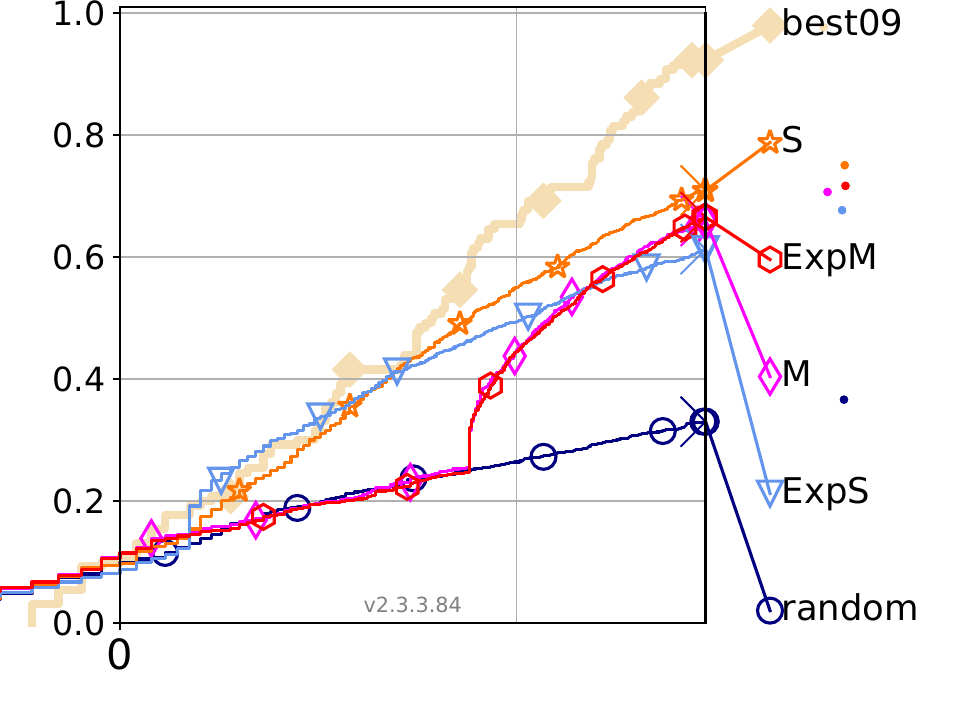}\\
\rotatebox[origin=l]{90}{\textbf{ Trend}} & 
\includegraphics[width=0.3\textwidth]{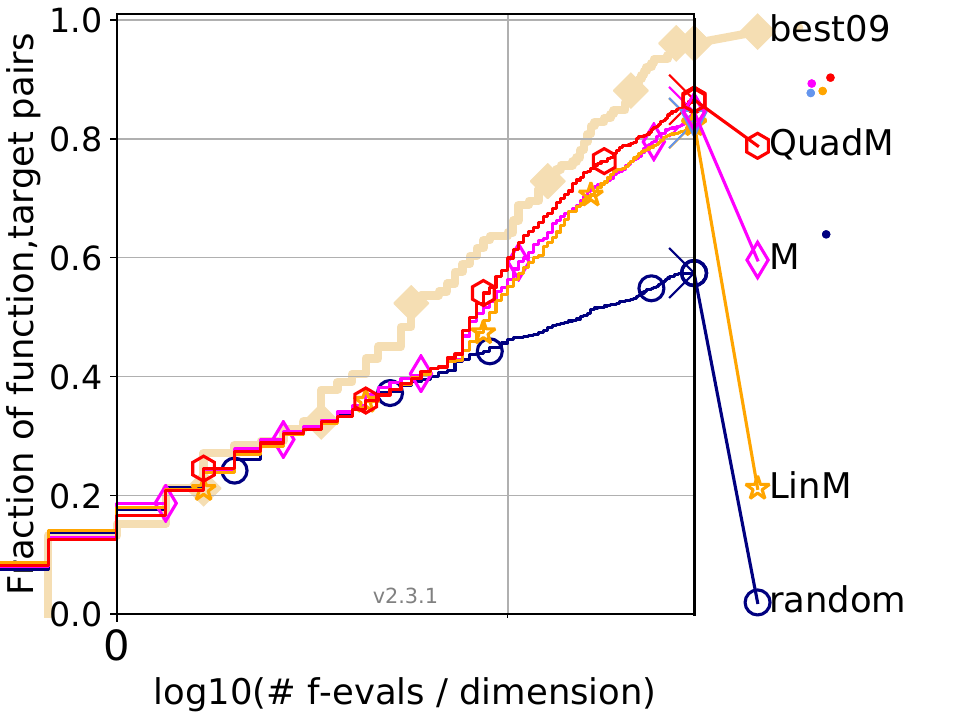}&
\includegraphics[width=0.3\textwidth]{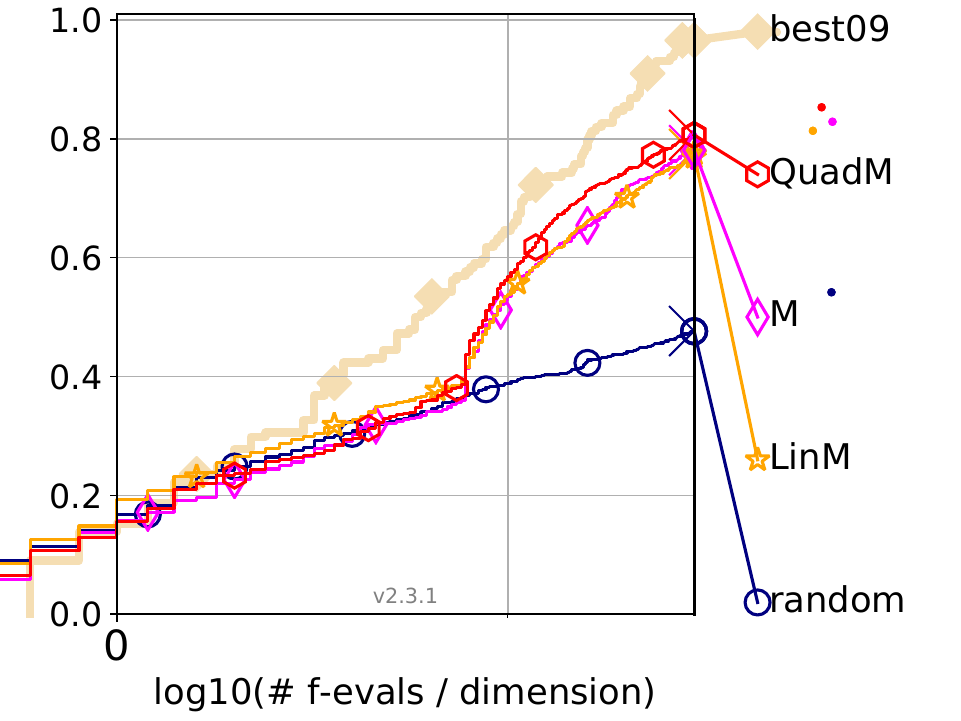}&
\includegraphics[width=0.3\textwidth]{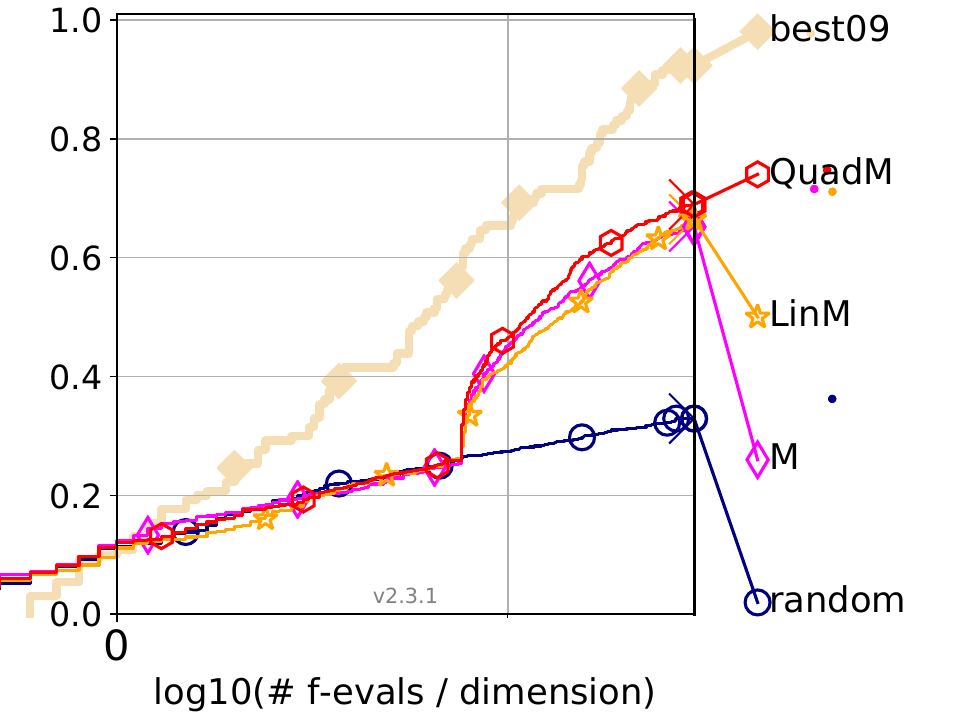}
\end{tabular}
\caption{
ERTD over all functions in $d=3,5$ and $10$ dimensions. 
The top row corresponds to varying initial DoE sizes : small (S, size = $d+4$), medium (M, $7.5 d$) and large (L, $20 d$). 
The middle row corresponds to the use of the exponential kernel (jointly with the initial DoE size): algorithms with ``Exp'' in the name use an anisotropic exponential kernel, the others use an anisotropic Mat\'ern 5/2. 
The bottom row corresponds to different trend functions: M (constant), LinM (linear), QuadM (quadratic). 
See Table~\ref{tab:allruns} for further details.
\label{fig-ERTD_all}
}
\end{figure}

Figure~\ref{fig-ERTD_all} shows the ERTD over all functions in $d=3,5$ and $10$ dimensions, for different configurations of the initial DoE size, kernel or trend.

\paragraph{Initial DoE size.}
The effect of the initial budget (top row) is clearly visible. During the initial LHS, the optimizer progresses like a random search (in blue). Once the Bayesian search starts, rapid progress is made, and the convergence speed is initially faster when the initial DoE is larger. 
Contrarily to a common belief that using a large portion of the budget to an initial random exploration leads to better models and eventually better performance\cite{forrester2008engineering}, there is no observed drawback for starting with a small initial DoE (size $d+4$) during the first $30d$ evaluations. Quite the opposite, there is a very substantial gain earlier in the search that creates a gap still visible at the end of the $30d$ budget.
Both S and M runs tend to the same performance at the end, which indicates that even for larger budgets than studied here ($>30d$), 
the initial DoE size may have a moderate impact.

\paragraph{Kernel choice.}
The kernel choice is analyzed in combination with the DoE size,
however the tests do not include large initial DoE sizes which have always been found not to yield competitive algorithms.
Figure~\ref{fig-ERTD_all} (middle row) shows a substantial difference in favour of the Mat\'ern kernel irrespective of the initial budget. 

\paragraph{Trend.}
Overall, M (constant trend) and LinM (linear trend) performances are indistinguishable, while QuadM (quadratic trend) slightly ouperforms them. 
This is quite unexpected, as most current practices (see Section \ref{sec:commonpractices}) do not include a trend at all.

We next detail our analysis by breaking the ERT distributions by functions groups. For simplicity we only present the results for the dimension 5,
which were found representative.

\subsubsection{Initial budget}

\begin{figure}
    \centering
\begin{tabular}{cc@{}c@{}c@{}}
& {Separable} & {Low conditioning} & {High conditioning} \\
\rotatebox[origin=l]{90}{\textbf{DoE size }} & 
\includegraphics[trim=0mm 0mm 10mm 0mm, clip, width=0.3\textwidth]{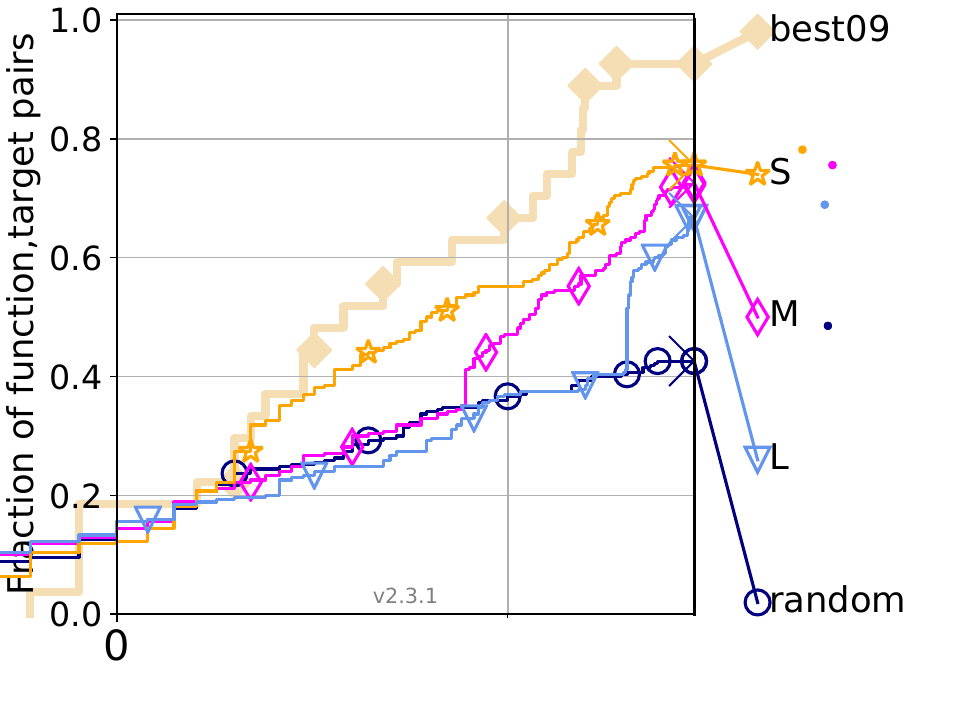} &
\includegraphics[trim=0mm 0mm 10mm 0mm, clip, width=0.3\textwidth]{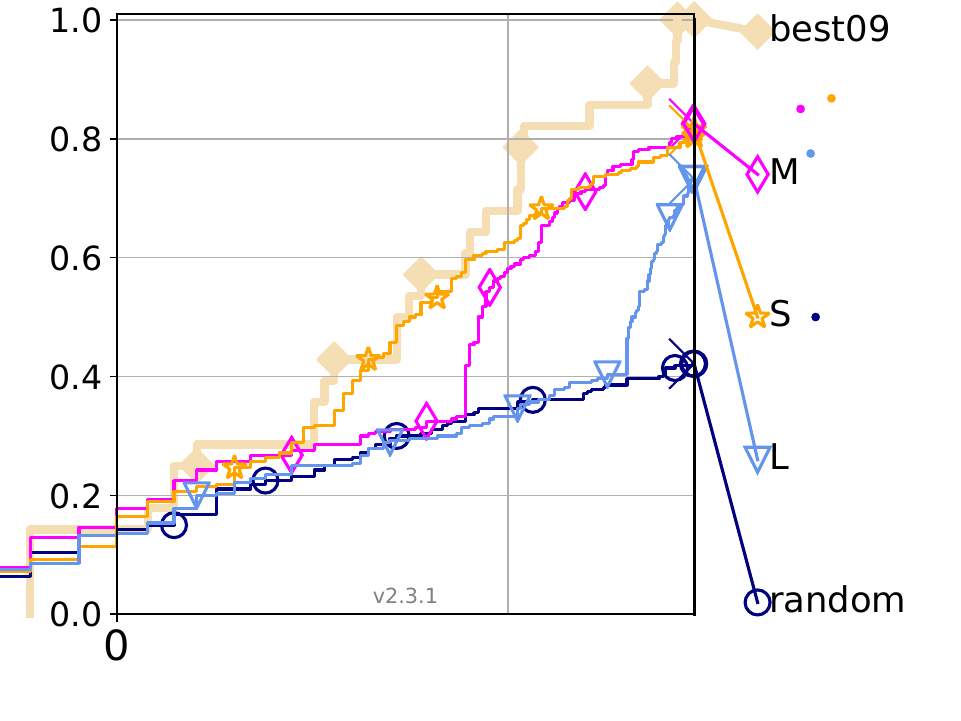}&
\includegraphics[trim=0mm 0mm 10mm 0mm, clip, width=0.3\textwidth]{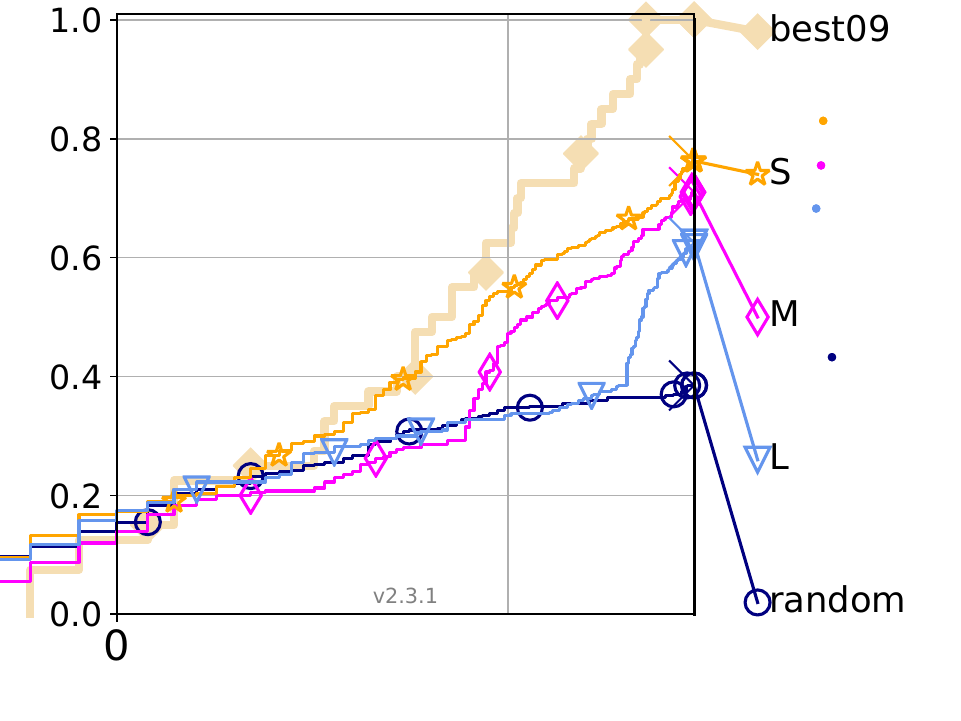} \\
& {Multimod., strong struct.} & {Multimod., weak struct.} & \\
&
\includegraphics[trim=0mm 0mm 10mm 0mm, clip, width=0.3\textwidth]{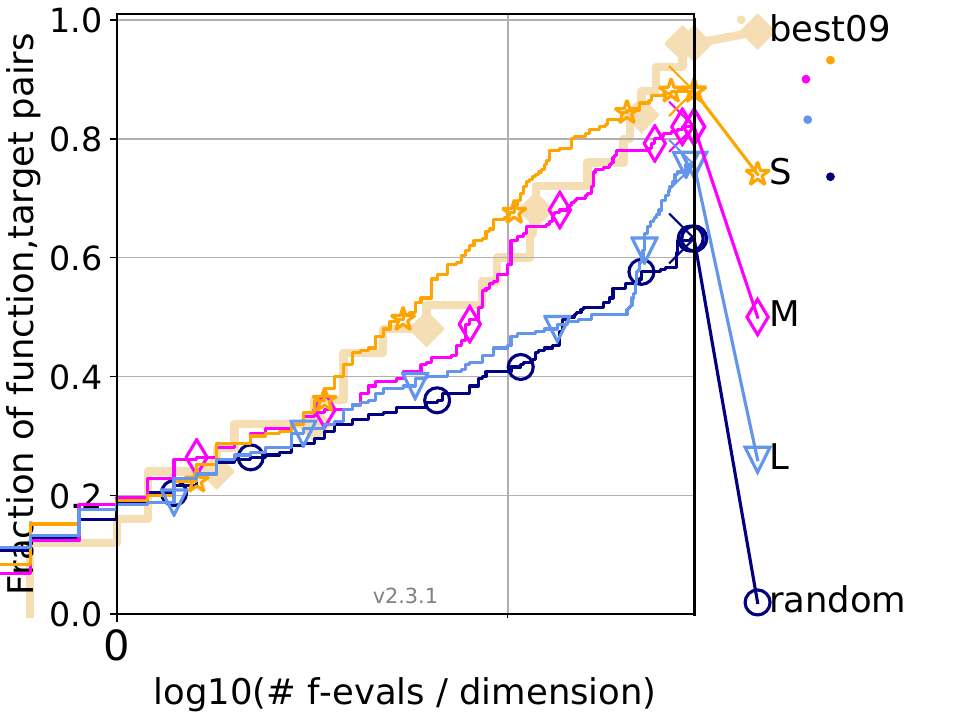}&
\includegraphics[trim=0mm 0mm 10mm 0mm, clip, width=0.3\textwidth]{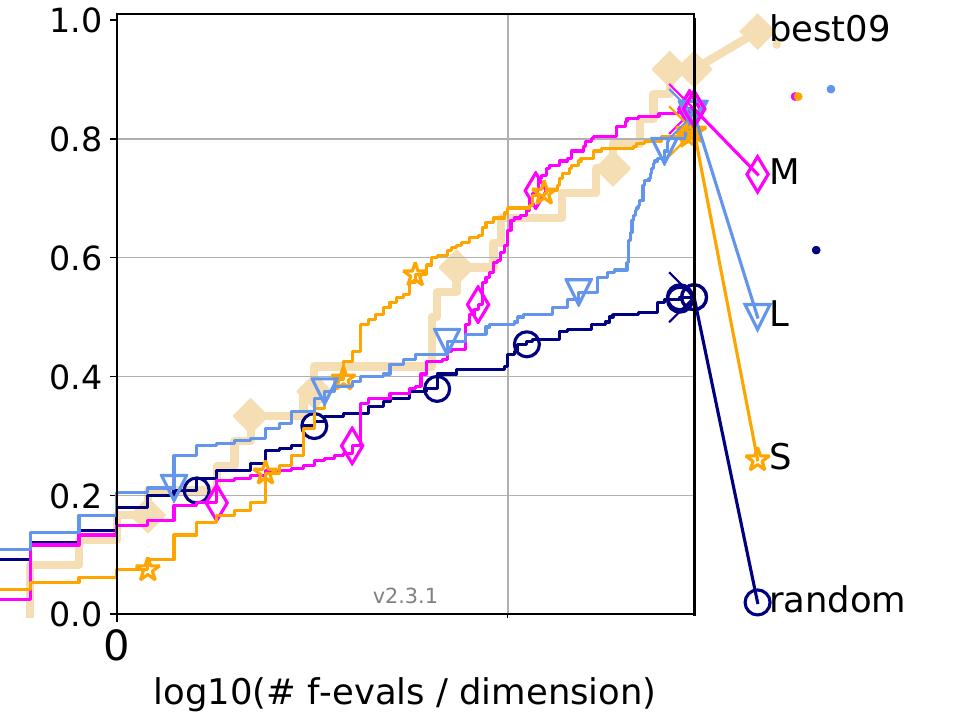}&
\end{tabular}
\caption{ERTD by function groups in $d=5$ dimensions with varying initial DoE sizes : small (S, size = $d+4$), medium (M, $7.5 d$) and large (L, $20 d$).
\label{fig-ERTD_budget_05D_groups}
}
\end{figure}

On Figure~\ref{fig-ERTD_budget_05D_groups},
we observe again the overall outperformance of the smallest initial budget, while the largest budget always has the worst performance. 
Note that for the last group, with highly multimodal functions and for which global exploration is most important, the medium initial budget outperforms the small one towards the last iterations, indicating a slightly better capacity to escape local minima. 

It is finally observed that EGO with small or medium initial DoE performs as well as, and sometimes better than, the utopian best09 reference for multimodal functions, in particular those with adequate global structure which can be learned by the GPs. 

The main individual exceptions, i.e., the functions that benefit from a medium or large size initial DoE, are the attractive sector (function 6) and the Lunacek bi-Rastrigin (function 24), cf. Fig.~\ref{fig-ERTD_budget_exceptions} in Appendix. 
Both functions have funnels that need sufficient global exploration to be located.

\subsubsection{Exponential kernel}

\begin{figure}
    \centering
\begin{tabular}{cc@{}c@{}c@{}}
 &{Separable} & {Low conditioning} & {High conditioning} \\
 \rotatebox[origin=l]{90}{\textbf{Exp kernel}} & 
\includegraphics[trim=0mm 0mm 10mm 0mm, clip, width=0.3\textwidth]{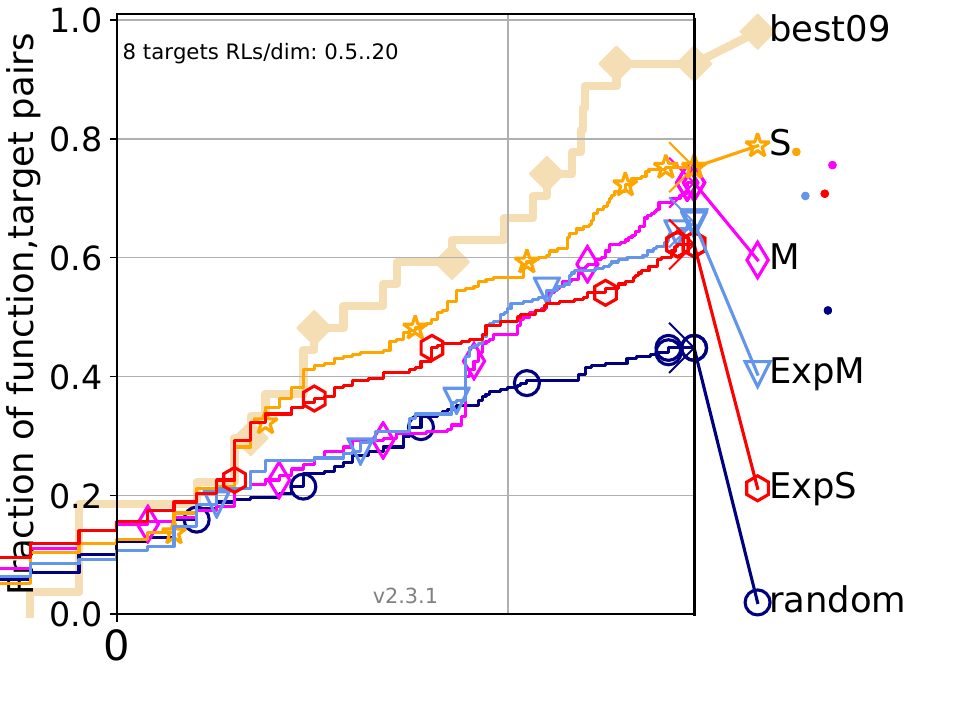} &
\includegraphics[trim=0mm 0mm 10mm 0mm, clip, width=0.3\textwidth]{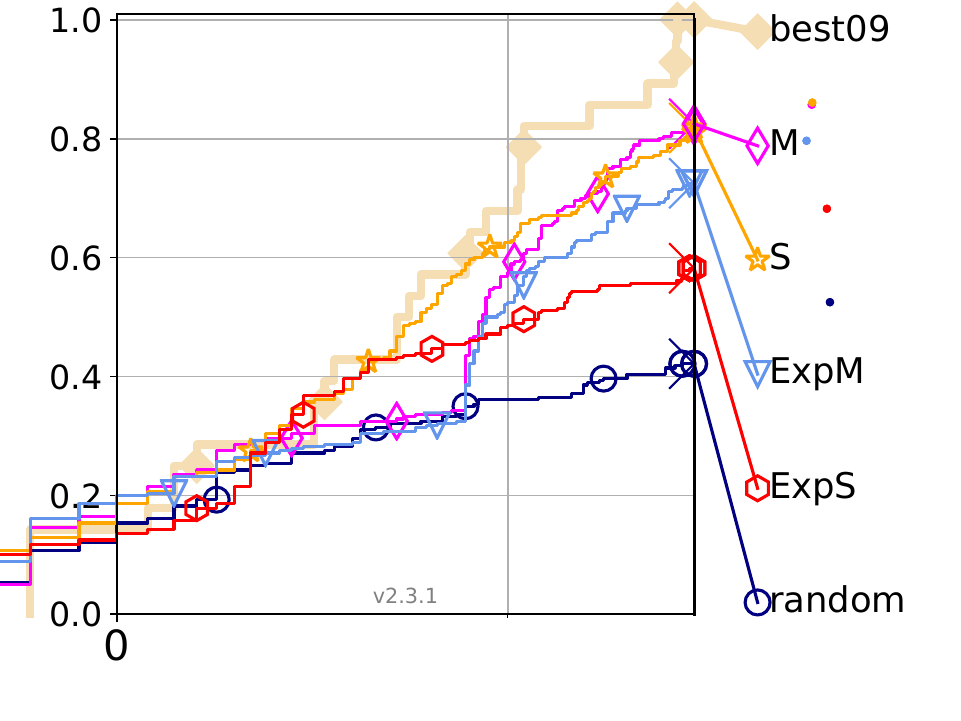}&
\includegraphics[trim=0mm 0mm 10mm 0mm, clip, width=0.3\textwidth]{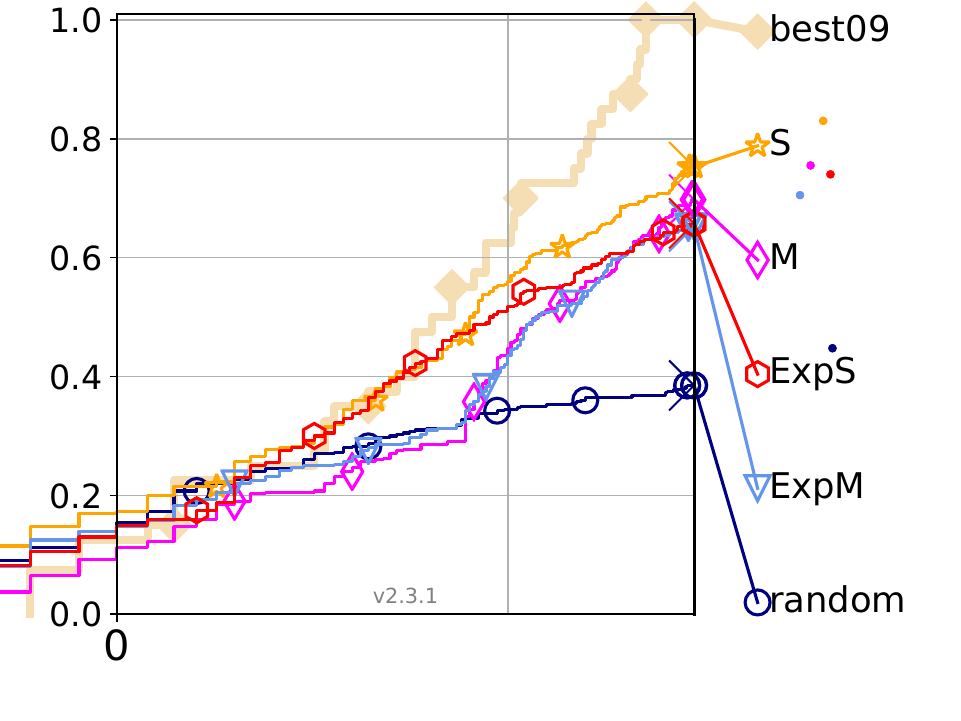} \\
& {Multimod., strong struct.} & {Multimod., weak struct.} & \\
&
\includegraphics[trim=0mm 0mm 10mm 0mm, clip, width=0.3\textwidth]{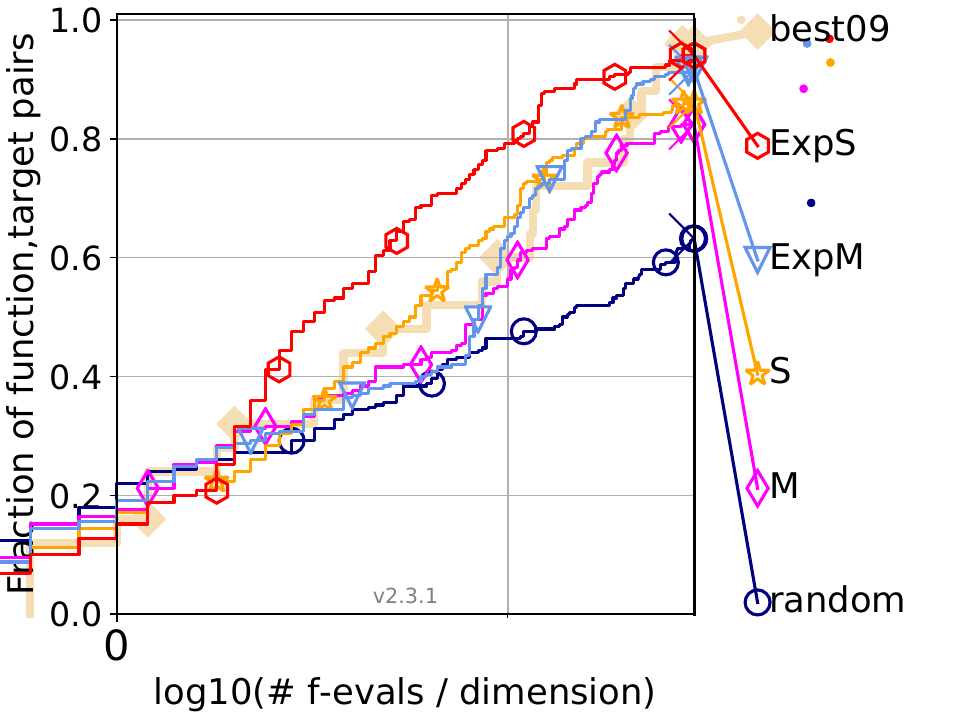}&
\includegraphics[trim=0mm 0mm 10mm 0mm, clip, width=0.3\textwidth]{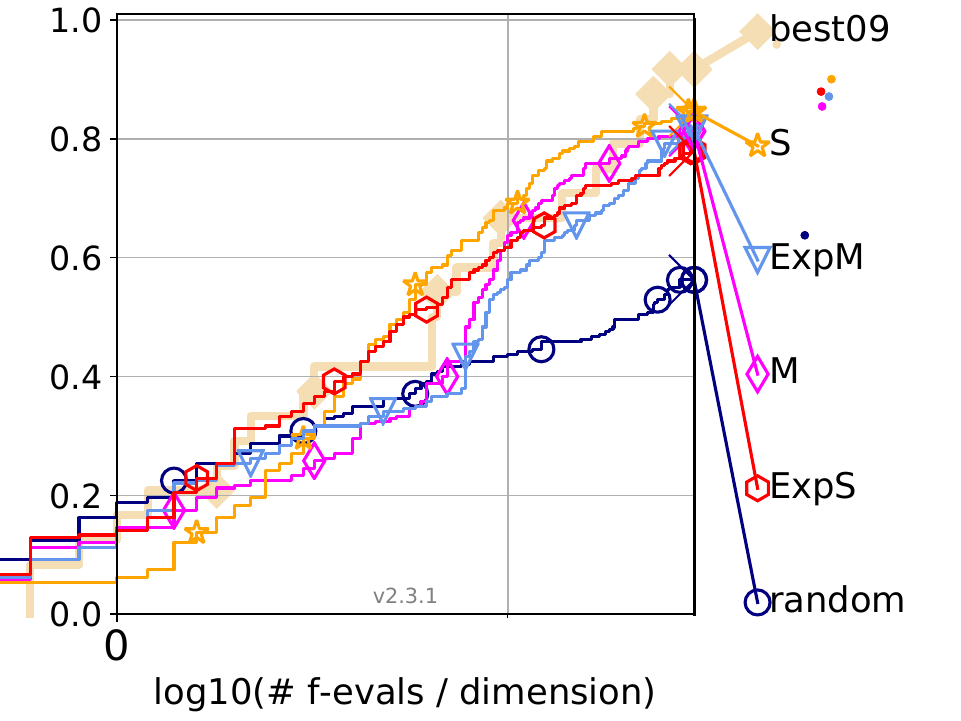}&
\end{tabular}
\caption{
ERTD by function groups in $d=5$ dimensions with exponential (Exp) and Mat{\'e}rn 5/2 kernels and small (S) or medium (M) DoE's.
\label{fig-ERTD_kernel_05D_groups}
}
\end{figure}

Figure~\ref{fig-ERTD_kernel_05D_groups} shows some surprisingly contrasted results. 
While on average \texttt{Exp} is dominated by Mat{\'e}rn, this largely depends on the function group.
\texttt{Exp} shows poor performance for the smooth functions (groups 1 and 2), amplified by a negative interaction with the small initial DoE. 
It is also dominated on the groups 3 and 5 but to a lesser extent.
However, the exponential kernel is beneficial when tackling multimodal functions with adequate global structure (4th group) irrespectively of the DoE size (ExpS and ExpM). On that function group and $d=5$, ExpS is significantly better than best09.

\begin{figure}
\centering
\hfill
\includegraphics[width=0.33\textwidth]{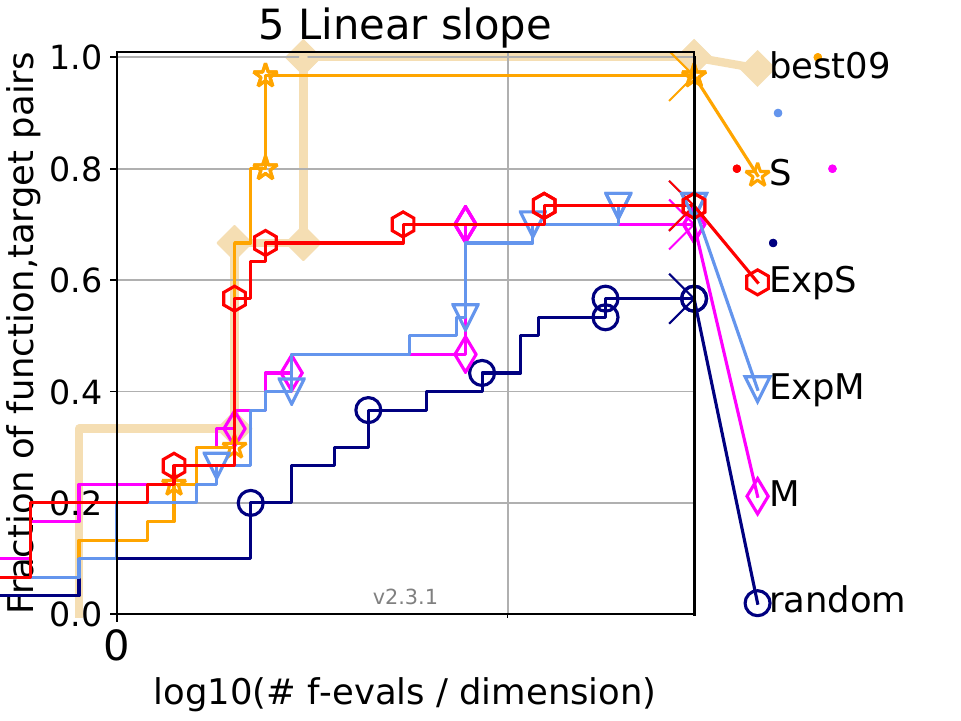}
\hfill
\includegraphics[width=0.33\textwidth]{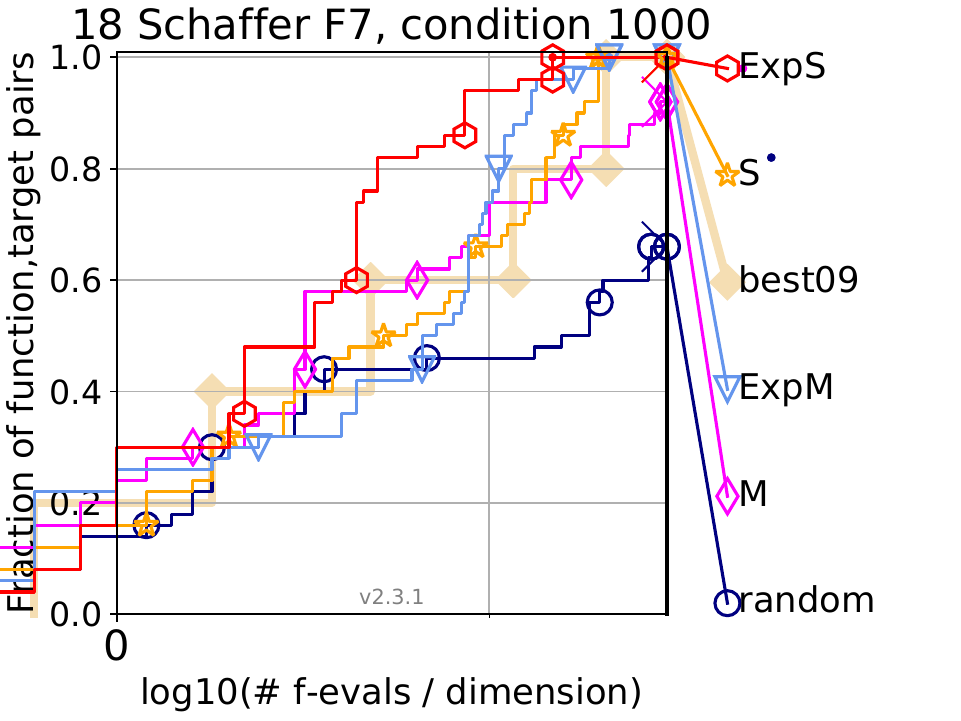}
\hfill
\caption{
ERTD for f5 (Linear Slope, left) and f18 (Shaffer F7 with a condition of 1000, right) in $d=5$. 
Algorithms with ``Exp'' in the name use an anisotropic exponential kernel, the others use an anisotropic Mat\'ern 5/2; S and M designate a small and a medium initial DoE. See Table~\ref{tab:allruns} for further details.
\label{fig-ERTD_expon_f5f18}
}
\end{figure}

Smooth functions, in particular the Linear Slope (f5) and the Ellipsoid (f2) 
are not well optimized through a GP with exponential kernel.
Vice versa, most multi-modal functions with adequate global structure (f15 to f19)
plus the Rastrigin function (f3) which is unimodal but wavy match well the exponential kernel. 
Illustrations are provided in Fig.~\ref{fig-ERTD_expon_f5f18} for the Linear Slope (f5, left) and Shaffer F7 with a moderate condition number (f18, right).

This behaviour of EGO with an exponential kernel is explained in part by the regression performance of GPs with exponential kernel for these functions, as it is detailed in Appendix \ref{sec-GP_Q2_optim}. 
The exponential kernel and its non-differentiable GP means does not competitively learn very smooth functions (ellipsoids, linear) but it is a very good regressor for many multimodal functions for which it also yields a high performing EGO (f16,f17,f18,f24). 
As analyzed in the Appendix, although there is a clear correlation between regression and optimization performances in particular when comparing Mat\'ern and exponential kernels, the regression ability is not the only contribution to optimization performance. 
EGO with an exponential kernel most often overperforms with respect to its regression qualities (on f3, f9, f12, f15 and f19) but it can underperform (on f21 and f22). 
Explanations other than the GP mean accuracy should be sought. For example, exponential kernels decrease the conditional variance in the vicinity of data points (at given length scales) more slowly than Mat\'ern kernels do, which might compensate for an inherent lack of exploration of the EI criterion on multimodal functions.

The interaction between the kernel, scaling, and warping was also studied and found not to be beneficial. To preserve the readability of the curves, it is reported in Appendix~\ref{sec-complement_expScalWarp}.

\subsubsection{Trend}
\begin{figure}
    \centering
\begin{tabular}{cc@{}c@{}c@{}}
 &{Separable} & {Low conditioning} & {High conditioning} \\
 \rotatebox[origin=l]{90}{\textbf{Trend}} & 
\includegraphics[trim=0mm 0mm 10mm 0mm, clip, width=0.3\textwidth]{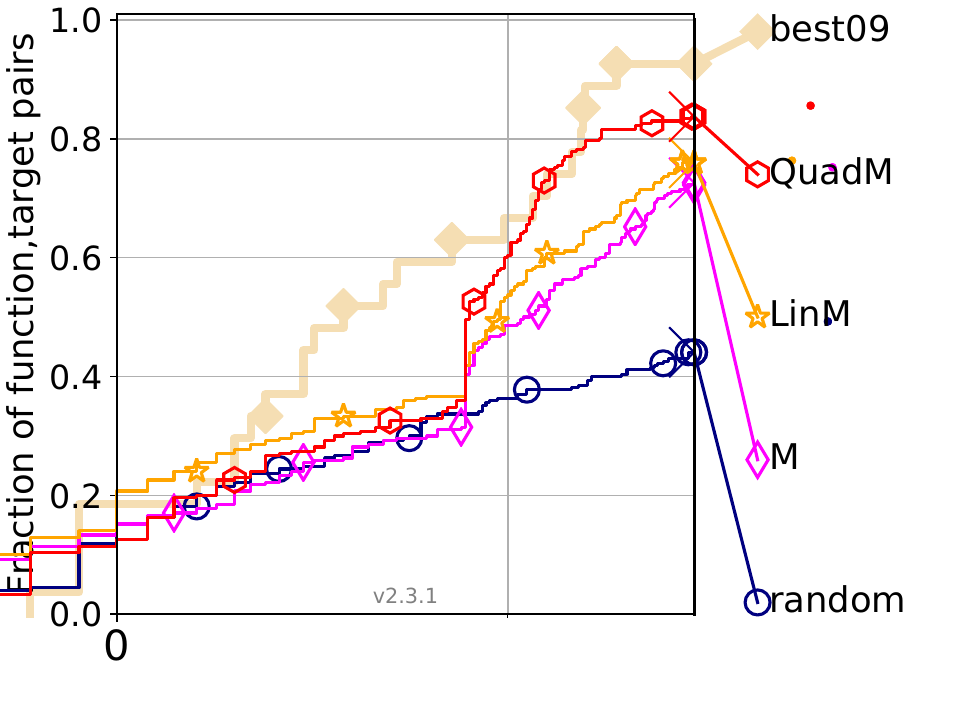} &
\includegraphics[trim=0mm 0mm 10mm 0mm, clip, width=0.3\textwidth]{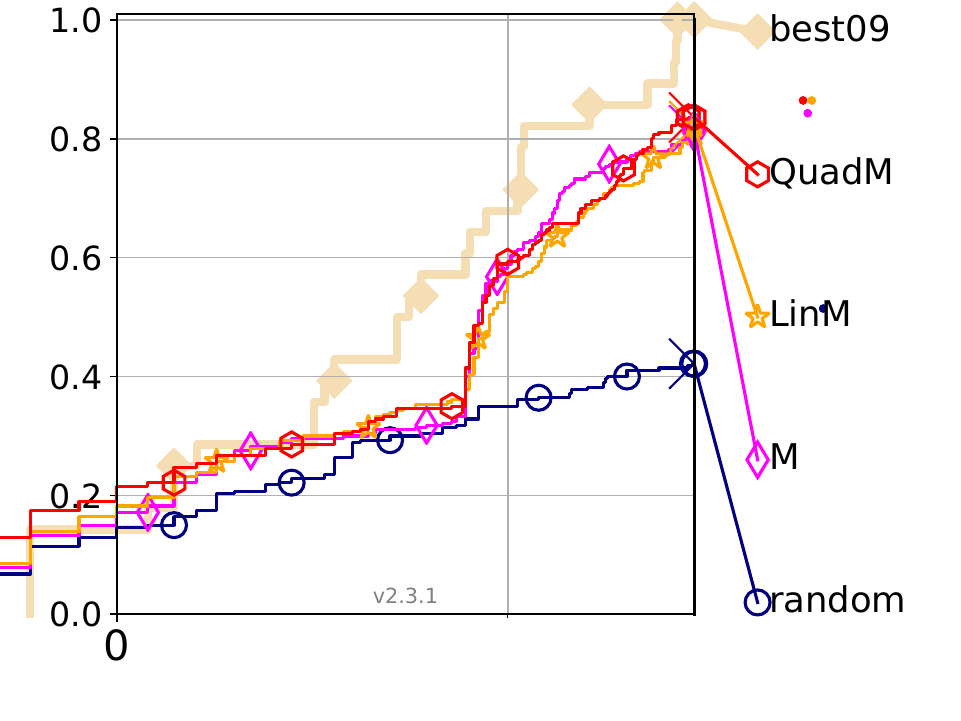}&
\includegraphics[trim=0mm 0mm 10mm 0mm, clip, width=0.3\textwidth]{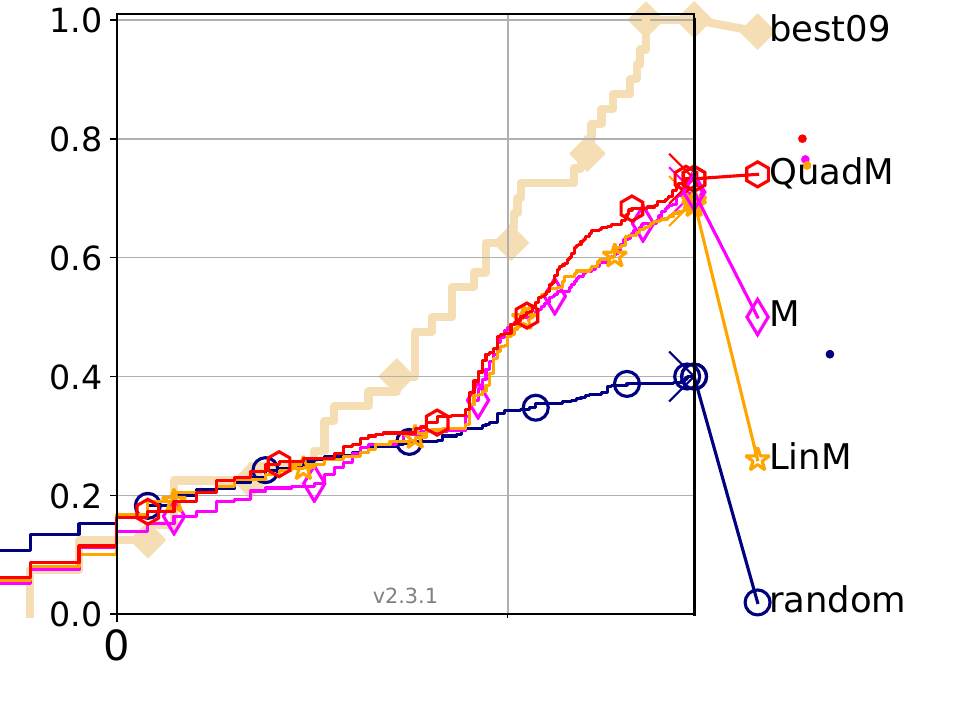} \\
& {Multimod., strong struct.} & {Multimod., weak struct.} & \\
&
\includegraphics[trim=0mm 0mm 10mm 0mm, clip, width=0.3\textwidth]{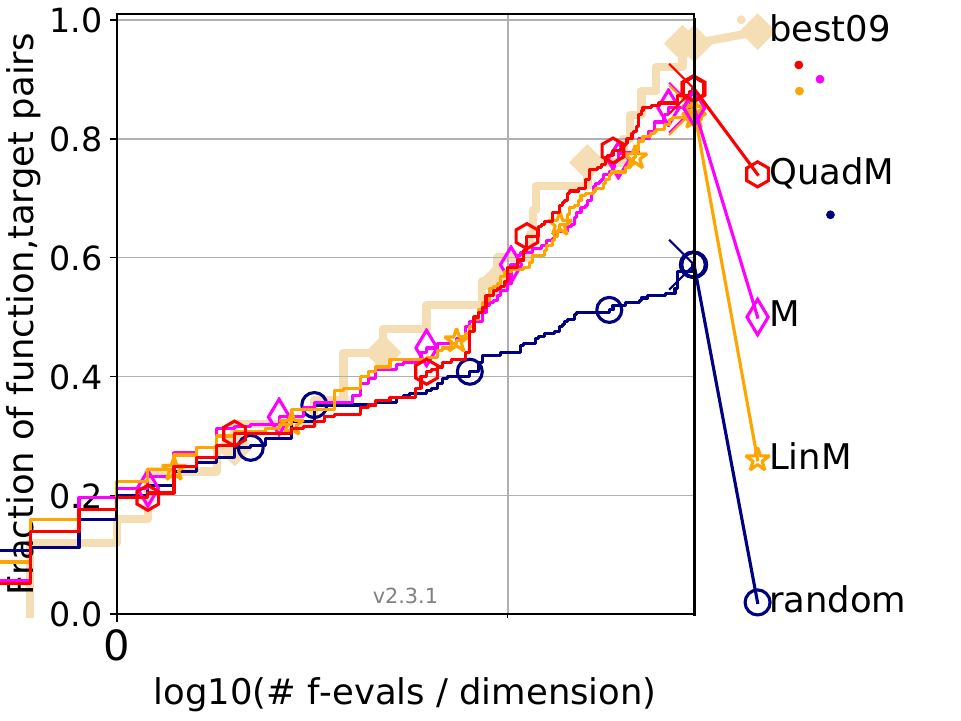}&
\includegraphics[trim=0mm 0mm 10mm 0mm, clip, width=0.3\textwidth]{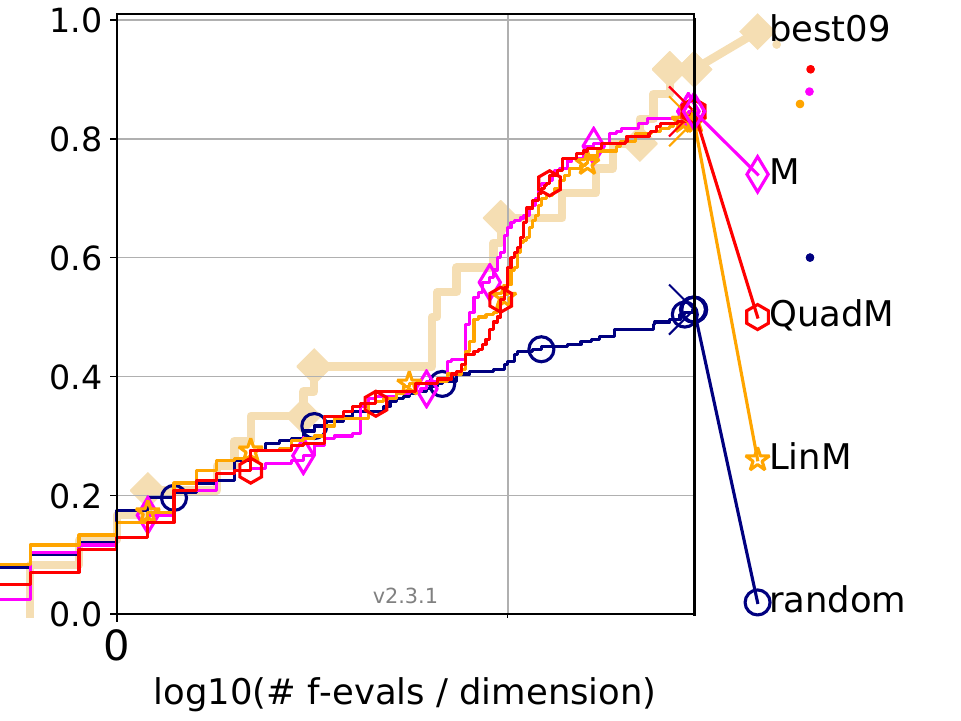}&
\end{tabular}
\caption{
ERTD by function groups in $d=5$ dimensions with varying trends : M (constant), LinM (linear), QuadM (quadratic).
\label{fig-ERTD_trend_05D_groups}
}
\end{figure}

Figure \ref{fig-ERTD_trend_05D_groups} shows that 
the quadratic trend largely outperforms alternatives on the first group but not on the others. 
As this group is not the main target for BO, it is understandable that most common practices do not include quadratic trends.
However, since the quadratic trend does not harm performance on the other groups, this appears as a very beneficial option, in particular for problems 
for which little prior information on their structure or difficulty is available.

\subsection{Testing EGO variants}

\subsubsection{Global analysis}

\begin{figure}
\begin{tabular}{cc@{}c@{}c@{}}
& \textbf{3D} & \textbf{5D} & \textbf{10D} \\
\rotatebox[origin=l]{90}{\textbf{ Input scaling}} & 
\includegraphics[width=0.3\textwidth]{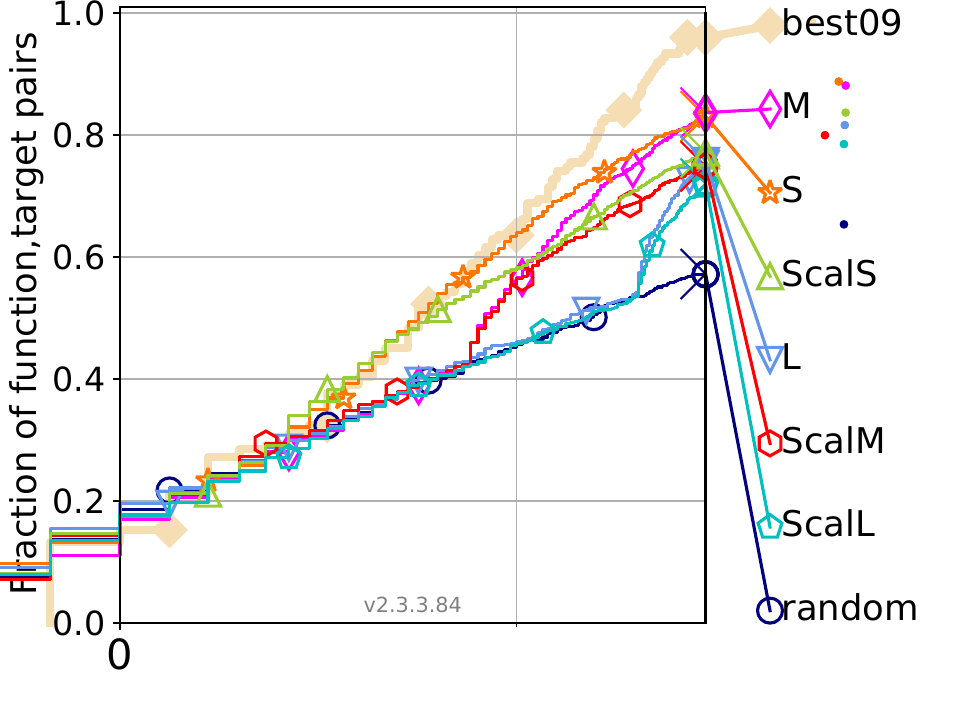}&
\includegraphics[width=0.3\textwidth]{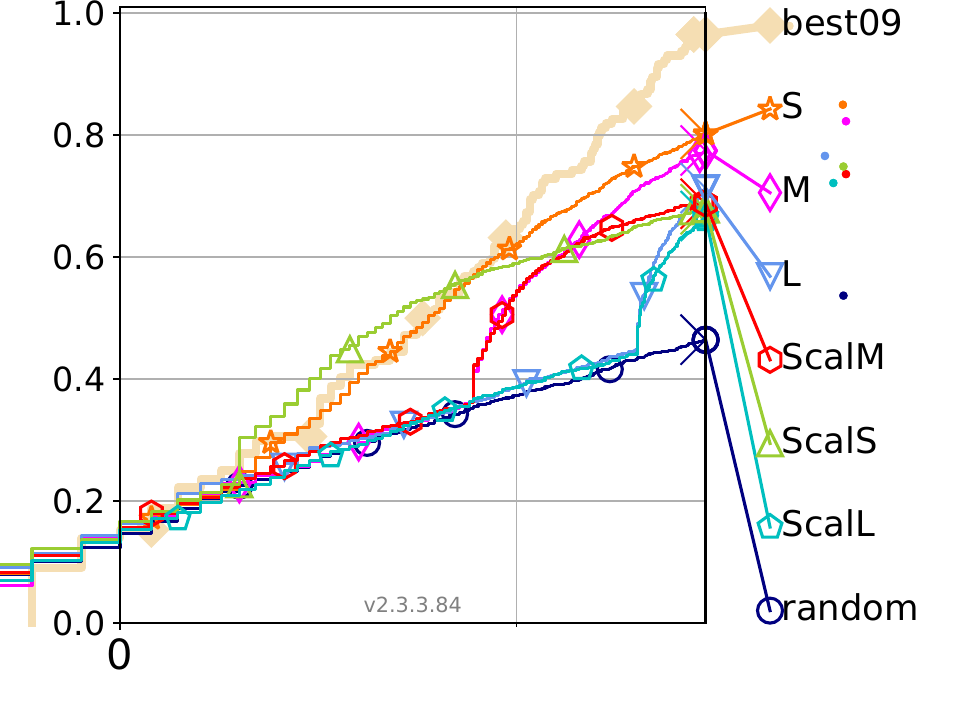}&
\includegraphics[width=0.3\textwidth]{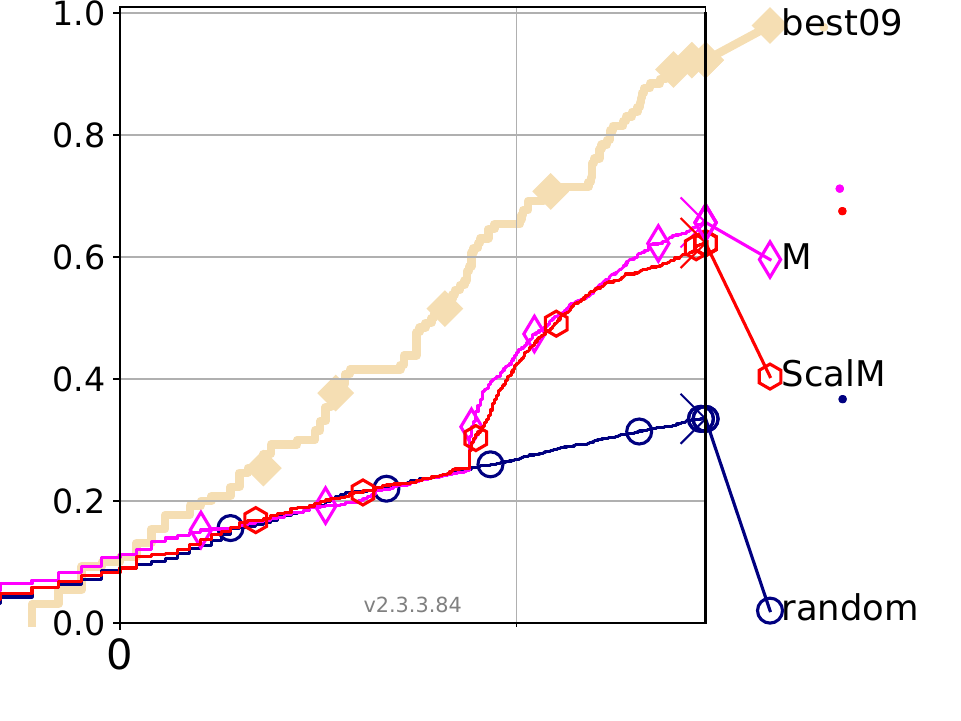}\\
\rotatebox[origin=l]{90}{\textbf{ Output warping}} & 
\includegraphics[width=0.3\textwidth]{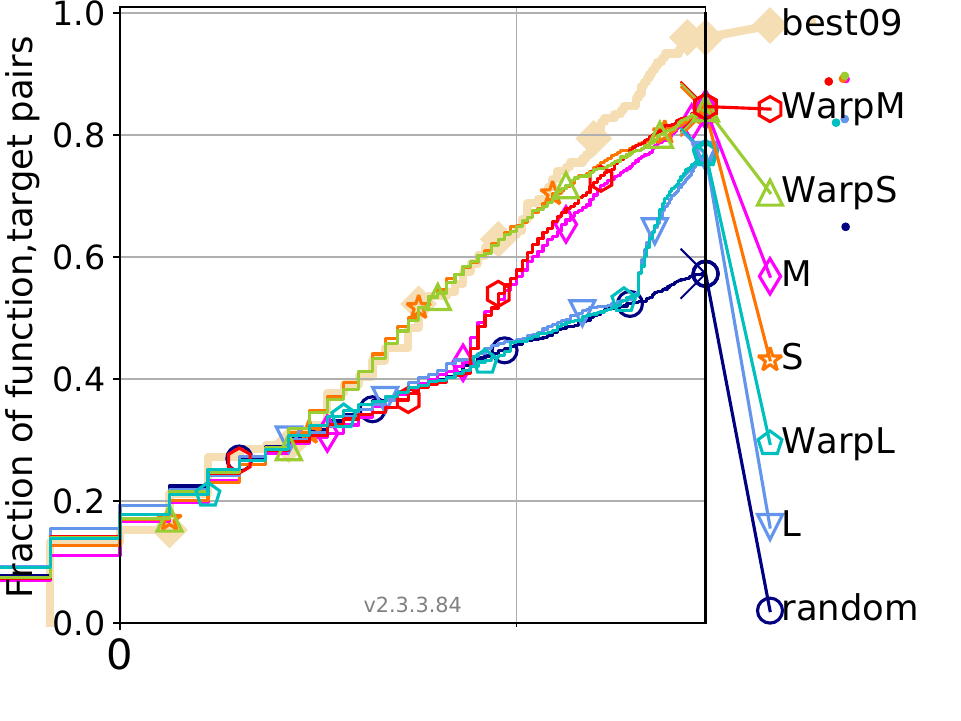}&
\includegraphics[width=0.3\textwidth]{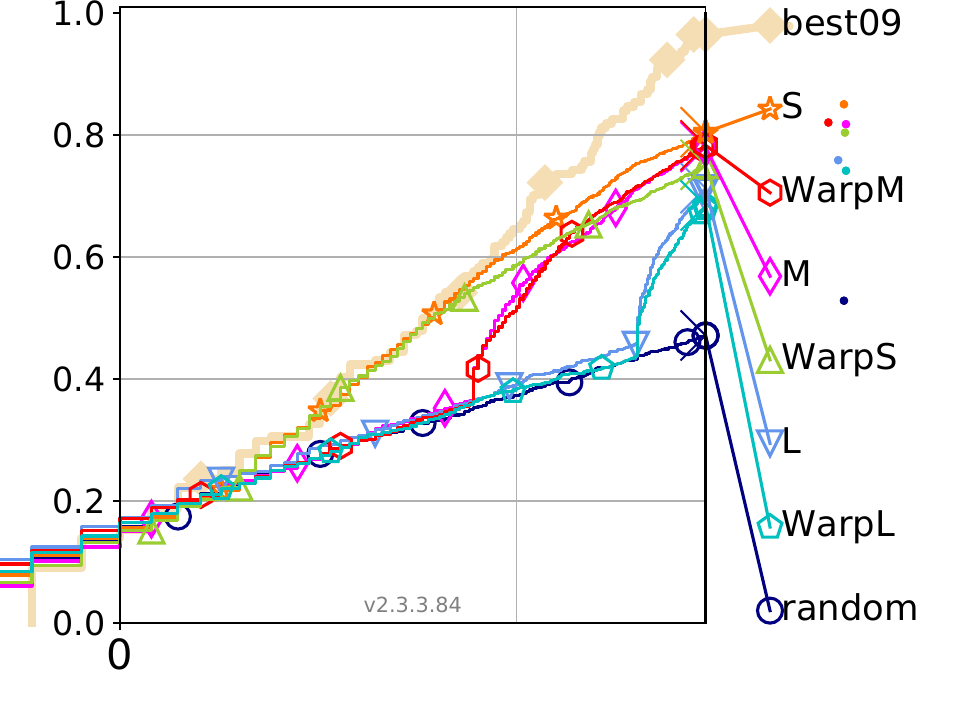}&
\includegraphics[width=0.3\textwidth]{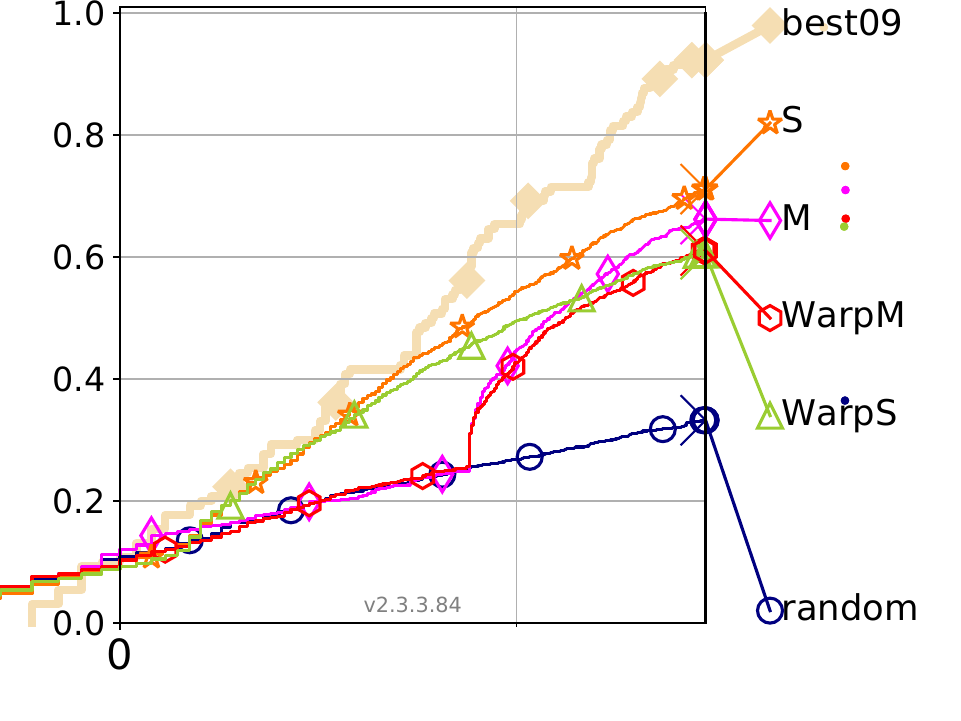}\\
\rotatebox[origin=l]{90}{\textbf{ GP mean}} & 
\includegraphics[width=0.3\textwidth]{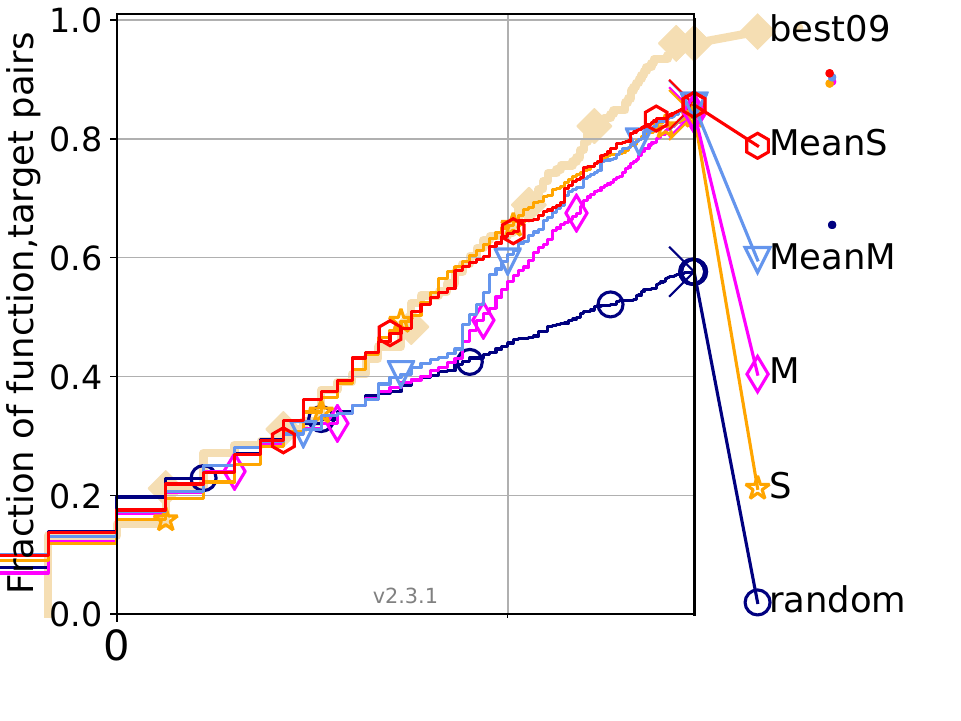}&
\includegraphics[width=0.3\textwidth]{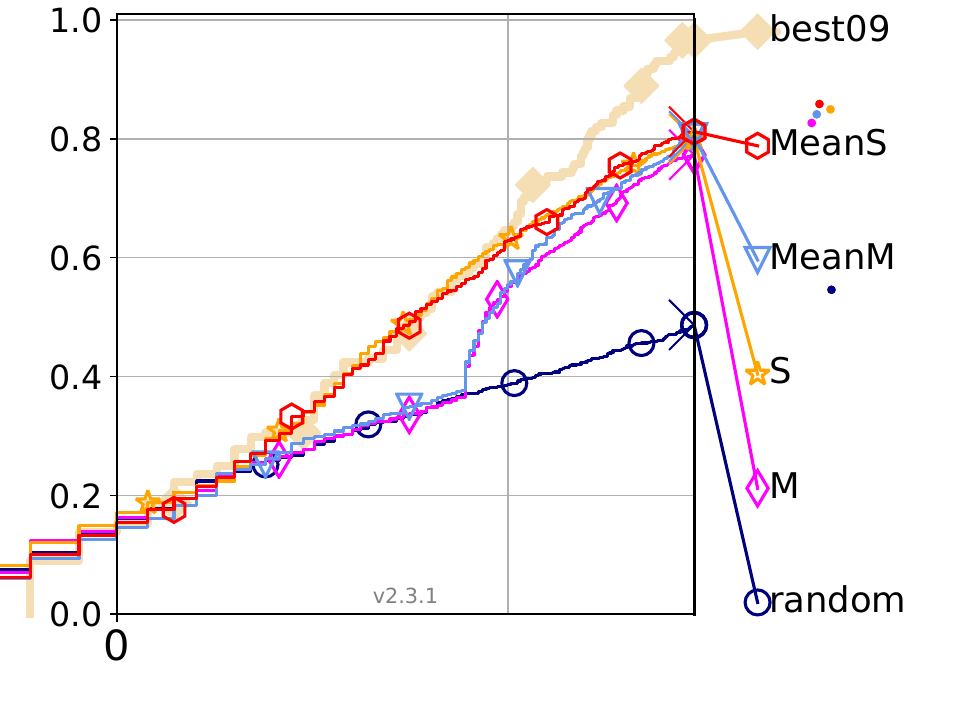}&
\includegraphics[width=0.3\textwidth]{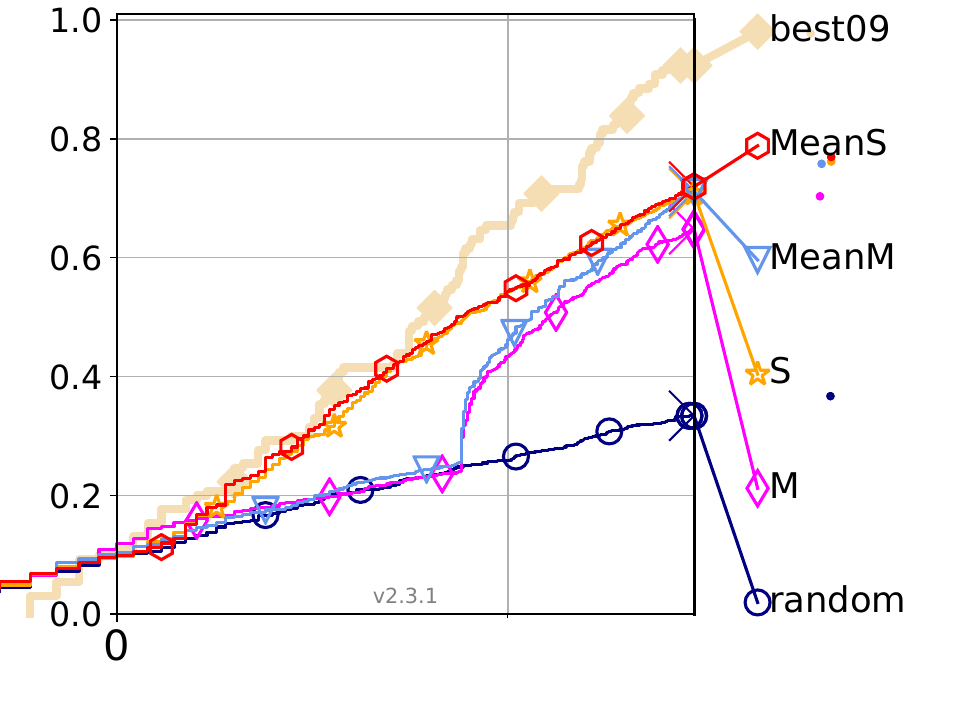}\\
\rotatebox[origin=l]{90}{\textbf{ EI optimization}} & 
\includegraphics[width=0.3\textwidth]{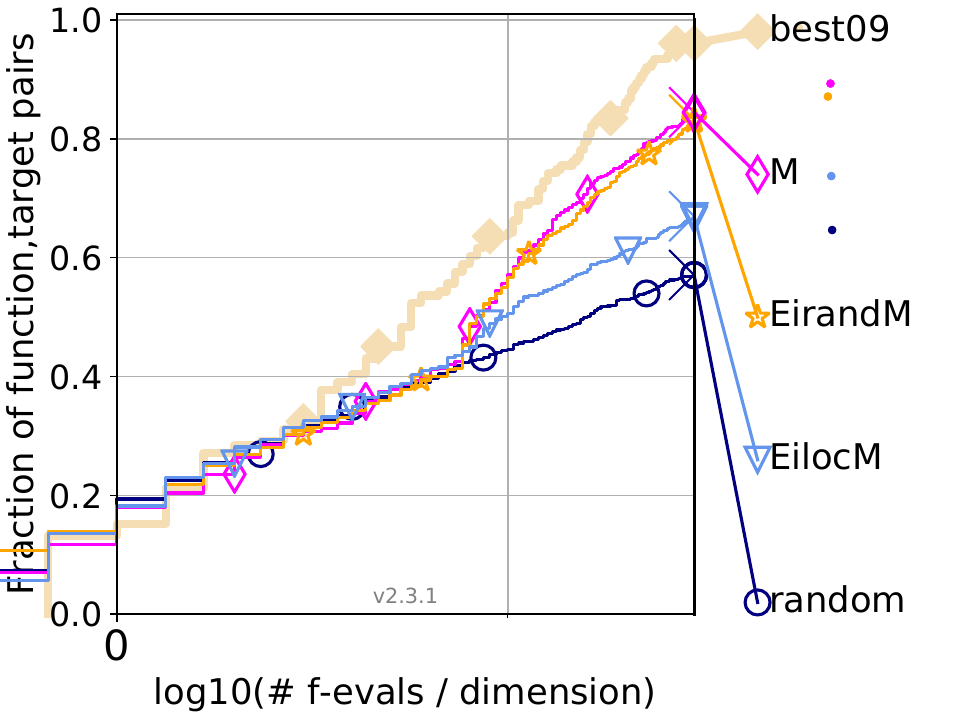}&
\includegraphics[width=0.3\textwidth]{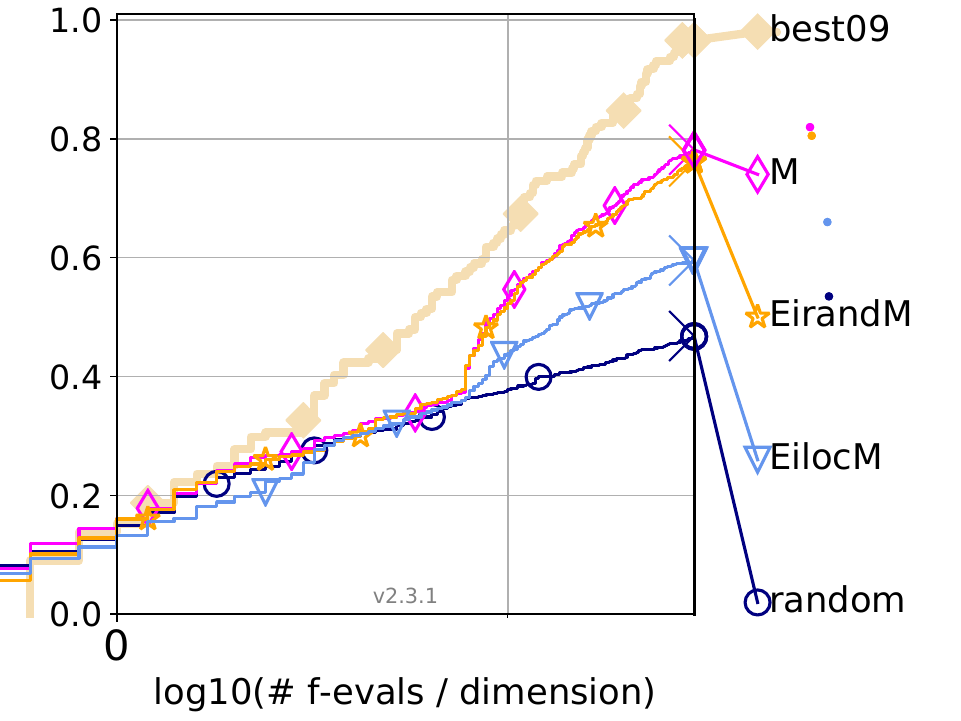}&
\includegraphics[width=0.3\textwidth]{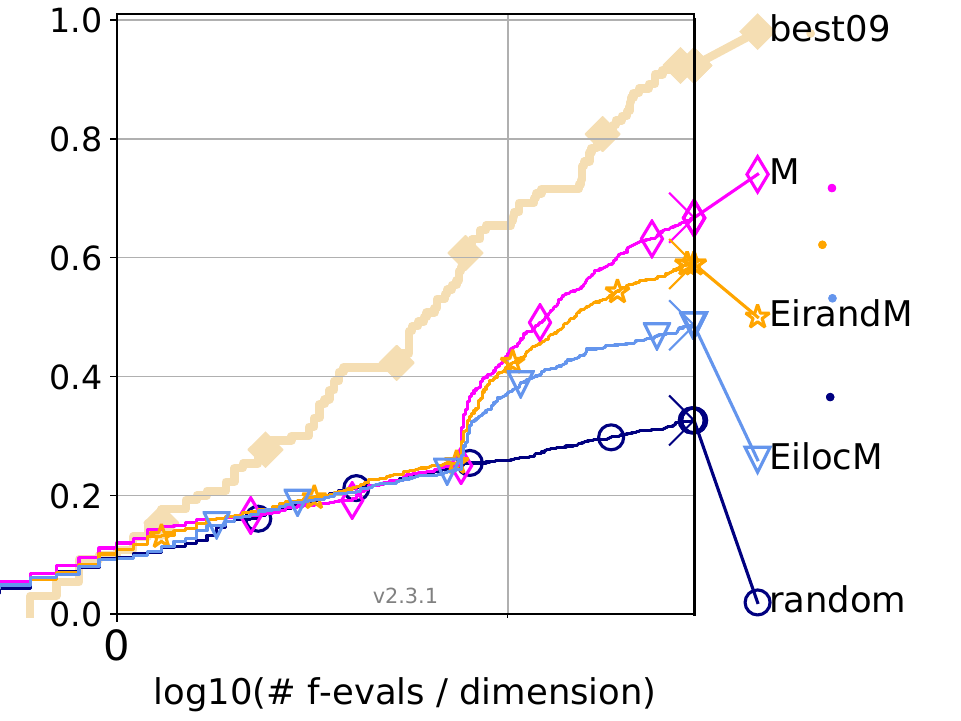}
\end{tabular}
\caption{ERTD over all functions in $d=3,5$ and $10$ dimensions.
Top row: with and without input scaling : ScalS/M/L (scaling and small/medium/large initial DoE). 
Second row: with and without output warping : WarpS/M/L (warping and small/medium/large initial DoE).
Third row: with and without the GP mean acquisition: MeanS/M (with small/medium initial DoE).
Bottom row: varying EI optimization schemes : EirandM by random search, EilocM by local search, M by repeated local searches from random initial points.
\label{fig-ERTD_ablation}
}
\end{figure}

Figure~\ref{fig-ERTD_ablation} shows the ERTD over all functions in $d=3,5$ and $10$ dimensions, for different configurations of the input scaling, output warping, 
GP mean acquisition or EI optimization.

\paragraph{Input scaling.}
Since it may benefit from large data sets, it is studied jointly with the initial DoE size for $d=3, 5$.
Overall, we can see a clear negative effect of scaling regardless of the dimension and initial sample size. 
Nevertheless, early in the search scaling with a small initial DoE is beneficial in 5 dimensions and marginally in 3 dimensions. 
This early asset of scaling and small population soon disappears and it ends up slowing down convergence.

\paragraph{Output warping.}
Similarly to scaling, it may benefit from large data sets, and it is studied jointly with the initial DoE size for $d=3, 5$.
No significant effect can be seen for $d=3$, but warping affects negatively performance for $d=5,10$, in particular for the small initial budget.

\paragraph{GP mean acquisition.}
Overall, this occasional GP mean acquisition slightly improves the results in all cases, in particular towards the end of the runs, in high dimension (5 and 10) and when the initial DoE has medium size. 

\paragraph{EI optimization.}
The single BFGS run (EilocM) clearly underperforms at all times for all dimensions. 
Surprisingly, random search (EirandM) is almost comparable to multi-start BFGS in dimensions 3 and 5.
Significant differences only appear for $d=10$.

As in the previous subsection, we next detail our analysis by breaking the ERT distributions by functions groups, in dimension 5 only.

\subsubsection{Scaling (input non-linear rescaling)}
\begin{figure}
    \centering
\begin{tabular}{cc@{}c@{}c@{}}
 &{Separable} & {Low conditioning} & {High conditioning} \\
 \rotatebox[origin=l]{90}{\textbf{Input scaling}} & 
\includegraphics[trim=0mm 0mm 10mm 0mm, clip, width=0.3\textwidth]{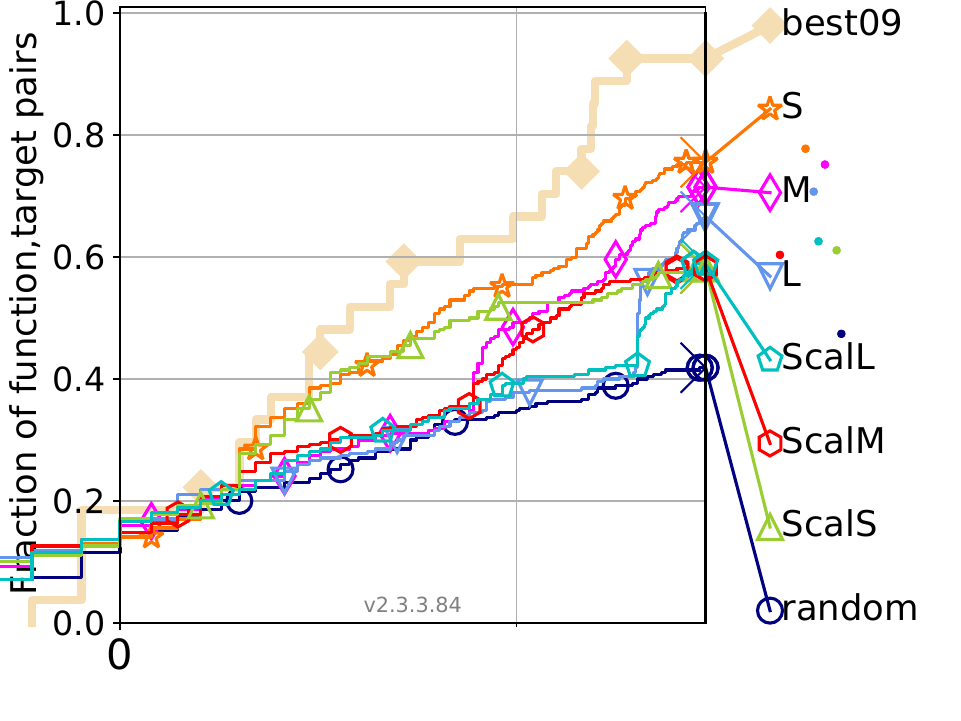} &
\includegraphics[trim=0mm 0mm 10mm 0mm, clip, width=0.3\textwidth]{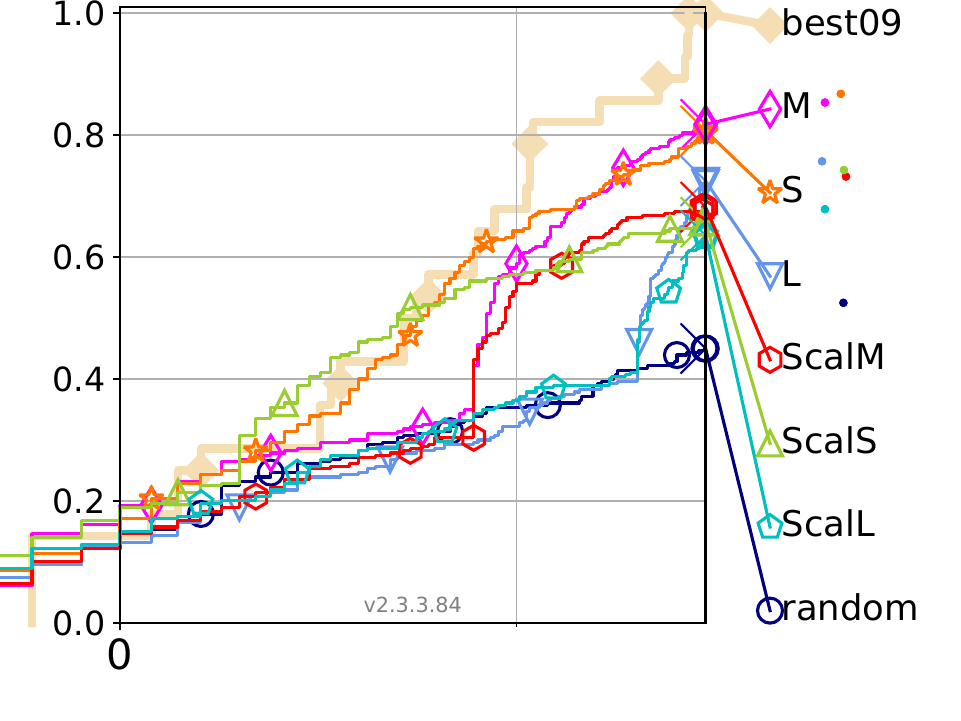}&
\includegraphics[trim=0mm 0mm 10mm 0mm, clip, width=0.3\textwidth]{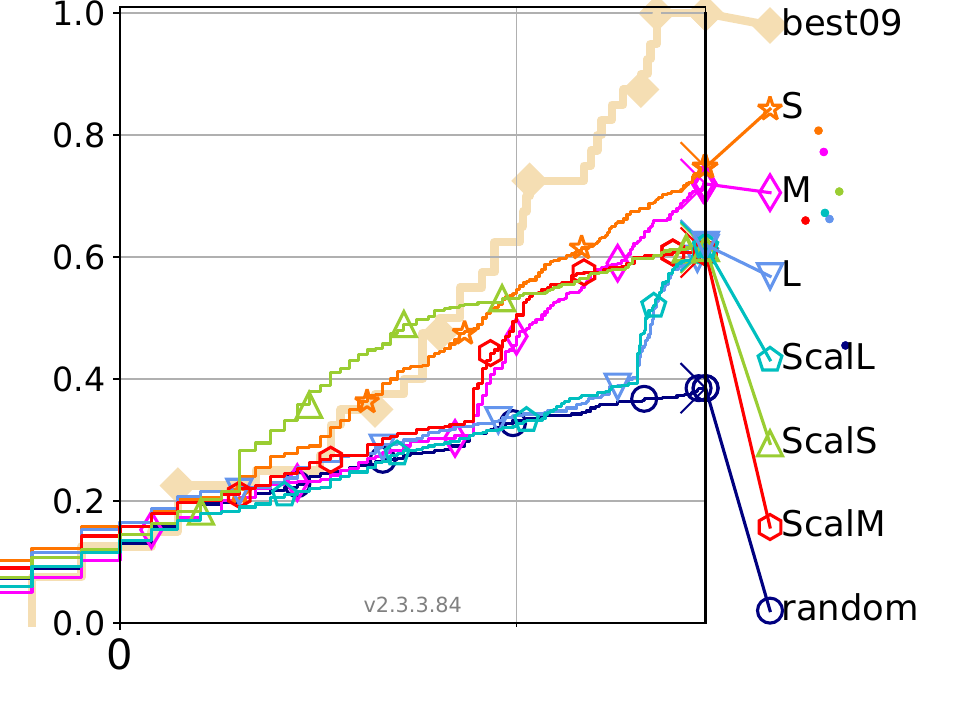} \\
& {Multimod., strong struct.} & {Multimod., weak struct.} & \\
&
\includegraphics[trim=0mm 0mm 10mm 0mm, clip, width=0.3\textwidth]{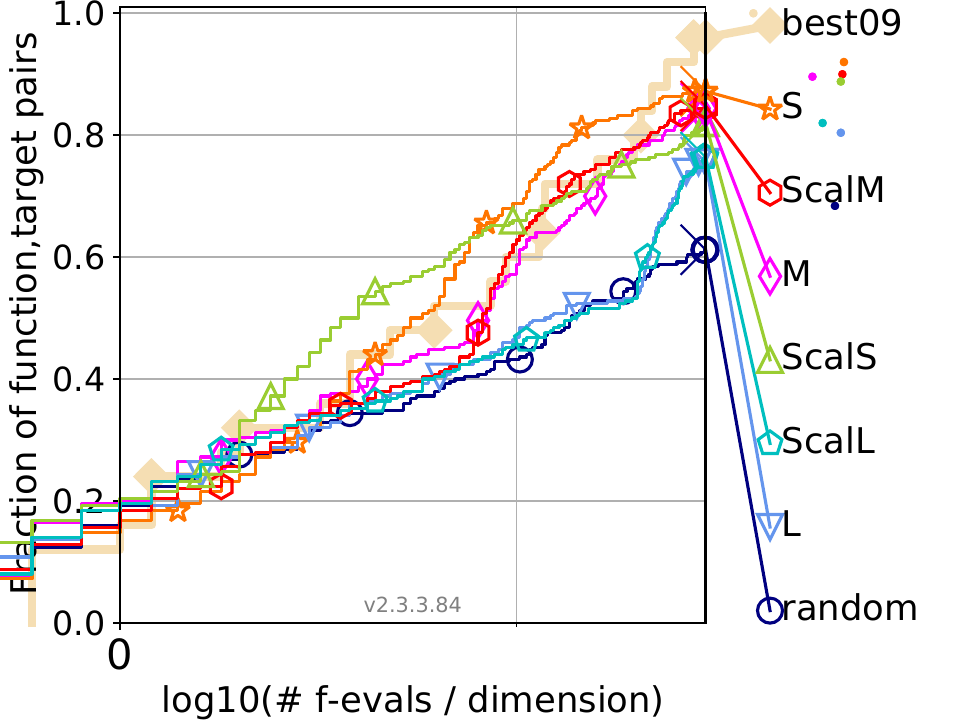}&
\includegraphics[trim=0mm 0mm 10mm 0mm, clip, width=0.3\textwidth]{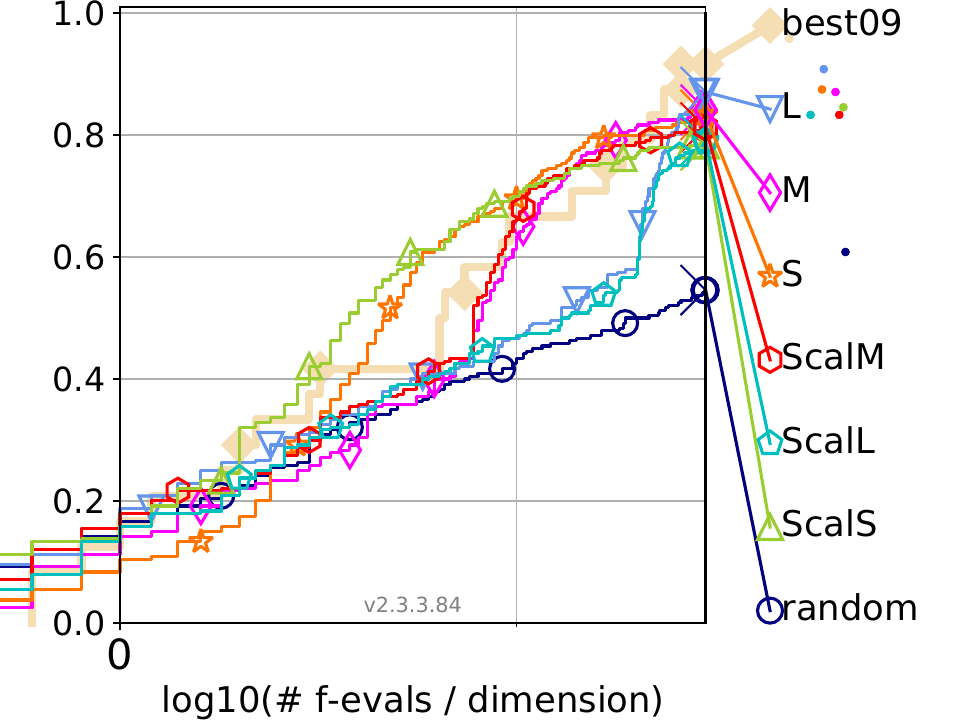}&
\end{tabular}
\caption{
ERTD by function groups in $d=5$ dimensions with and without scaling : ScalS/M/L = scaling with small/medium/large initial DoE.
\label{fig-ERTD_scaling_05D_groups}
}
\end{figure}

Figure \ref{fig-ERTD_scaling_05D_groups} shows an early benefit of scaling with a small initial DoE in all groups but the separable functions.
However, scaling then considerably slows down convergence, so that it becomes rapidly less efficient that the regular runs. 
Hence, it seems that scaling helps improving the GP quality thus progress at the early steps. 
However, later in the search, scaling prevents the algorithm from locating small, high-performance, level sets, possibly due to overfitting issues. 
For larger budgets, for which one might expect scaling to be more stable and efficient, the effect is either null or negative. 
In general, this approach is not beneficial on this benchmark and might be limited to specific cases.



\begin{figure}
\centering
\begin{tabular}{c@{}c@{}}
 {f19 in 3D} & {f13 in 5D} \\
\includegraphics[trim=0mm 0mm 0mm 8mm, clip, width=0.33\textwidth]{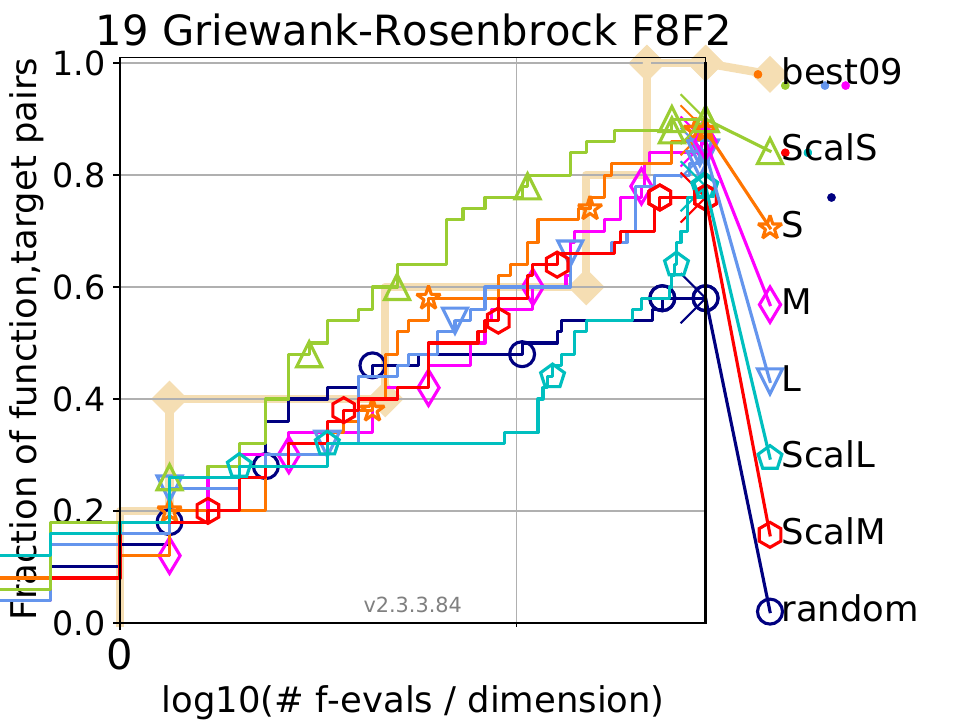} &
\includegraphics[trim=0mm 0mm 0mm 8mm, clip, width=0.33\textwidth]{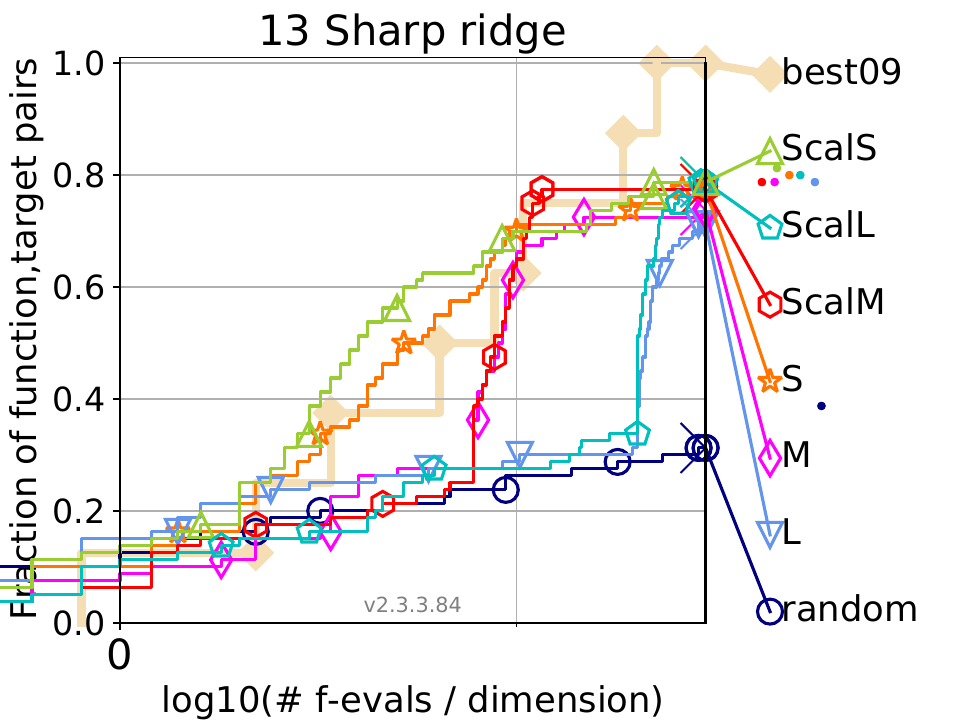}
\end{tabular}
\caption{
Examples of functions where scaling helps.
\label{fig-ERTD_scaling_2fcts}
}
\end{figure}

Two functions where scaling is beneficial are detailed in Fig.~\ref{fig-ERTD_scaling_2fcts}. The first function is f19, the composite Griewank-Rosenbrock function, in 3 dimensions. This function resembles the Rosenbrock function in a multimodal way. ScalS, the variant of the algorithm with scaling and a small initial DoE, is the best algorithm. The second function (right plot) is f13 in 5D i.e., the sharp ridge function. In both functions, the good points are distributed in a bent valley which is not aligned with the axes.


\subsubsection{Warping (output non-linear rescaling)}
\begin{figure}
    \centering
\begin{tabular}{cc@{}c@{}c@{}}
 &{Separable} & {Low conditioning} & {High conditioning} \\
 \rotatebox[origin=l]{90}{\textbf{Output warping}} & 
\includegraphics[trim=0mm 0mm 10mm 0mm, clip, width=0.3\textwidth]{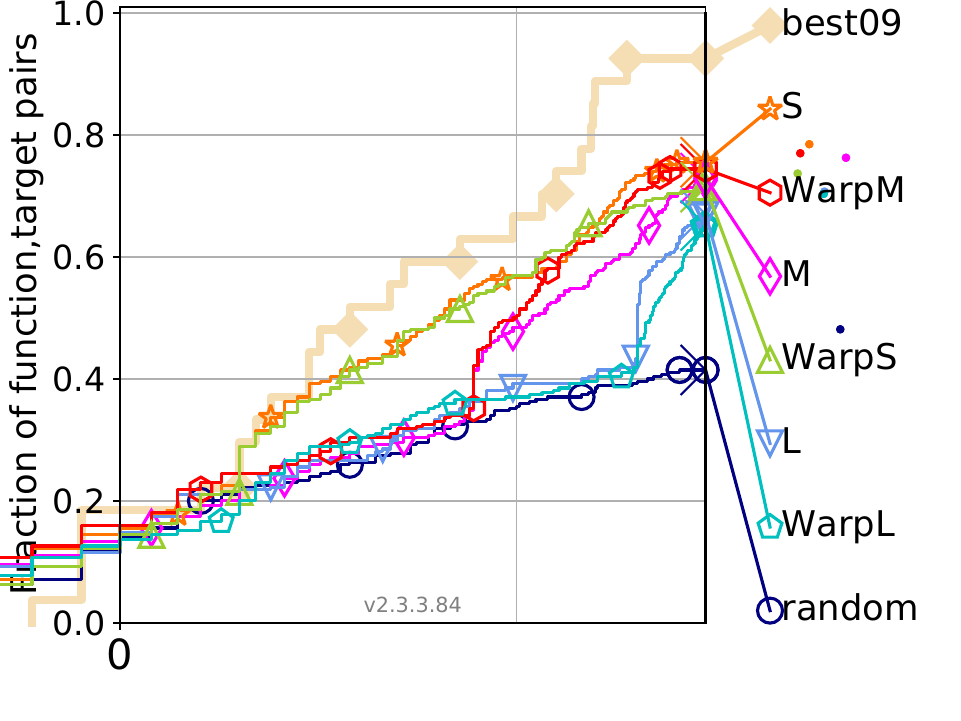} &
\includegraphics[trim=0mm 0mm 10mm 0mm, clip, width=0.3\textwidth]{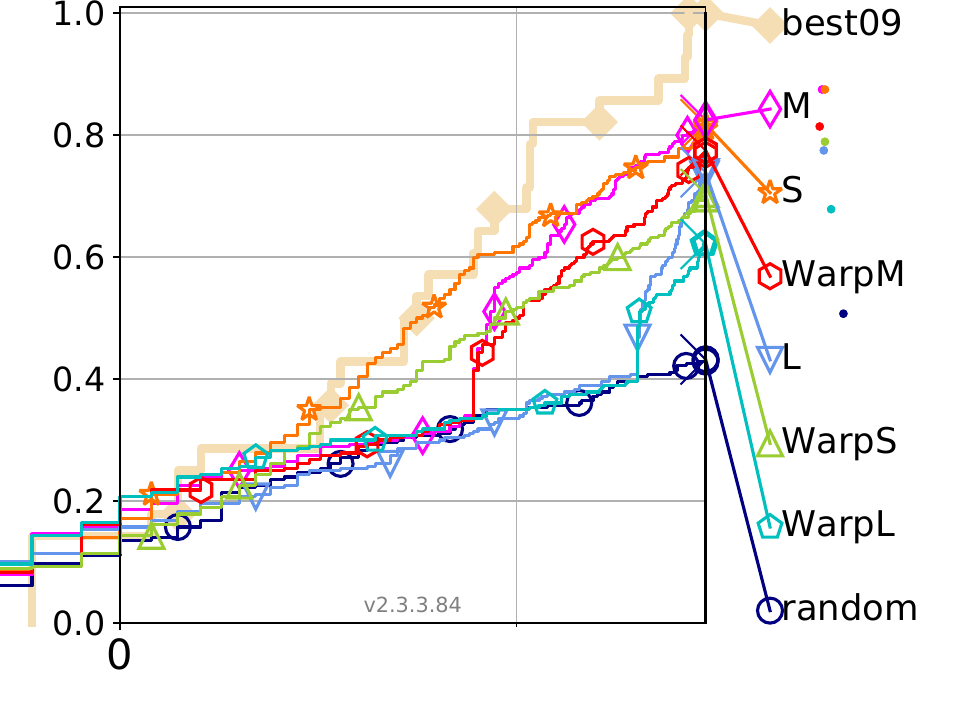}&
\includegraphics[trim=0mm 0mm 10mm 0mm, clip, width=0.3\textwidth]{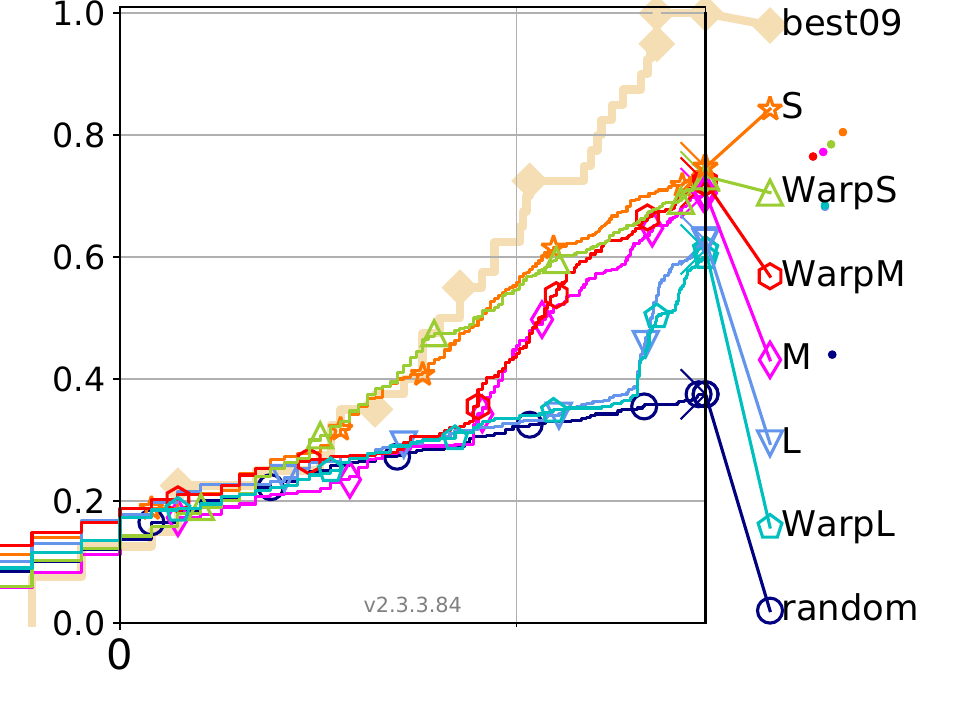} \\
 & {Multimod., strong struct.} & {Multimod., weak struct.} & {Multimod., weak struct., 10D} \\
&
\includegraphics[trim=0mm 0mm 10mm 0mm, clip, width=0.3\textwidth]{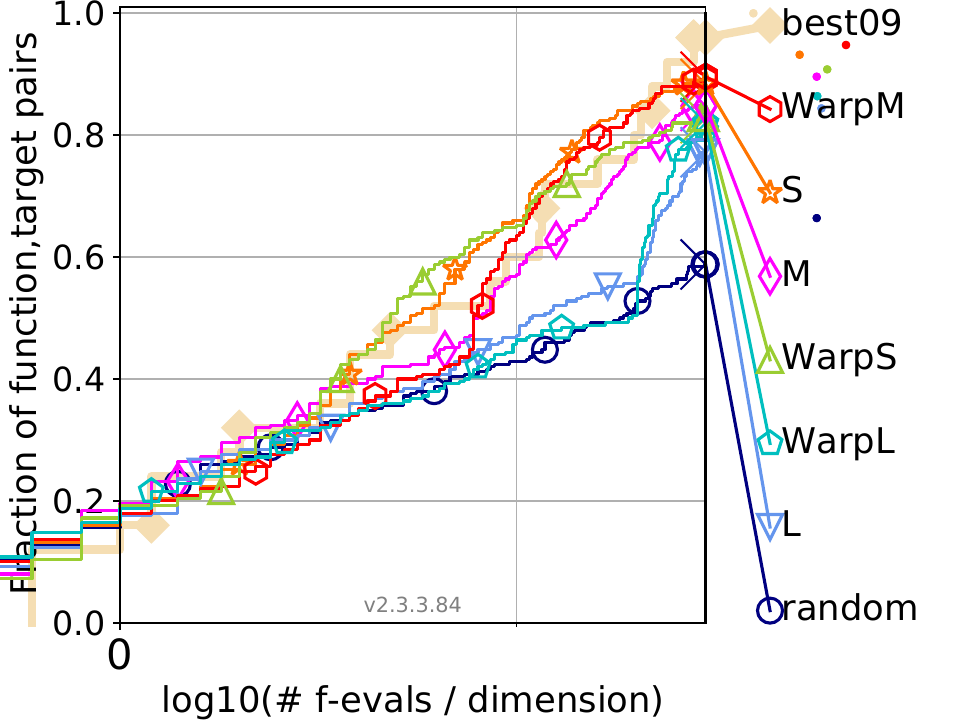}&
\includegraphics[trim=0mm 0mm 10mm 0mm, clip, width=0.3\textwidth]{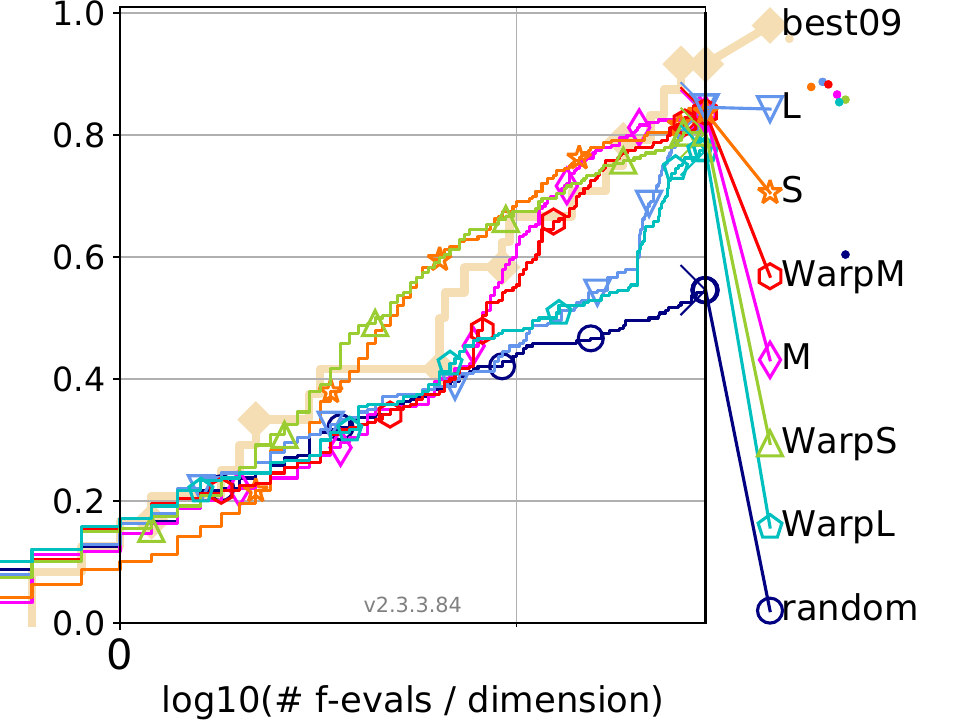}&
\includegraphics[trim=0mm 0mm 10mm 0mm, clip, width=0.3\textwidth]{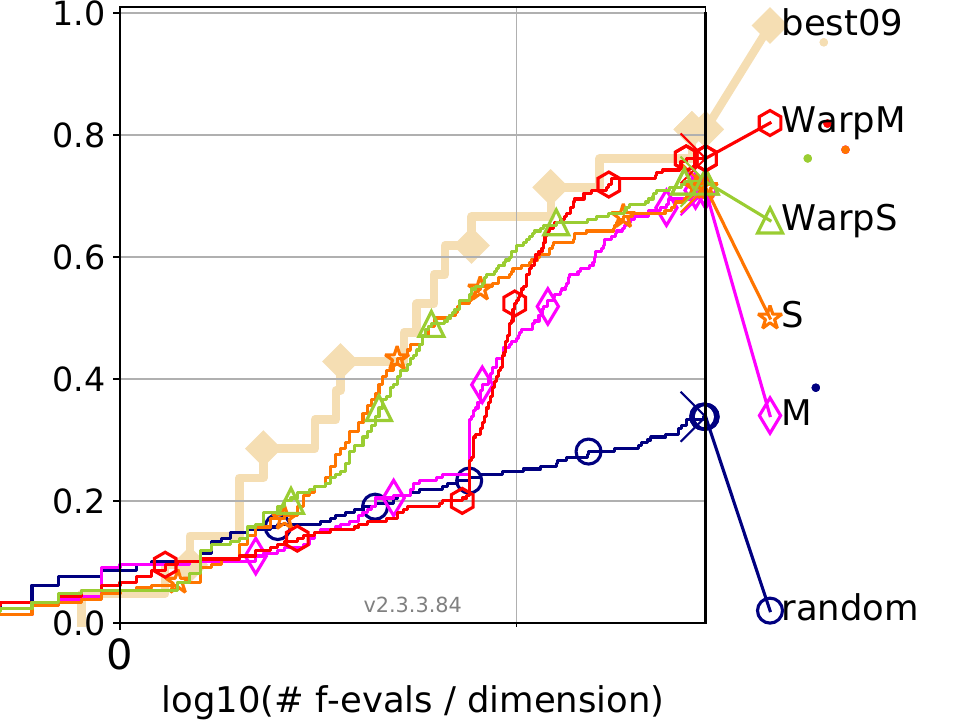}
\end{tabular}
\caption{
ERTD by function groups in $d=5$ dimensions with and without warping : WarpS/M/L = warping with small/medium/large initial DoE. Bottom right: group 5 in 10D.
\label{fig-ERTD_warping_05D_groups}
}
\end{figure}

Figure \ref{fig-ERTD_warping_05D_groups} shows that the overall mixed effects of warping are driven by a clear decline in performance on function group 2, i.e., the unimodal functions with low conditioning. In this case, there is no need to transform the outputs. 
The figure also confirms that warping does not accomodate well a small initial population: warpS is always inferior to S in 5 dimensions.
However, the medium initial budget (warpM) is a good setting for multimodal functions, in particular for those with low global structure in 10 dimensions (see bottom right of the Figure).

\begin{figure}
\centering
\begin{tabular}{c@{}c@{}}
 {f17 in 5D} & {f22 in 10D} \\
\includegraphics[trim=0mm 0mm 0mm 8mm, clip, width=0.33\textwidth]{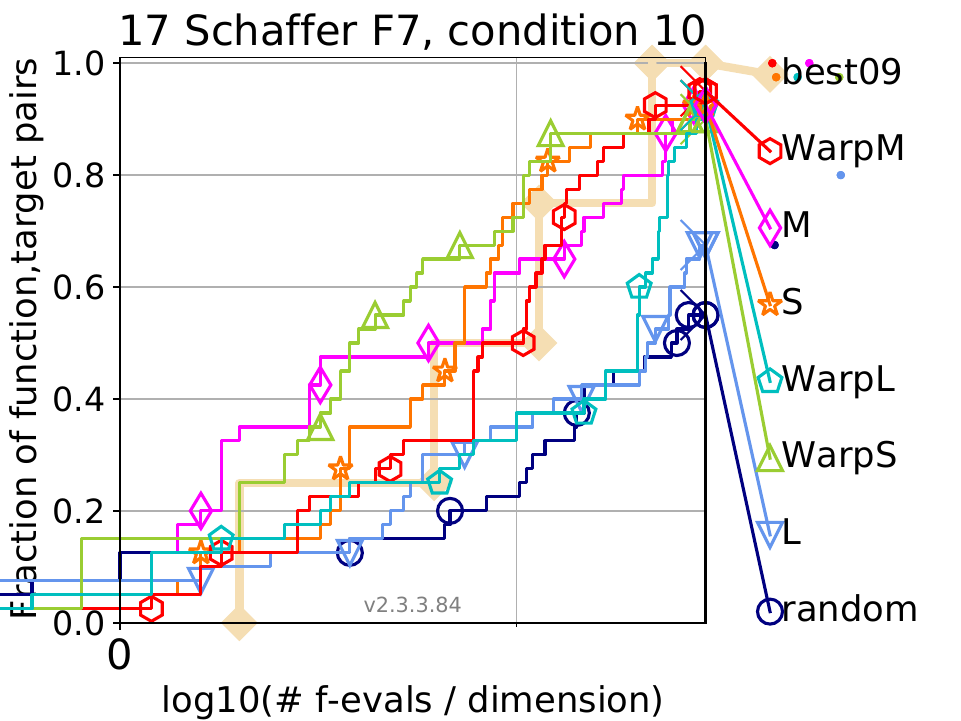} &
\includegraphics[trim=0mm 0mm 0mm 8mm, clip, width=0.33\textwidth]{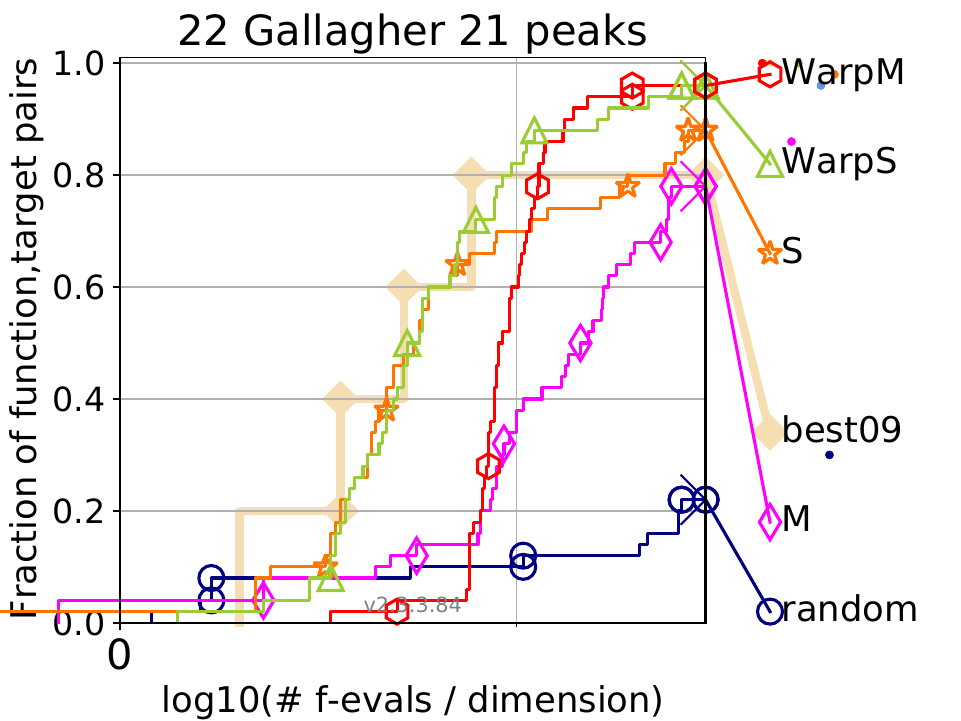} 
\end{tabular}
\caption{
Examples of functions for which warping is beneficial: f17 -- Schaffers F7 function -- in 5D (left) and f22 -- Gallagher's Gaussian 21 peaks -- in 10D (right).
\label{fig-ERTD_warp_good_fs}
}
\end{figure}

Examples of positive contributions of warping are reported in Fig.~\ref{fig-ERTD_warp_good_fs}.
The functions are the Schaffers F7 function (f17 in 5D) and Gallagher's Gaussian 21 Peaks function (f22 in 10D). Both functions are multimodal and non stationary.
While warpS, warping with a small DoE, is good early on, a medium DoE, warpM, is the best choice around 30$\times$d evaluations. 
Notice that on f22 in 10D warpM is significantly better than best09.

\subsubsection{GP mean acquisition }
\begin{figure}
    \centering
\begin{tabular}{cc@{}c@{}c@{}}
 &{Separable} & {Low conditioning} & {High conditioning} \\
 \rotatebox[origin=l]{90}{\textbf{GP mean acq.}} & 
\includegraphics[trim=0mm 0mm 10mm 0mm, clip, width=0.3\textwidth]{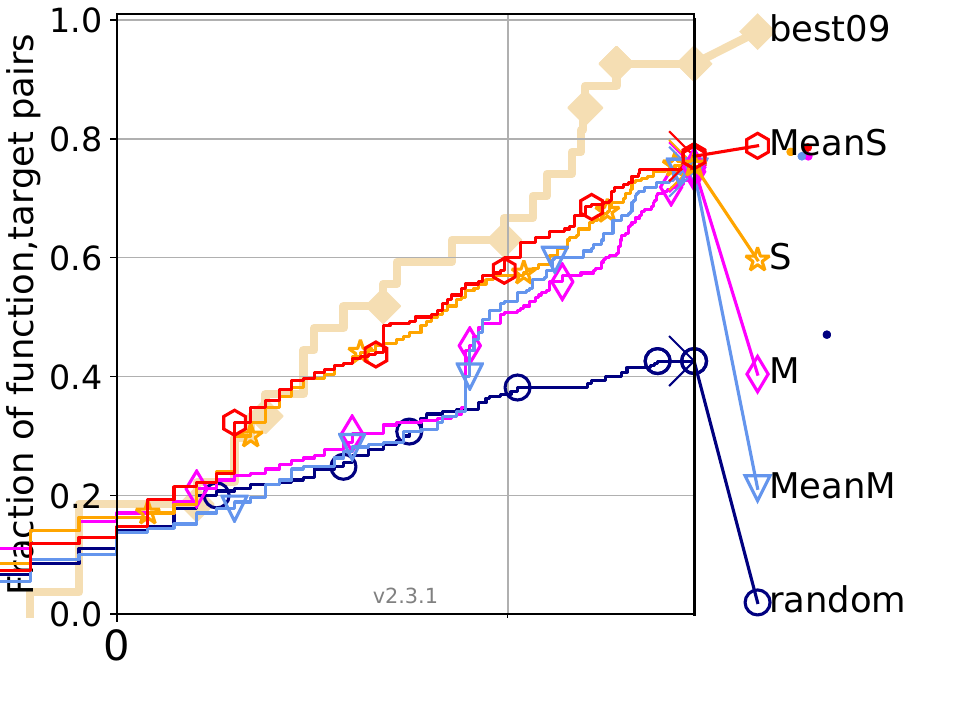} &
\includegraphics[trim=0mm 0mm 10mm 0mm, clip, width=0.3\textwidth]{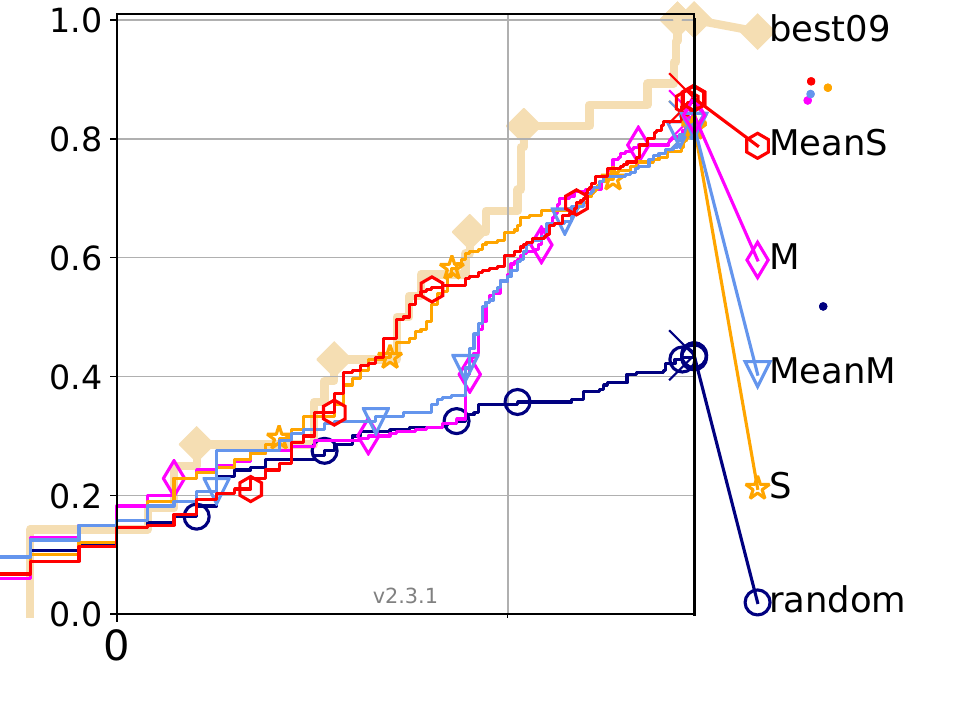}&
\includegraphics[trim=0mm 0mm 10mm 0mm, clip, width=0.3\textwidth]{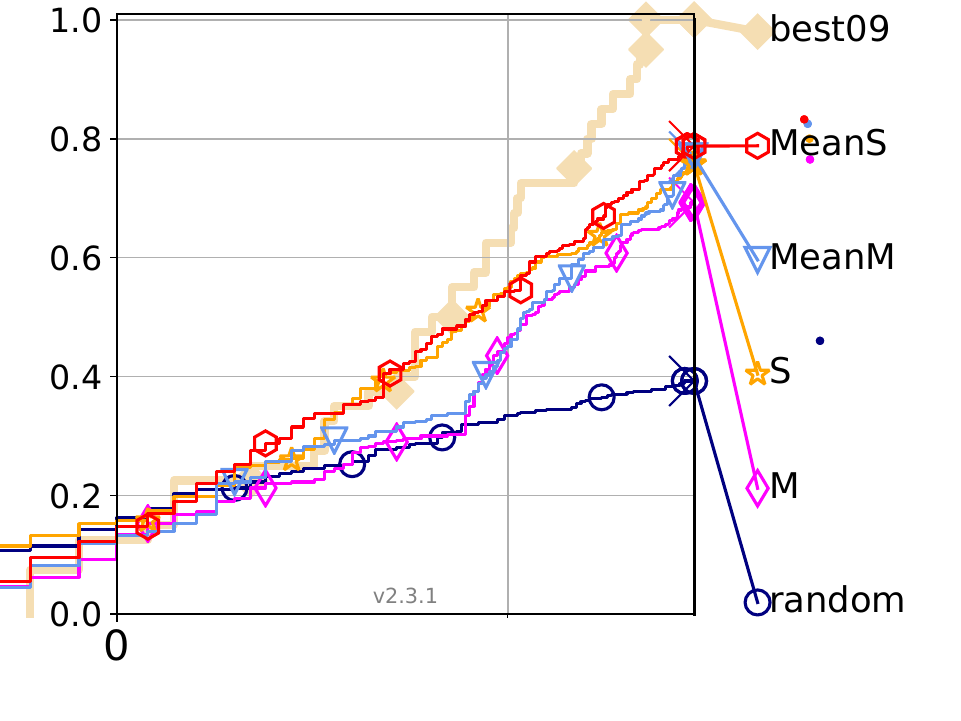} \\
 & {Multimod., strong struct.} & {Multimod., weak struct.} & ~ \\
&
\includegraphics[trim=0mm 0mm 10mm 0mm, clip, width=0.3\textwidth]{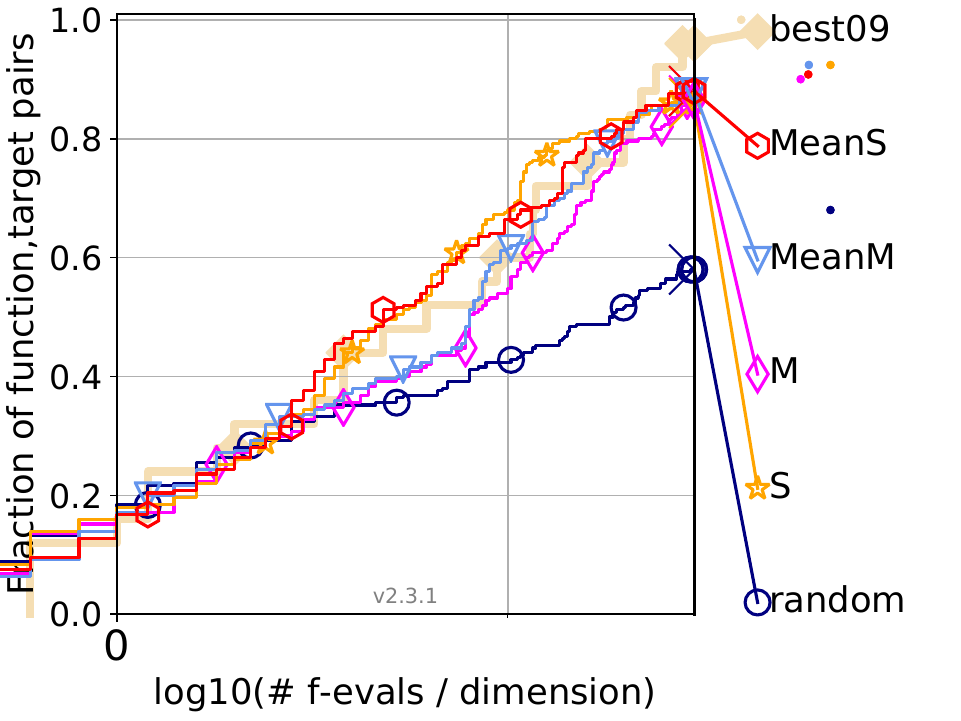}&
\includegraphics[trim=0mm 0mm 10mm 0mm, clip, width=0.3\textwidth]{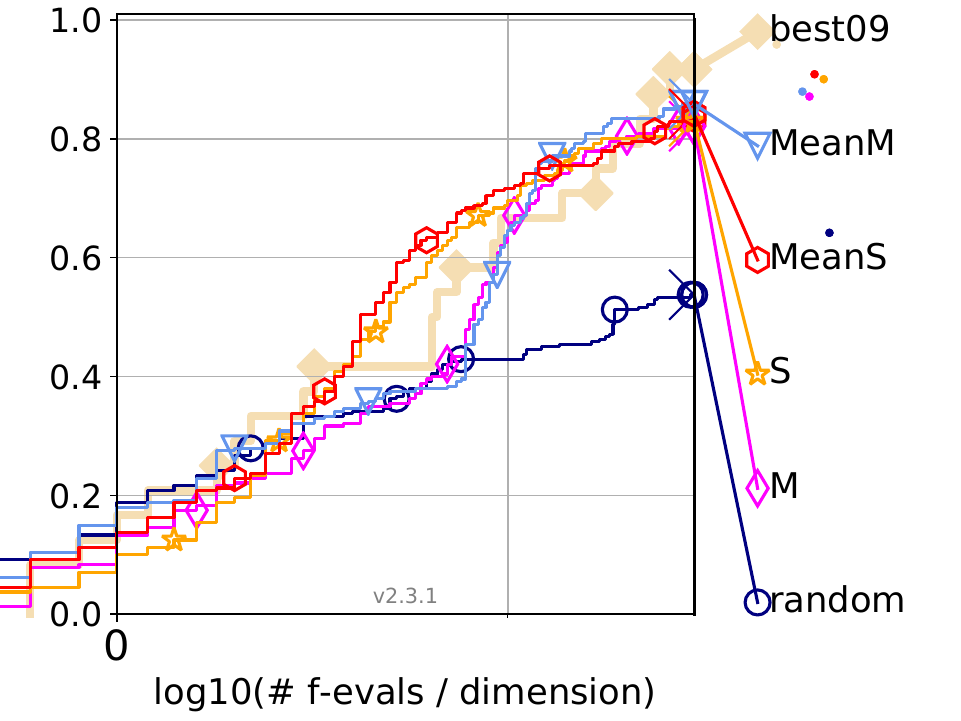}&
~
\end{tabular}
\caption{
ERTD by function groups in $d=5$ dimensions with and without GP mean acquisition : MeanS/M = GP mean acquisition with small/medium initial DoE. 
\label{fig-ERTD_gpmean_05D_groups}
}
\end{figure}

Considering functions groups (Fig.~\ref{fig-ERTD_gpmean_05D_groups}), it is observed that this strategy has a slightly positive effect on most cases, the most significant being on the third group (unimodal functions with high conditioning) towards late iterations.

\begin{figure}
\centering
\begin{tabular}{c@{}c@{}}
 {f13 in 5D} & {f16 in 10D} \\
\includegraphics[trim=0mm 0mm 0mm 8mm, clip, width=0.33\textwidth]{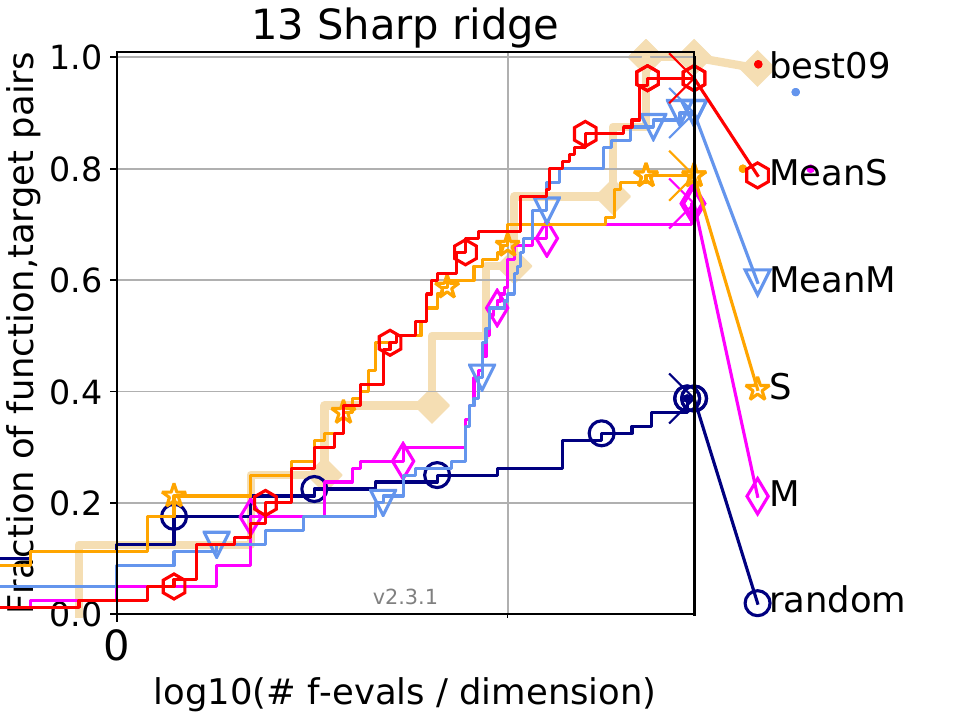} &
\includegraphics[trim=0mm 0mm 0mm 8mm, clip, width=0.33\textwidth]{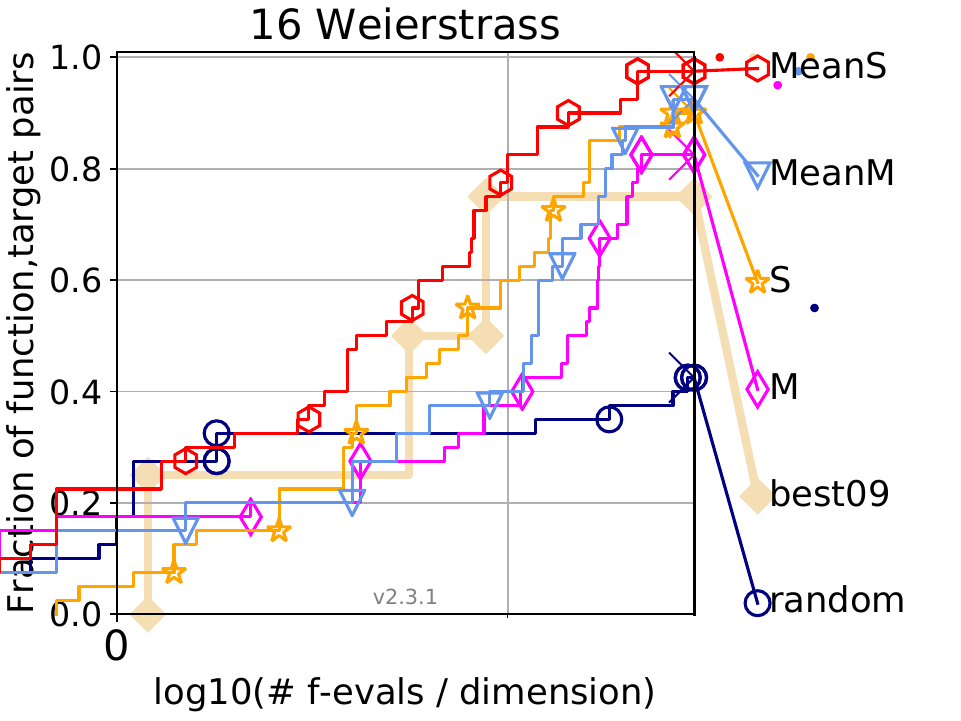} 
\end{tabular}
\caption{
Examples of functions for which the GP mean acquisition is beneficial.
\label{fig-ERTD_gpmean_f1322}
}
\end{figure}

Fig.~\ref{fig-ERTD_gpmean_f1322} shows two functions for which forcing exploitation was beneficial: the mean is useful when optimizing the Sharp Ridge function (f13) or the Weierstrass function (f16) but it is detrimental to a search on Gallagher's Gaussian 21 Peaks (f22).

Overall, EGO is sometimes found to struggle to exploit sufficiently, and using an occasional GP mean acquisition seems to be slightly beneficial, in particular towards late iterations. 
These observations point to the same overall tendency as the ones of another \COCO study \cite{rehbach2020expected} that additional exploitation helps BO searches.
More advanced approaches that were not tested here, e.g. \cite{mcleod2018optimization,de2021greed}, could prove beneficial in this scenario.


\subsubsection{Effect of EI maximization}
\begin{figure}
    \centering
\begin{tabular}{cc@{}c@{}c@{}}
 &{Separable} & {Low conditioning} & {High conditioning} \\
 \rotatebox[origin=l]{90}{\textbf{EI max. tech.}} & 
\includegraphics[trim=0mm 0mm 10mm 0mm, clip, width=0.3\textwidth]{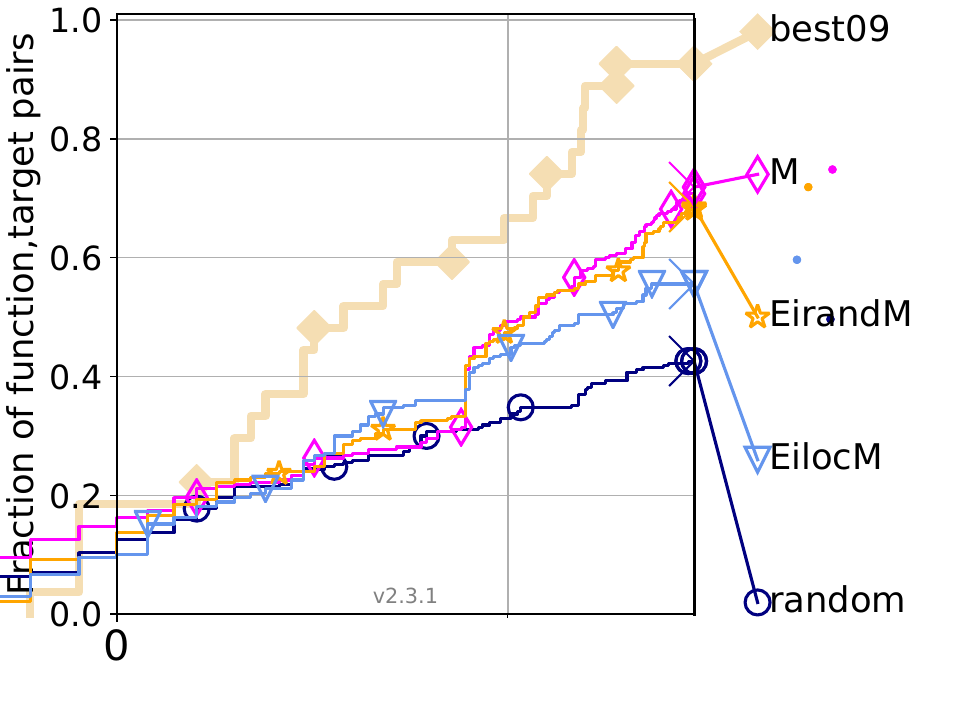} &
\includegraphics[trim=0mm 0mm 10mm 0mm, clip, width=0.3\textwidth]{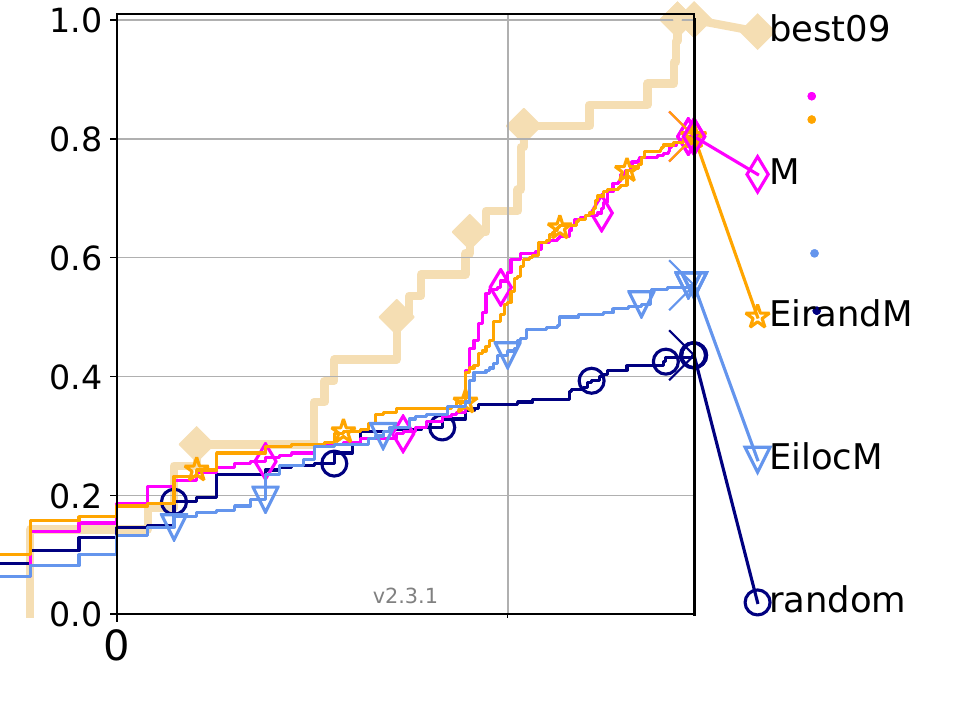}&
\includegraphics[trim=0mm 0mm 10mm 0mm, clip, width=0.3\textwidth]{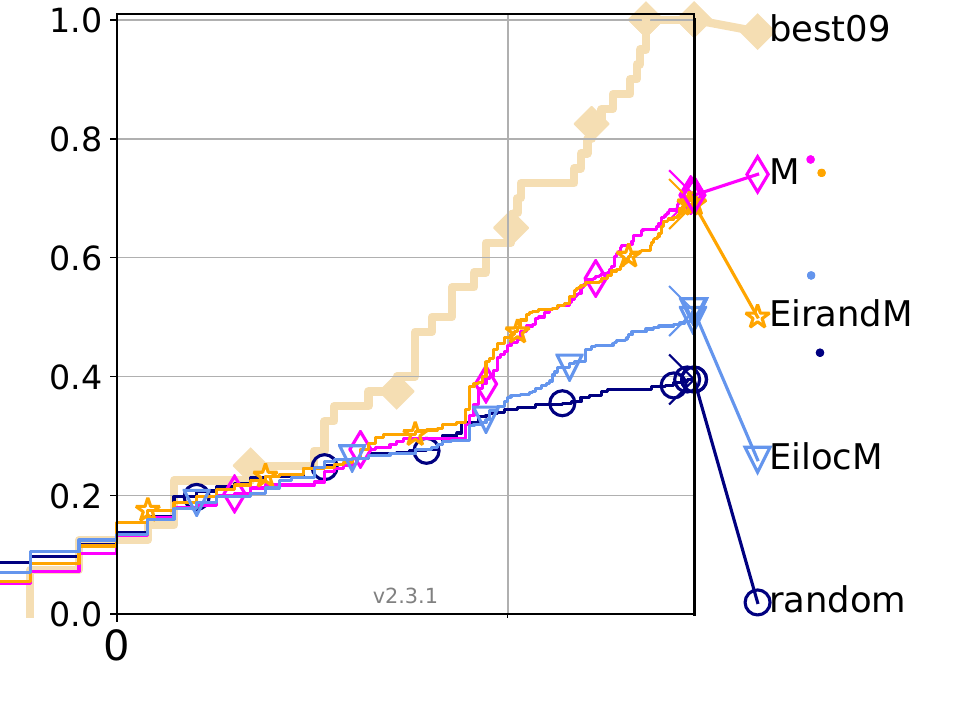} \\
 & {Multimod., strong struct.} & {Multimod., weak struct.} & ~ \\
&
\includegraphics[trim=0mm 0mm 10mm 0mm, clip, width=0.3\textwidth]{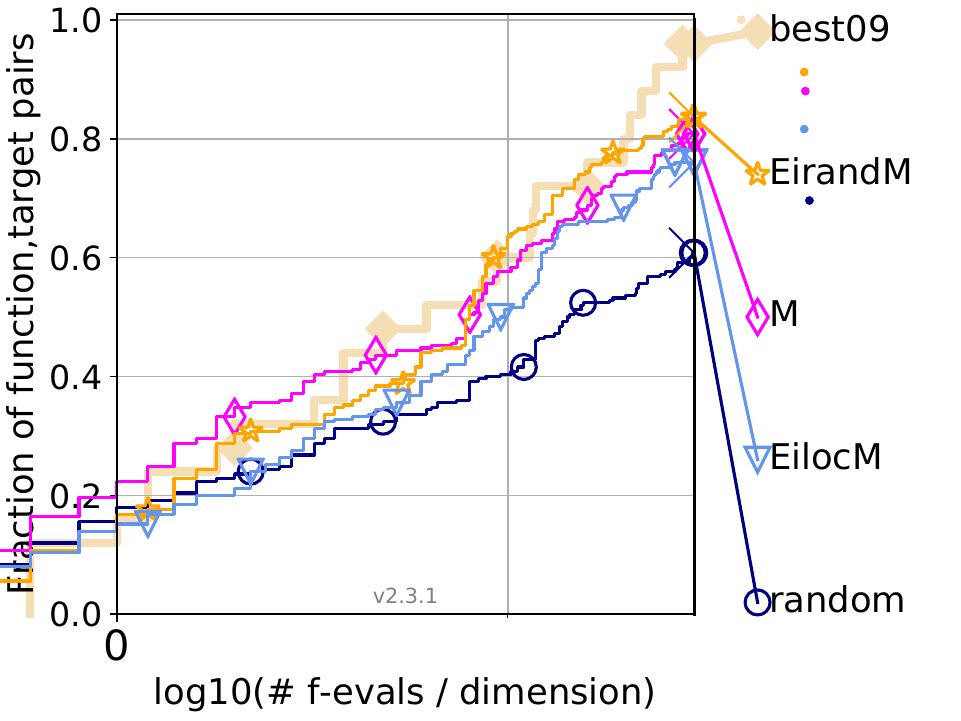}&
\includegraphics[trim=0mm 0mm 10mm 0mm, clip, width=0.3\textwidth]{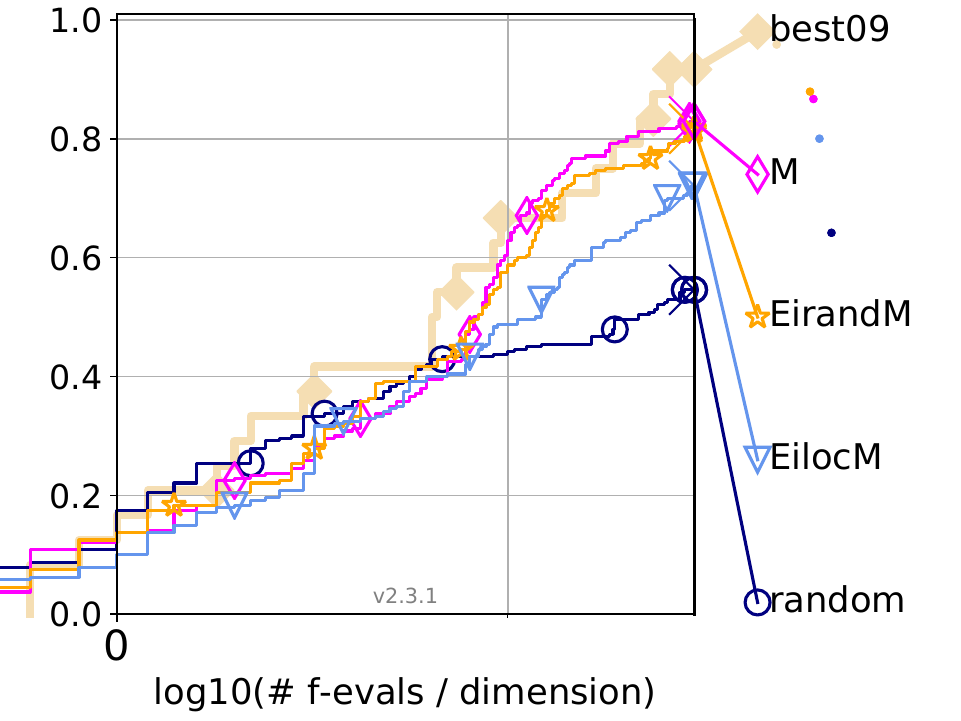}&
~
\end{tabular}
\caption{
ERTD by function groups in $d=5$ dimensions when varying the EI maximization technique: Eirand is a random search, Eiloc is a local search (BFGS), the default (M) is a restarted local search. All versions use a medium size initial DoE. 
\label{fig-ERTD_eiopt_05D_groups}
}
\end{figure}

Figure~\ref{fig-ERTD_eiopt_05D_groups} confirms that globally optimizing EI is much more important than locally optimizing it, regardless of the objective function topology. 
The effect is strongest on unimodal functions (which was not expected as EI is less multimodal on those functions, but EI remains multimodal in all cases).
For $d=5$ on multimodal functions with strong structure, Eirand slightly outperforms our baseline, showing that locally optimizing EI is not needed at all (not that both approaches lead to similar results in smaller dimensions, see supplementary material \cite{supplement_BO_COCO}).

It is necessary to optimize EI sufficiently well to benefit from the information captured by the GP.
The single BFGS is in any case a very poor choice because of the multimodality of EI.
The very limited benefit of local optimization with restart compared to a random search in small dimensions can be explained by the fact the GP hyperparameters, which affect the shape of EI, are uncertain, and may change from one iteration to the next. Hence, in small dimension, random search suffices. 
With $d=10$, random search become too inefficient for the space volume and it makes sense to leverage gradients to find points with high EI values, at the condition of restarting the gradient-based search (here BFGS). 


\subsection{Comparison with other algorithms for expensive functions}

\begin{figure}
\begin{tabular}{c@{}c@{}c@{}}
 \textbf{3D} &\textbf{5D} & \textbf{10D} \\
\includegraphics[width=0.32\textwidth]{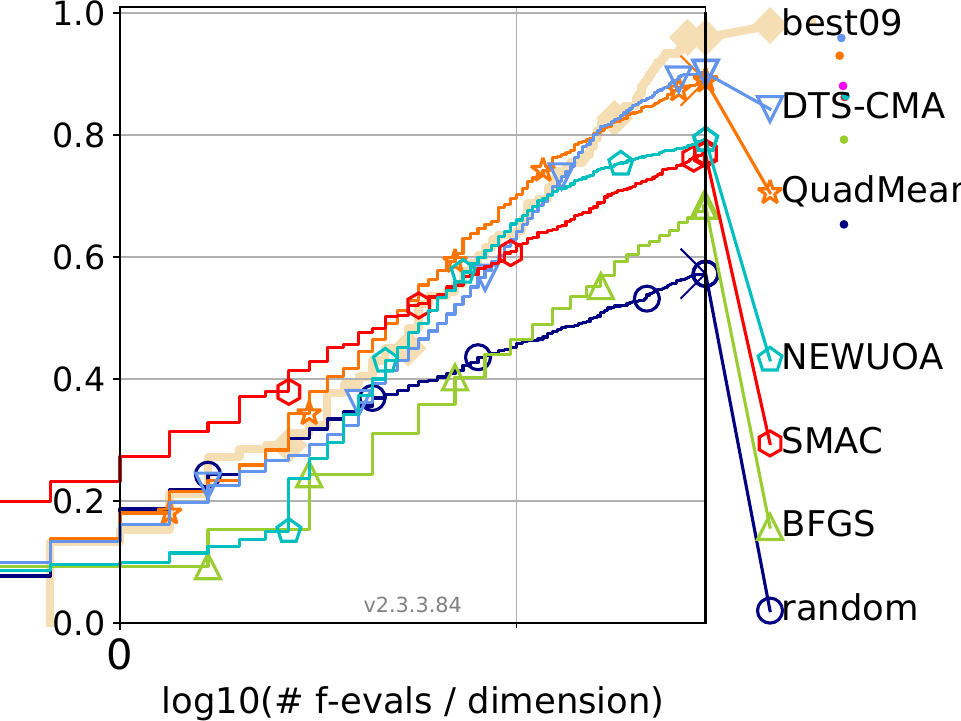}&
\includegraphics[width=0.32\textwidth]{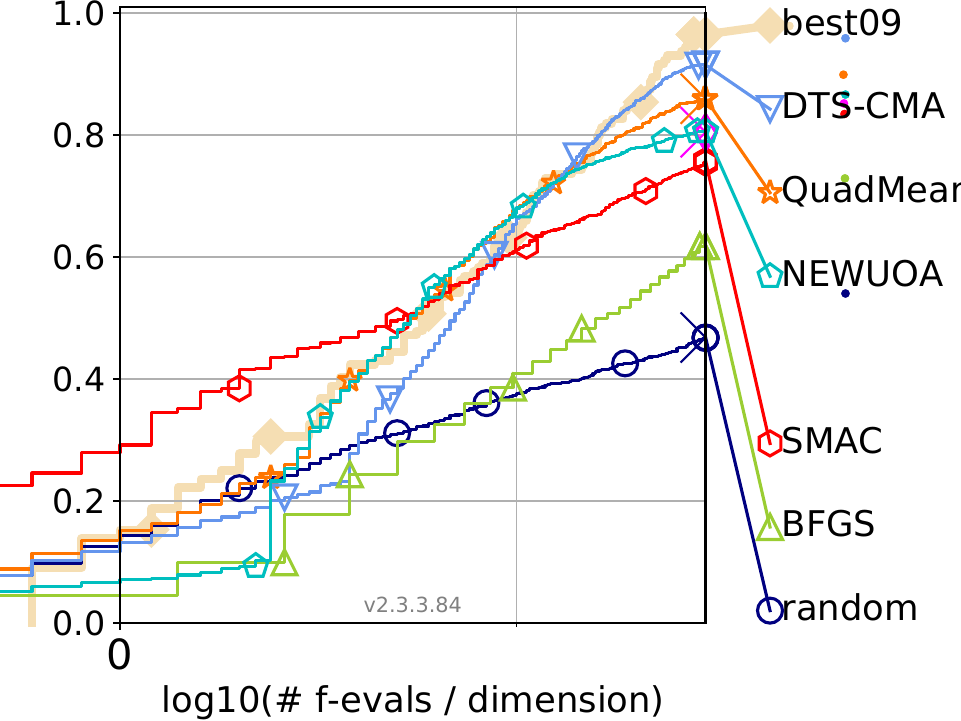}&
\includegraphics[width=0.32\textwidth]{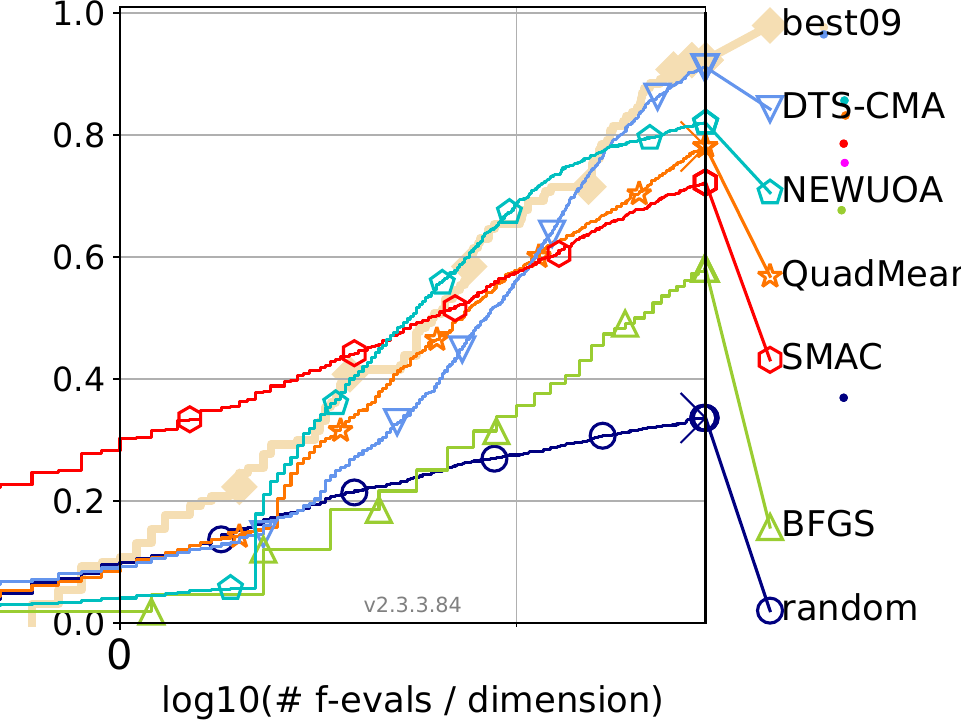}\\
\end{tabular}
\caption{Comparison of the ERTD over all functions in $d=3,5$ and $10$ dimensions of the best EGO variant with the most competitive algorithms on the expensive noiseless COCO testbed.
\label{fig-ERTD_compare}
}
\end{figure}

\begin{figure}
    \centering
\begin{tabular}{c@{}c@{}c@{}}
 {Separable} & {Low conditioning} & {High conditioning} \\
 \includegraphics[trim=0mm 0mm 0mm 0mm, clip, width=0.33\textwidth]{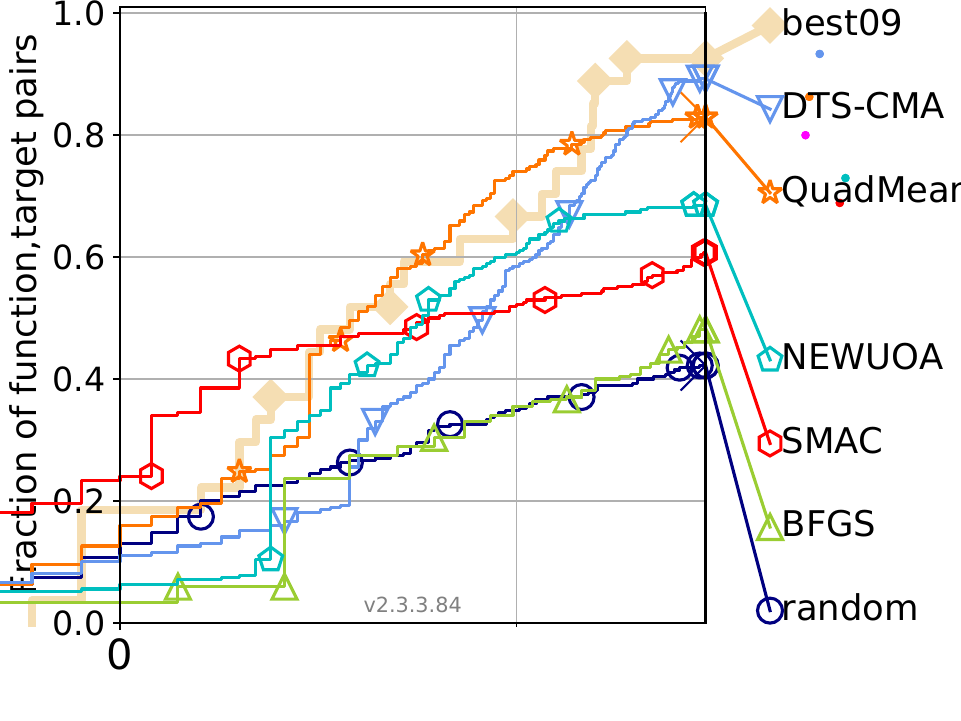} &
 \includegraphics[trim=0mm 0mm 0mm 0mm, clip, width=0.33\textwidth]{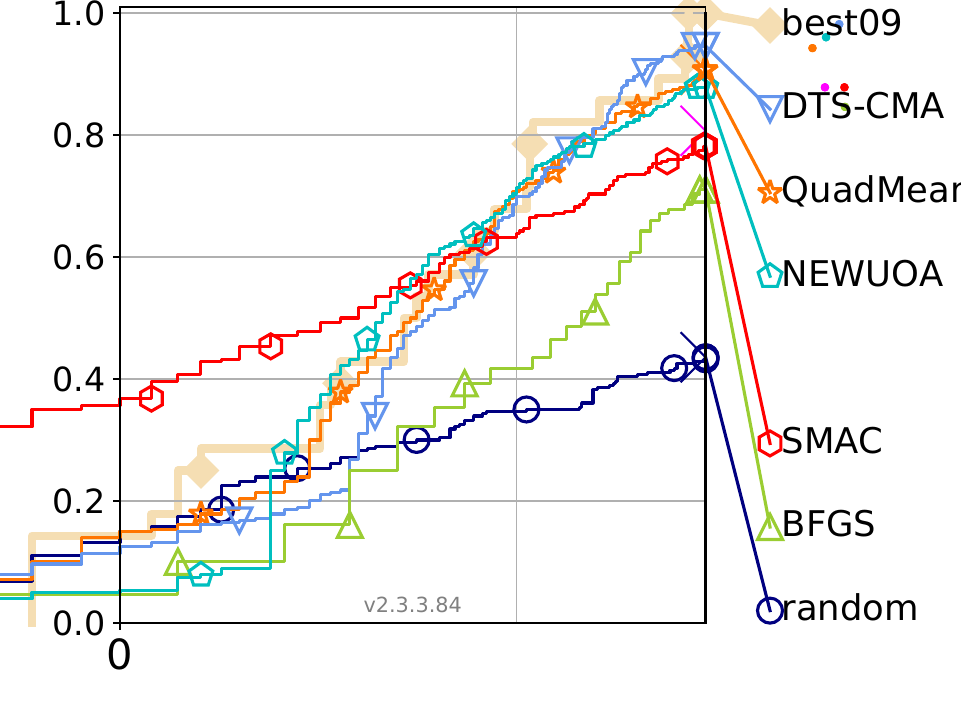}&
 \includegraphics[trim=0mm 0mm 0mm 0mm, clip, width=0.33\textwidth]{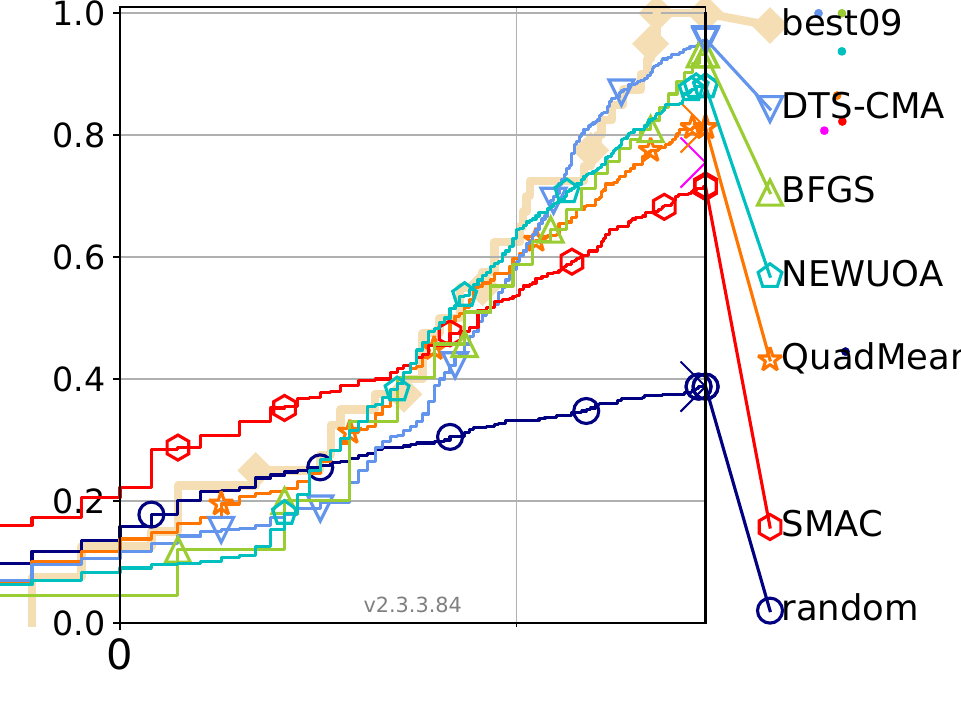} \\
 {Multimod., strong struct.} & {Multimod., weak struct.} & \\
 \includegraphics[trim=0mm 0mm 0mm 0mm, clip, width=0.33\textwidth]{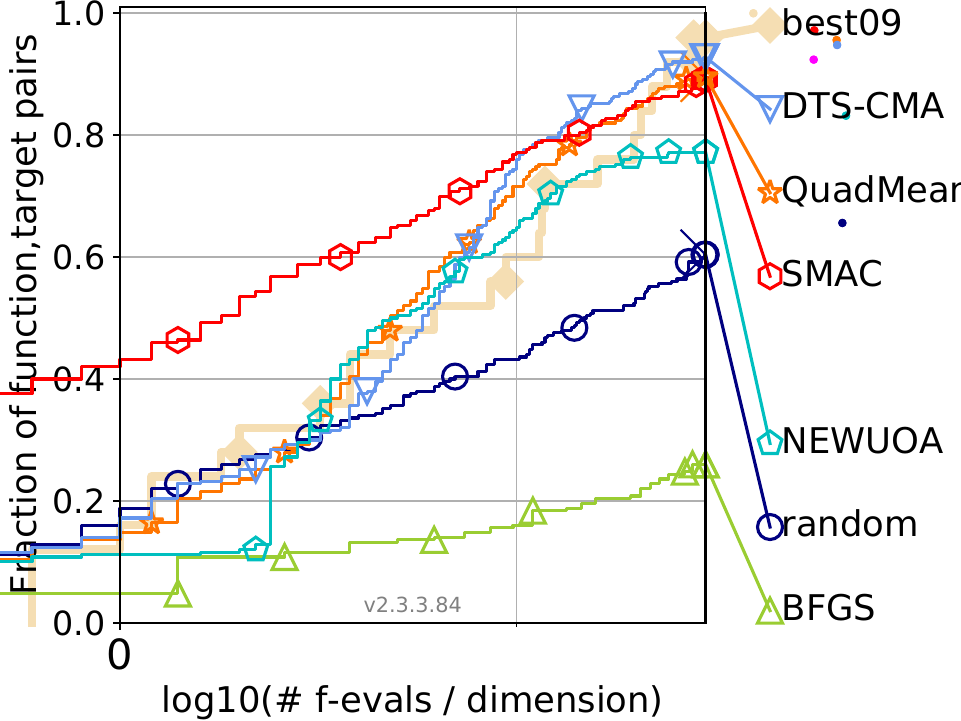}&
 \includegraphics[trim=0mm 0mm 0mm 0mm, clip, width=0.33\textwidth]{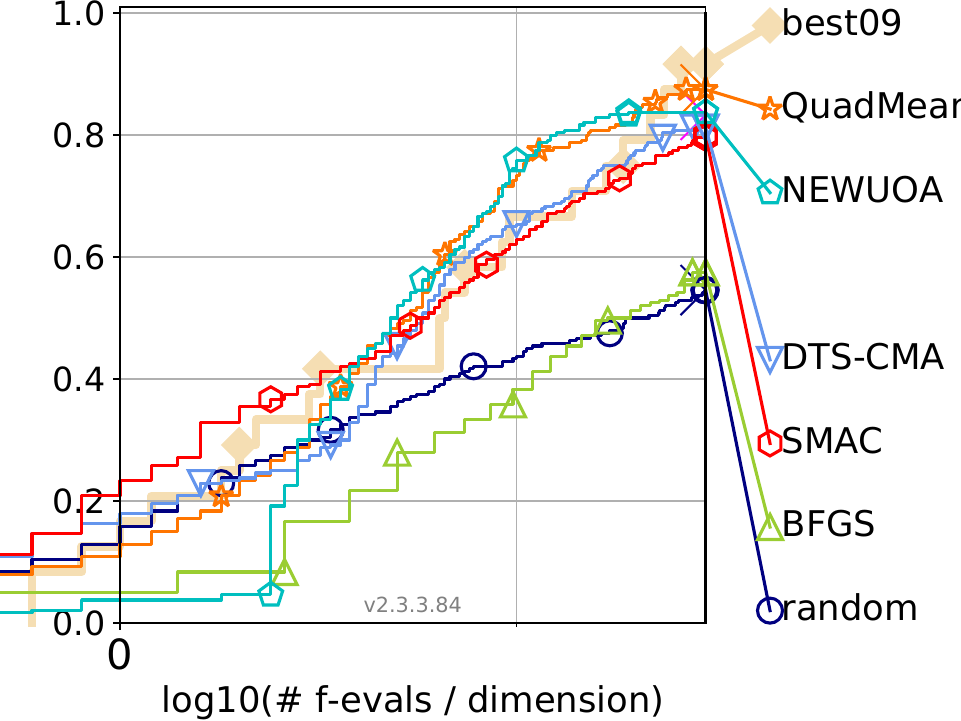}&
\end{tabular}
\caption{
ERTD by function groups in $d=5$ dimensions of the best EGO variant with the most competitive algorithms on the expensive noiseless COCO testbed.
\label{fig-ERTD_compare_05D_groups}
}
\end{figure}

The tests, reported in Figs.~\ref{fig-ERTD_compare} and \ref{fig-ERTD_compare_05D_groups}, compare our best classical EGO implementation to three state-of-the-art competitors for expensive functions: NEWUO, SMAC and DTA-CMA-ES. We use best09, L-BFGS and random search as additional baselines.

DTS-CMA ~\cite{pitra2017comparison} is a surrogate-assisted evolution strategy based on a combination of the CMA-ES algorithm and Gaussian process surrogates.
The DTS-CMA solver is known to be very competitive compared to the state-of-the-art black-box optimization solvers particularly on some classes of multimodal test problems. 
SMAC~\cite{hutter2013evaluation} is another GP-based BO solver; the main differences between SMAC and our approach are the use of an isotropic covariance function,
and the ability to start with a DoE that consists of a single point. SMAC is also performing very well compared to the state-of-the-art blackbox optimizers. 
NEWUOA~\cite{powell2006newuoa} is a trust-region based algorithm that uses interpolating quadratic response surfaces, that showed very good performance, in particular in high dimension~\cite{hansen2010comparing}. 
Based on our previous findings, our EGO algorithm uses a small initial DoE, a Matern52 kernel, a quadratic trend, and the occasional GP mean acquisition function. It is reported as the QuadMean algoritm in the figures.

On average over all functions (Fig.~\ref{fig-ERTD_compare}), we see that EGO (i.e., QuadMean) is competitive for $d=3$ and $5$, but falls behind DTS-CMA-ES and NEWUOA for $d=10$.
Regardless of dimension, SMAC has a very early start and already optimizes the functions while the other methods are still initializing. 
However the rate of convergence of SMAC is slower and this method is outperformed by the other competitors for budgets larger than $10 \times d$.
DTA-CMA-ES makes the better use of large budgets and outperforms all other methods for budgets larger than $20 \times d$. 
The areas of expertise of EGO and NEWUOA are somewhere in the middle, with best performance for intermediate budgets, EGO being best in small dimension and NEWUOA in larger dimension. 
EGO proposes an attractive trade-off, as it is never the worst algorithm.

One may notice that all algorithms are quite competitive and often beat the best09 reference, largely outperforming a baseline algorithm such as L-BFGS.

Looking now at performance by function groups (Fig.~\ref{fig-ERTD_compare_05D_groups}), we see that EGO excels on separable functions, clearly surpassing the other algorithms (except DTS-CMA for the latest iterations and best09, the 2009 upperbound which contains algorithms adapted to separable functions). 
EGO and NEWUOA also significantly outperform the other algorithms on group 5 (multimodal functions with weak global structure). 
EGO underperforms comparatively on the third group (unimodal functions with high conditioning). 



\begin{figure}
    \centering
\begin{tabular}{c@{}c@{}c@{}}
 {f7, $d=5$} & {f17, $d=5$} & {f21, $d=5$} \\
 \includegraphics[trim=0mm 0mm 0mm 8mm, clip, width=0.33\textwidth]{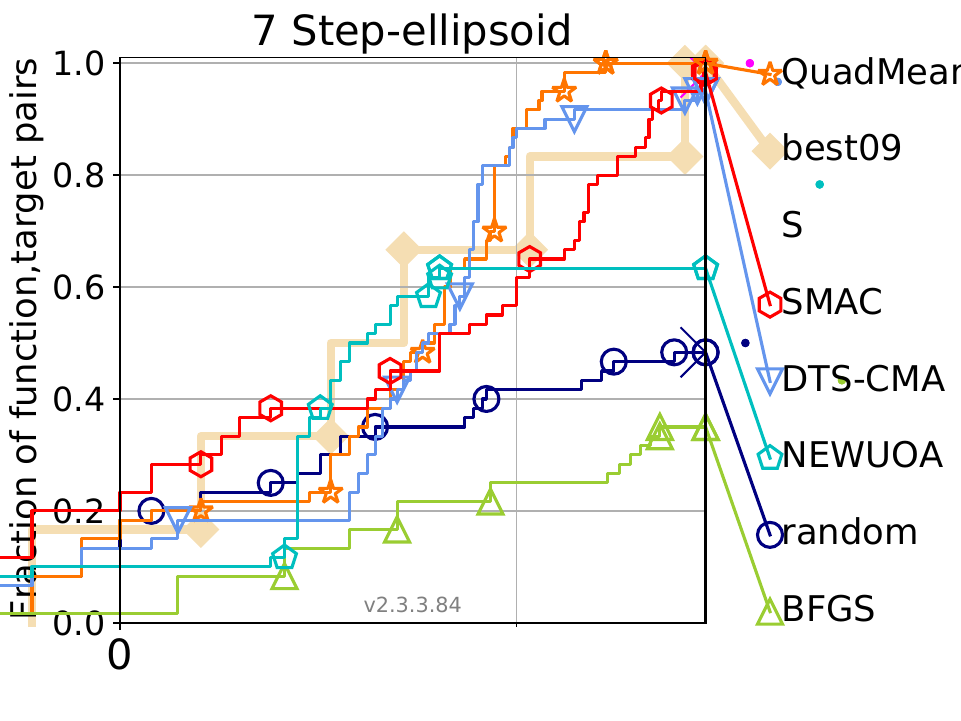} &
 \includegraphics[trim=0mm 0mm 0mm 8mm, clip, width=0.33\textwidth]{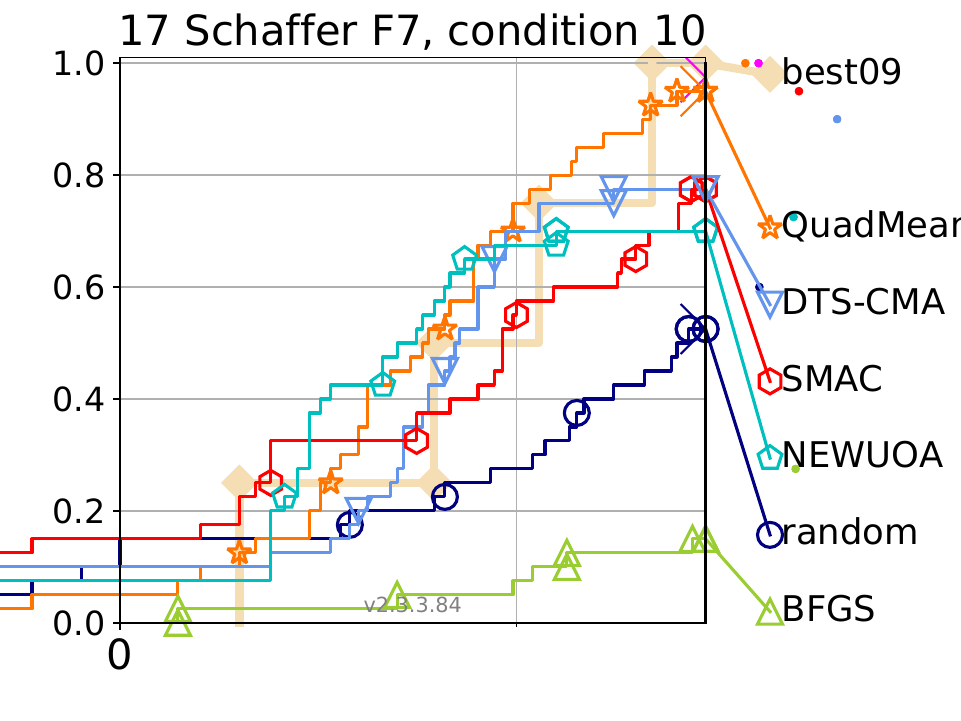}&
 \includegraphics[trim=0mm 0mm 0mm 8mm, clip, width=0.33\textwidth]{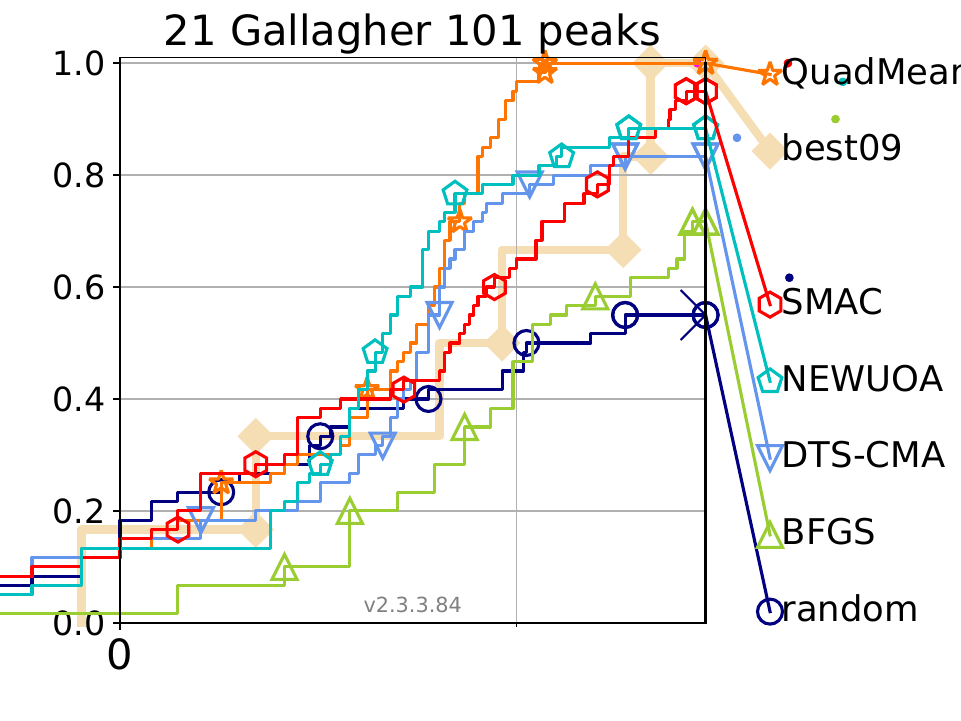} \\
 {f12, $d=5$} & {f11, $d=10$} & {f23, $d=10$} \\
 \includegraphics[trim=0mm 0mm 0mm 8mm, clip, width=0.33\textwidth]{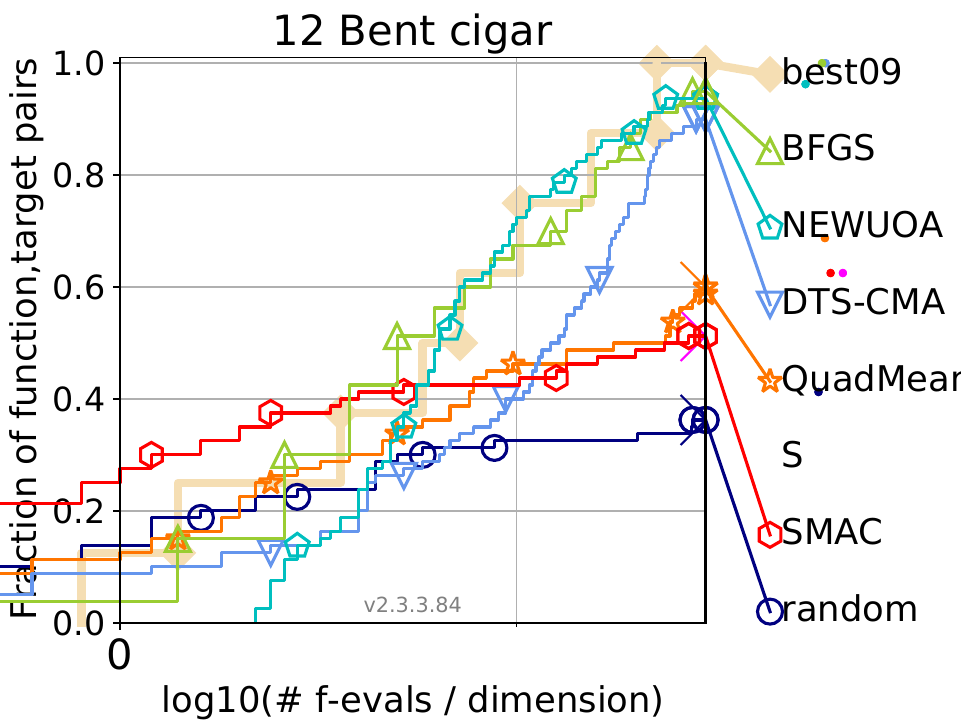}&
 \includegraphics[trim=0mm 0mm 0mm 8mm, clip, width=0.33\textwidth]{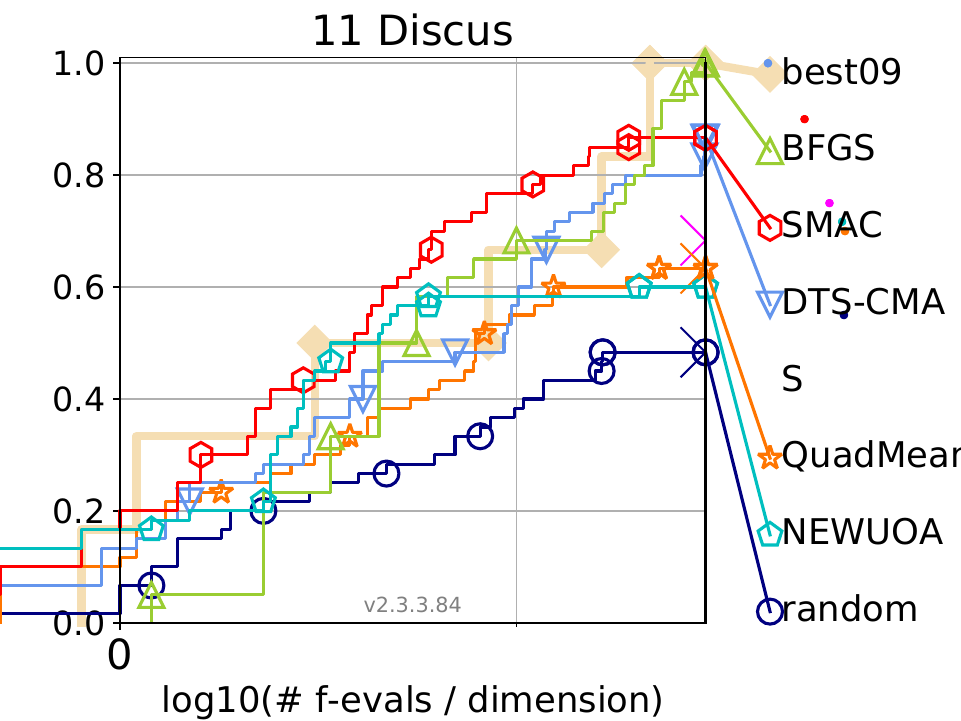}&
 \includegraphics[trim=0mm 0mm 0mm 8mm, clip, width=0.33\textwidth]{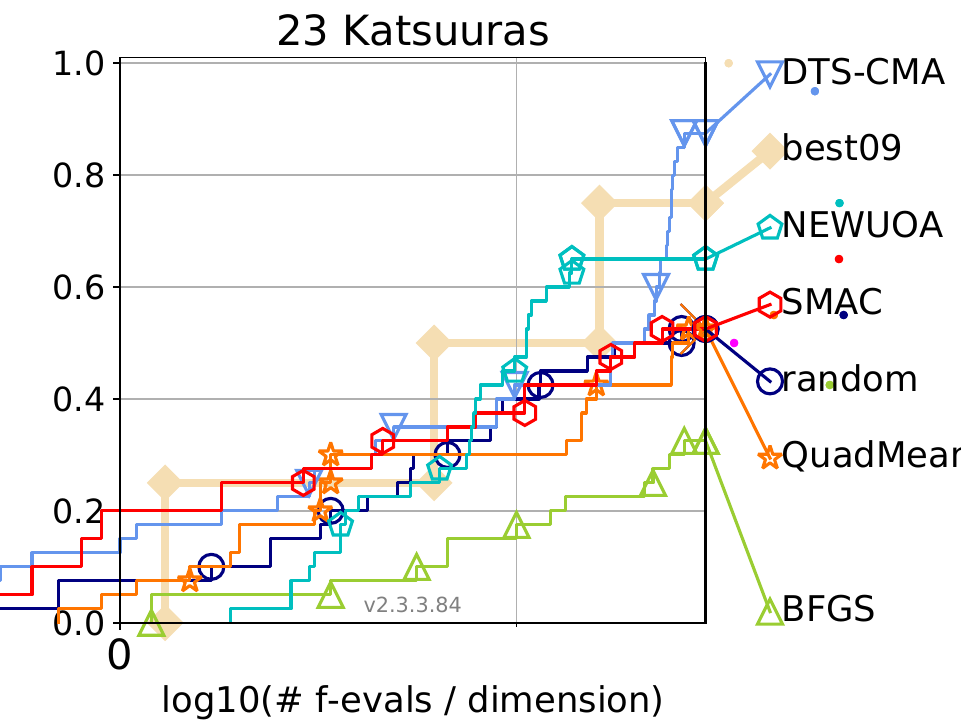}
\end{tabular}
\caption{
Examples of functions for which EGO performs well, top row, or badly, bottow row, when compared to best09 and other state-of-the-art optimizers.
\label{fig-ERTD_functions_3good_3bad}
}
\end{figure}

Figure~\ref{fig-ERTD_functions_3good_3bad} shows examples of individual functions 
which are particularly advantageous (top row) or detrimental (bottom row) to EGO relatively to the other competitors.
f7, f17 and f21 in $d=5$ dimension make the favorable set. They all have a low conditioning. The functions which are most difficult for EGO (f12 in $d=5$, f11 and f23 in $d=10$) have a high conditioning or a high dimension. In addition, as it is explained in Appendix~\ref{sec-GP_Q2_optim}, f23, Katsuura's function, is hard to model with a GP (left of Figure~\ref{fig-Q2}) and even more so with a quadratic trend (right of the Figure).

\section{Concluding comments}
This article provides a comprehensive experimental validation of well and less-known 
design considerations of Bayesian optimization algorithms. 

Our results validate most common practices: the Mat\'ern 5/2 kernel, no data transformation and no non-stationary GPs by default, and a small initial DoE. 
The only slightly disruptive result we observed is that a quadratic trend helps and should be used as default. 
We have uncovered the importance of the EI internal optimization, for which global search is much more important than local one.

When compared to other algorithms, EGO excels in small dimension ($d \le 5$) with moderate budgets (10 to 20 times $d$). It is best where we expect it to be, that is for multimodal functions, and not so good on poorly conditioned unimodal functions. 
EGO shows surprisingly good performance on some convex functions. 
Overall, EGO is a fast and reliable optimization method for small budgets and dimensions in the sense that it is never the worst among the strong competitors we have considered.
However, EGO is particularly sensitive to the curse of dimensionality. 

This work indicates the following directions for improving BO algorithms:
\begin{itemize}
	\item BO algorithms should better comply with poorly conditioned unimodal functions.
	\item They would gain from mimicking the early start of the SMAC method. 
	\item Transformations (scaling, warping) offers occasional gains but were detrimental on average: we believe that work is still needed to provide fast, automated and reliable methods. 
	\item BO algorithms should better scale with the dimension $d$. This is already a popular topic of research, with attempts to identify latent spaces of reduced dimension (e.g., linear in \cite{wang2016bayesian}, nonlinear in \cite{oh2018bock,gaudrie2020modeling}), via ad-hoc structure \cite{kandasamy2015high} or with an additional control on the search volume through a trust region (e.g., \cite{Turbo_2019,TREGO_2020}).
\end{itemize}

In any case, an extensive testing of the algorithm, such as what has been done here, is 
the main track to validate new ideas. This is difficult with BO algorithms which are computationally intensive. There is therefore also a need for BO methods which would either be more frugal or make a better use of hardware (as is done in \cite{balandat2019botorch}).

\subsubsection*{Acknowledgments}
We would like to warmly thank Dimo Brockhoff for his precious comments and help with the \COCO benchmark software. This study is a late addition to the PGMO/FMJH AESOP (Algorithms for Expensive Simulation-Based Optimization) project (\url{https://www.fondation-hadamard.fr/fr/lmh-projets-collaboratifs/projets-collaboratifs-maths-ia}).


%
%

\bibliographystyle{plain} 
\bibliography{BO_COCO}   

\appendix

\section{From GP representation ability to optimization performance}
\label{sec-GP_Q2_optim}
The ability to build a reliable representation of the function 
is likely to be a critical ingredient of the optimization method.
In this Section, we measure the regression accuracy of GP variants on the functions of the BBOB noiseless testbed and investigate how it can be predictive of an EGO algorithm's performance.

In accordance with the rest of the article (cf. in particular Section \ref{sec-factors} for more explanations), 5 GPs are compared: the \emph{default}, made of Matérn 5/2 kernels and constant trend; the \emph{quadratic} version, which is the default with the trend made quadratic; the \emph{scaling} variant where the inputs are transformed; the \emph{warping} version where the outputs are transformed; and the \emph{exponential} version where the kernel is exponential, hence creating a less regular GP than the other versions.

\subsection*{Regression performance of the different GPs on the functions testbed}

For each of the 24 test functions, a maximin LHS design of size $30 d$ is created, the corresponding function is evaluated there and the investigated variant of the GPs is learned. This GP is then tested on another maximin LHS of the same size. 
The procedure is repeated on 15 instances of the functions, i.e., random translations and rotations of the functions.

The regression accuracy is measured by two criteria. First, the Q2 criterion which quantifies the prediction (GP mean) accuracy on a test set as a normalized quadratic distance to the true function values with values ranging from 1 (perfect) to -1 (opposite variations), 0 being the best possible constant prediction, the function mean. 
Second, a normality test is carried out on the normalized residuals, $(y^{(i)}-\mu(x^{(i)}))/\sigma(x^{(i)})$. Neglecting the spatial correlation of the process, this measures the validity of the assumption that the normalized test residuals should be $\mathcal N(0,1)$, in other terms that both the predicted mean and variance work well together to satisfy the normality assumption.
Normality of the residuals is measured by the p-value of a Kolmogorov-Smirnov test.
The p-value is the probability of observing a sample at least as extreme as the one observed under the hypothesis of a normal distribution (with proper centering and normalizing). 
Small p-values (typically less than 0.05) yield to rejection of the normality assumption.
Because there are 15 instances of each function, the mean and standard deviation of the Q2 and p-values are calculated. In both cases, the larger (closer to 1), the better.

\begin{figure}
\mbox{
\begin{minipage}{0.5\textwidth}
\centering
\includegraphics[width=\textwidth]{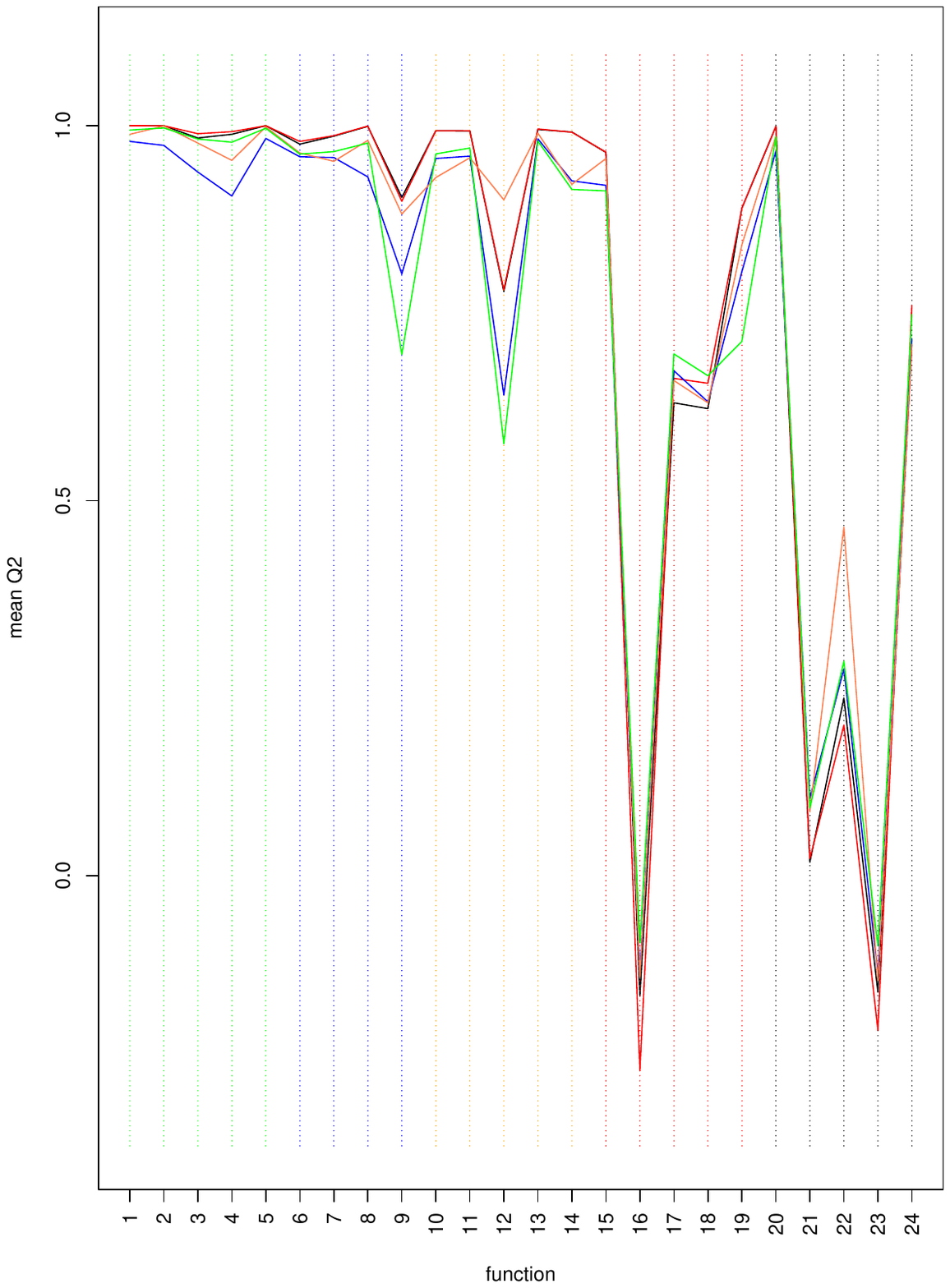} \\
a) mean
\end{minipage}
\begin{minipage}{0.5\textwidth}
\centering
\includegraphics[width=\textwidth]{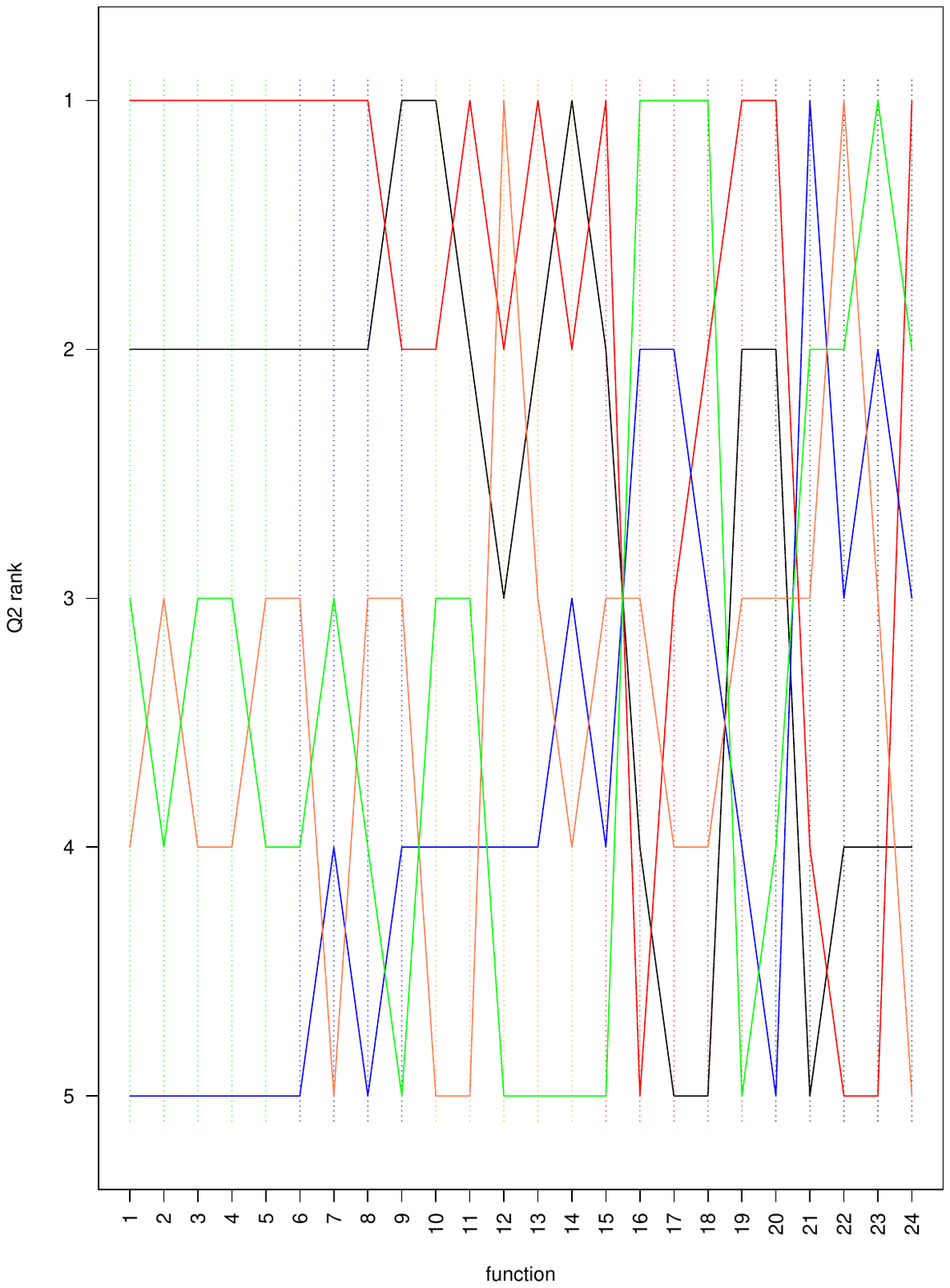} \\
b) ranked by mean Q2
\end{minipage}
}
\caption{
Q2 statistics for the 5 GP versions on the BBOB noiseless functions: 
a) Q2 means b) ranks of the GPs in terms of mean Q2 (1 is best, 5 worst).
Black = default GP, red = quadratic trend, green = exponential kernel, blue = with scaling, coral = with warping.
The colours of the vertical dotted lines reflect the function group.
\label{fig-Q2}
}
\end{figure}

Q2 statistics are plotted in Fig.~\ref{fig-Q2}: 
Fig.~\ref{fig-Q2} a) is the mean Q2 achieved by a specific GP per function, 
b) is the rank of a GP in terms of Q2 mean (1 is best, 5 worst). It describes how variants perform relatively to each other.
The colours of the vertical dotted lines designate the function group: green for separable functions, blue for unimodal functions with low or moderate conditioning, orange for unimodal functions with high conditioning, red for multimodal funtions with adequate global structure and black with weak global structure.
The Q2 standard deviation were not found to provide additional information with respect to the mean: the Q2 standard deviation varies in the opposite direction to the Q2 mean; large variations occur for functions with low Q2 means and the best performing GPs in terms of Q2 mean have lowest standard deviation.

It is observed that all GPs have Q2's that vary together and have roughly the same order of magnitude for the different functions. 
The largest differences occur for f12 (Bent Cigar) and f22 (Gallagher's Gaussian 21-hi Peaks) with Q2 interval of 0.3 but on most functions the spread is smaller than 0.1. 
A general pattern for all GPs is that the functions 15 to 24, which belong to the multimodal families, are more difficult to predict (with the exception of Rastrigin f15 and Schwefel f20). Vice versa, the separable and unimodal functions f1 to f14 are well predicted apart from the Rotated Rosenbrock f9 and the Bent Cigar f12. 
When the functions are well predicted (again, unimodal or separable plus f15 and f20), the best GPs are the quadratic or the default. This is best seen in Fig.~\ref{fig-Q2} b) where the quadratic (red) and default (black) GPs take up all of ranks 1 and 2. 
When the functions are difficult to predict (f12, f16 to f24 but f20), the exponential and quadratic GPs are often good choices. 
The appropriateness of the exponential kernel for modeling the more irregular functions was expected as it is the only kernel with non differentiable trajectories in our methods. The good performance of the quadratic trend is more surprising and may reveal a small bias of the 
BBOB noiseless testbed in favor of functions with a quadratic global structure.
The scaling and warping GPs are rarely among the best modeling options. Out of the 24 functions, scaling and warping GPs have ranks 1 or 2 only 4 and 2 times, respectively, against 19, 16 and 7 times for the quadratic, the default and the exponential GPs. 
Yet, as a possible illustration to the No Free Lunch theorem in machine learning \cite{wolpert1996lack}, the GP with warping is the best for f12 and f22 (the functions with least consensus among GPs) and scaling the inputs yields the best GP for Gallagher's Gaussian 101-me Peaks (f21).

\subsection*{Link between GP regression quality and EGO performance}
We now discuss how the regression performance of the GP is correlated to the optimization performance. 
It is logical that a metamodel (a GP here) with a good regression performance achieves a good optimization performance in a metamodel based optimizaton algorithm. 
However, this is not a necessary condition. A metamodel which correctly predicts the lowest regions of the function up to any rank preserving transformation of the function will lead the optimization algorithm to the minimum (at the condition that the algorithm converges of course).
In addition, optimization is a dynamic, iterative, procedure. A metamodel with too low a complexity to fit well the function at the final budget may still be a good guide earlier in the search because it is fast to learn. 
For these reasons, a GP with a bad Q2 may still work well for optimization. 
The relationship between Q2 and optimization performances is complex. 

As an additional metric, we introduce a scalar for an optimization run called \emph{relative optimization performance}. It is defined as
\begin{equation} 
\Popt(\text{algo}) ~ \coloneqq~ \frac{ \text{ERTD}(\text{algo} @ 30d) - \text{ERTD}(\text{random} @ 30d)} { \text{ERTD}(\text{best09} @ 30d) - \text{ERTD}(\text{random} @ 30d)}  
\label{eq-Popt}
\end{equation} 
Above, $\text{ERTD}(\text{algo} @ 30d)$ is the ERTD of the ``algo'' optimizer after $30d$ function evaluations in dimension $d=5$. 
The denominator is the difference in (ERTD) optimization performance between best09 and a random optimization, which measures function easiness: a small value of the denominator means that a sophisticated algorithm does hardly better than a random search and vice versa.
$\Popt$ is a normalized measure. Values above 0 do better than a random sampling and vice versa. Values above 1 do better than best09. 
$\Popt$ is a performance relative to a random search and best09.

In our experiments, the mean Q2s and p-values achieved with a DoE of size $30 d$ are compared to the ERTD of EGO at the same budget of $30 d$ with the same GP. More precisely, the ERTD is normalized into the relative optimization performance $\Popt$ defined in Eq.~\ref{eq-Popt}.

\begin{figure}
\mbox{
\begin{minipage}{0.5\textwidth}
\centering
\includegraphics[width=\textwidth]{q2_means_lines.pdf} \\
a) regression (Q2)
\end{minipage}
\hfill
\begin{minipage}{0.5\textwidth}
\centering
\includegraphics[width=\textwidth]{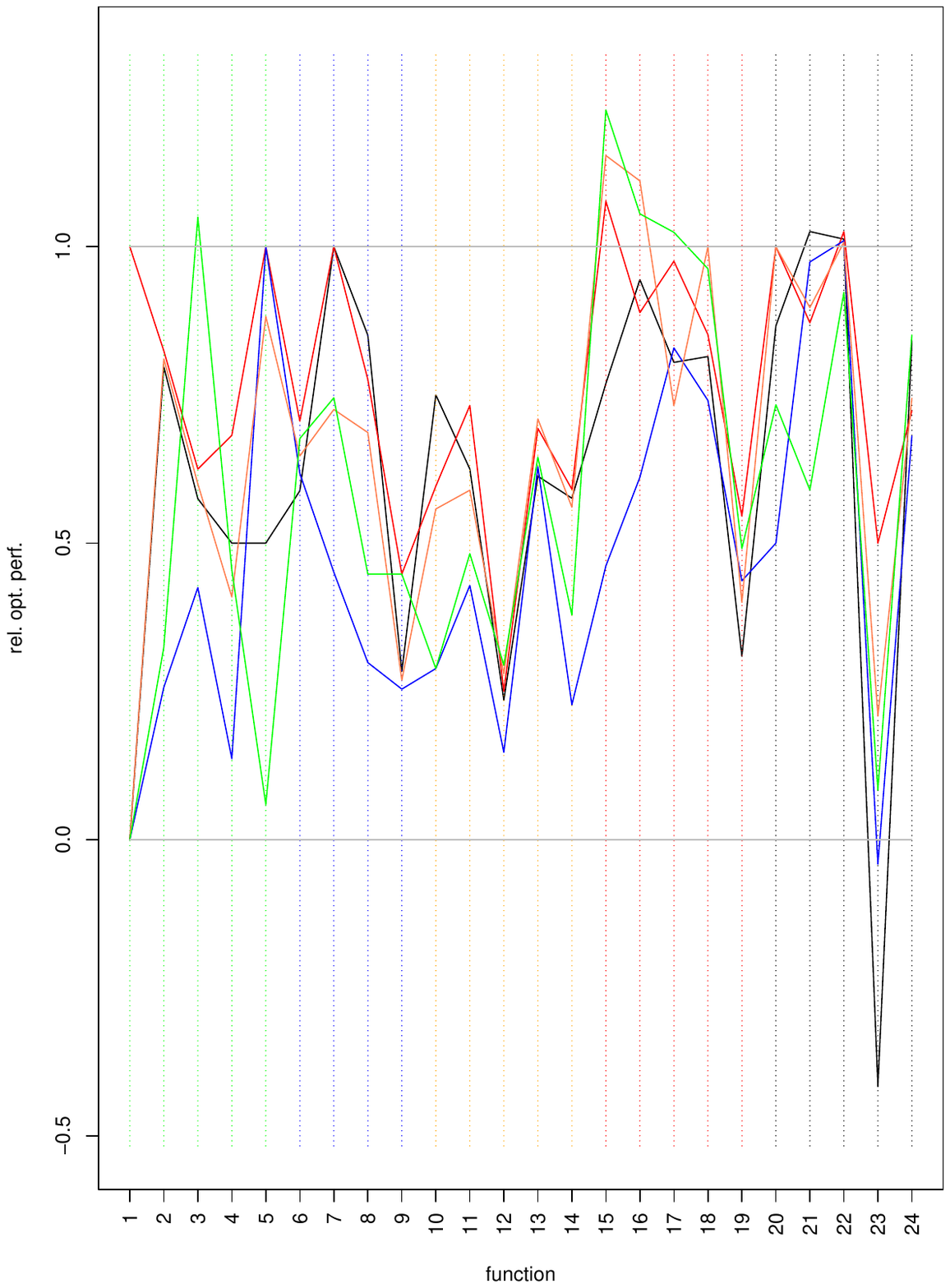} \\
b) rel. optim. perf. ($\Popt$)
\end{minipage}
}
\caption{
Comparison of regression (left) and optimization (right) performances of the GP variants expressed in terms of a) Q2 , b) relative optimization performance, $\Popt$.
Black = default GP, red = quadratic trend, green = exponential kernel, blue = with scaling, coral = with warping.
\label{fig-Q2optim}
}
\end{figure}

Two types of links are studied. 1) the link between the magnitude of Q2 and that of $\Popt$, cf. Fig.~\ref{fig-Q2optim}.
2) the link between relative sensitivities of the variants. More precisely, we investigate if an ordering of the Q2s is indicative of an ordering of the $\Popt$.

In Fig.~\ref{fig-Q2optim}, it is seen that bad values of Q2 relatively to the other functions 
in the same group are often indicative of a bad performance later on in optimization, for all EGO variants. This is the case for functions f9 (Rotated Rosenbrock), f12 (Bent Cigar) and f23 (Katsuura). When this occurs, the GP is not able to sufficiently learn the functions to guide an optimization search. 
The Q2s of the GPs on f16 and f21 are also low relative to their groups, however some of the associated EGO variants have high $\Popt$. In these cases, although the functions are not accuratly predicted, some of the EGO algorithms perform well relatively to best09 and random searches. These high optimization performances are artifacts of the definition of $\Popt$ which is relative to random and best09.

\begin{figure}
\mbox{
\begin{minipage}{0.5\textwidth}
\centering
\includegraphics[width=\textwidth]{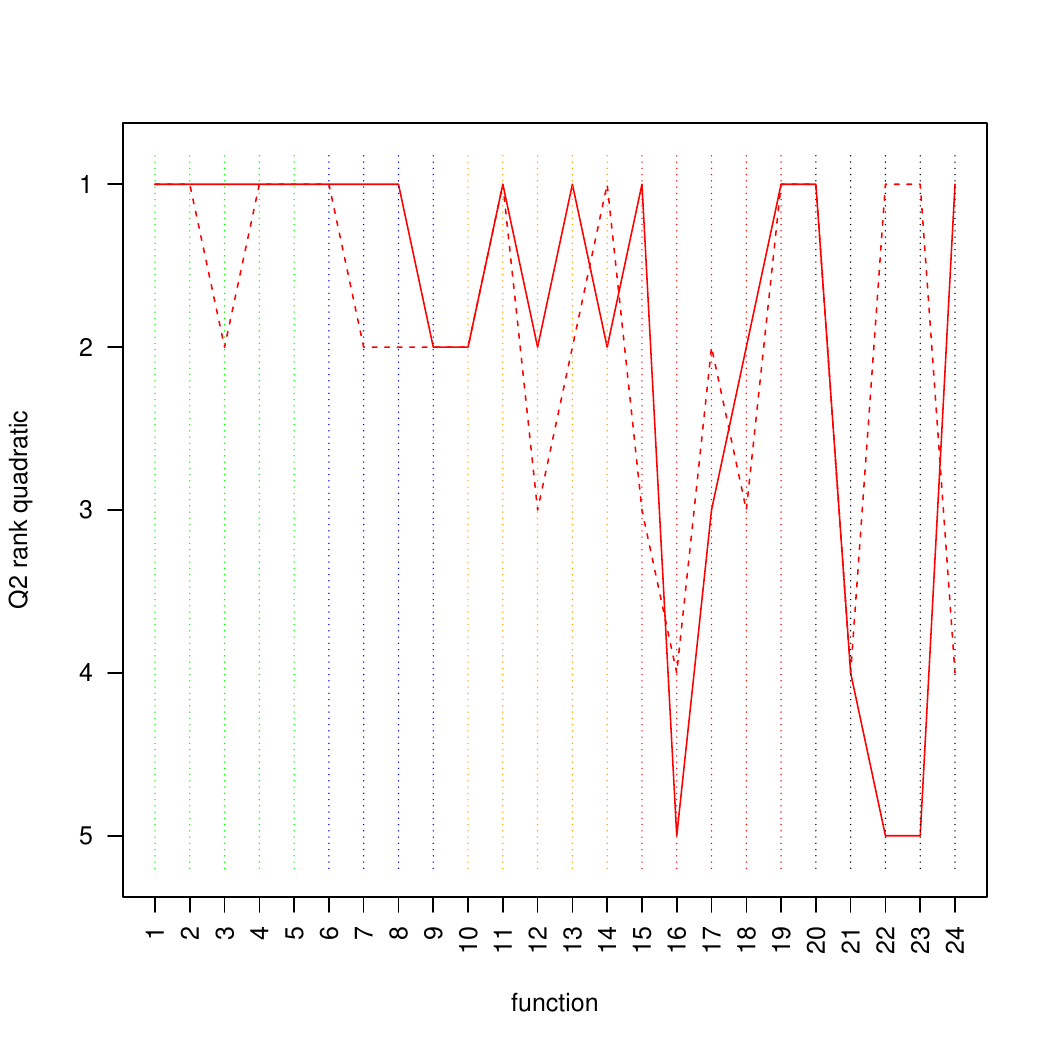} \\
a) quadratic
\end{minipage}
\hfill
\begin{minipage}{0.5\textwidth}
\centering
\includegraphics[width=\textwidth]{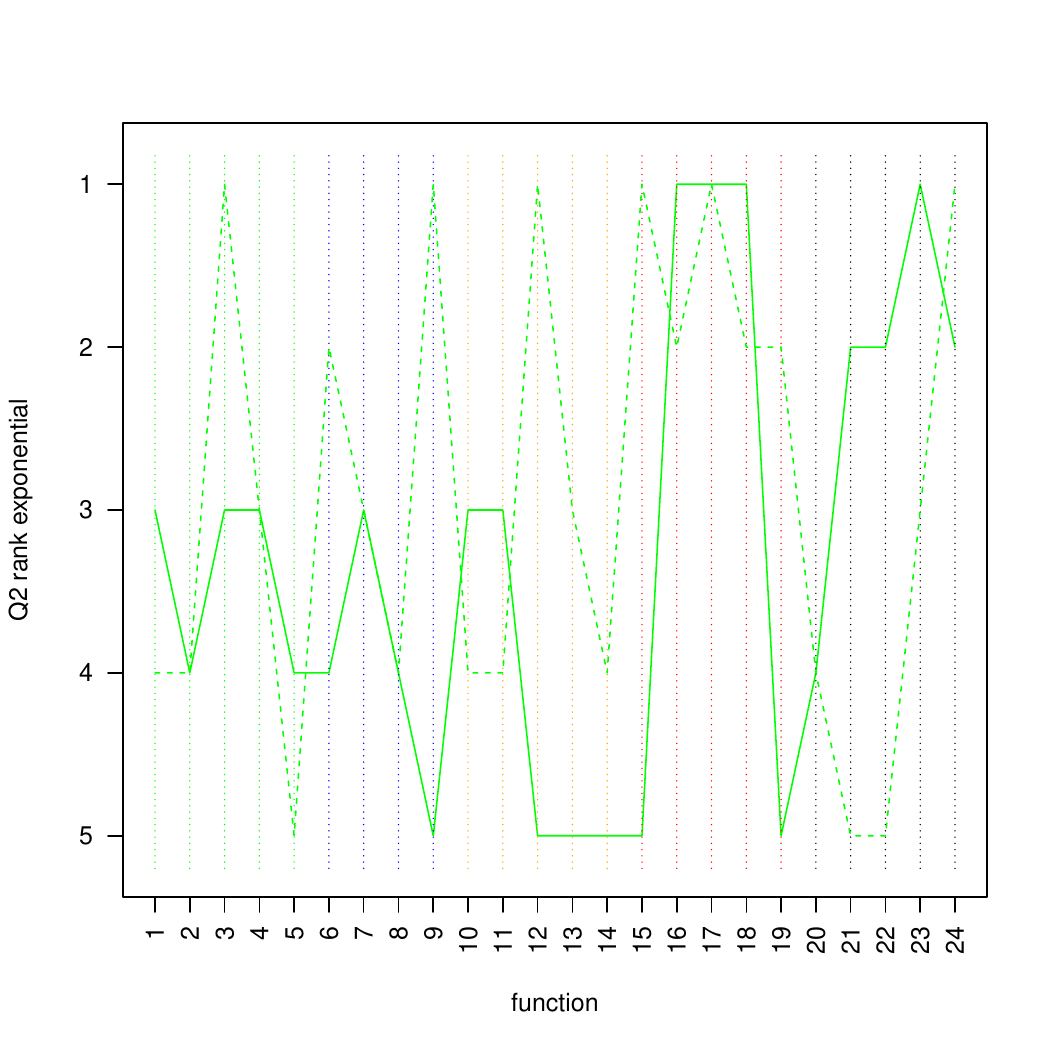} \\
b) exponential
\end{minipage}
}
\caption{
Comparison of regression and optimization performances of two GP variants expressed in terms of ranks (1 is best, 5 worst).
Solid lines are Q2 ranks, dashed lines final optimization (ERTD) ranks. a) quadratic (red), b) exponential (green) GPs.
\label{fig-Q2optimRanks}
}
\end{figure}

In Fig.~\ref{fig-Q2optimRanks}, all performance measures are normalized in terms of ranks ranging from 1 (best) to 5 (worst).
These plots allow to compare the Q2 means (solid lines) and the optimization performances as ERTD (dashed lines) of a GP/EGO variant relatively to the other GPs/EGOs.
A distance of less than 1 rank is a good agreement. If the dashed line is above the continuous line, the optimization performance was better than expected from the Q2 performance, and vice versa.
Excepted for f15, f22, f23 and f24, there is a good agreement between the Q2 and the ERTD ranks of the quadratic GP. 
The exponential kernel's optimization rank is often underestimated when inferred from its Q2 rank: out of the 11 functions where the ranks differ by more than 1, 7 are underestimations.

\begin{figure}
\centering
\includegraphics[width=0.5\textwidth]{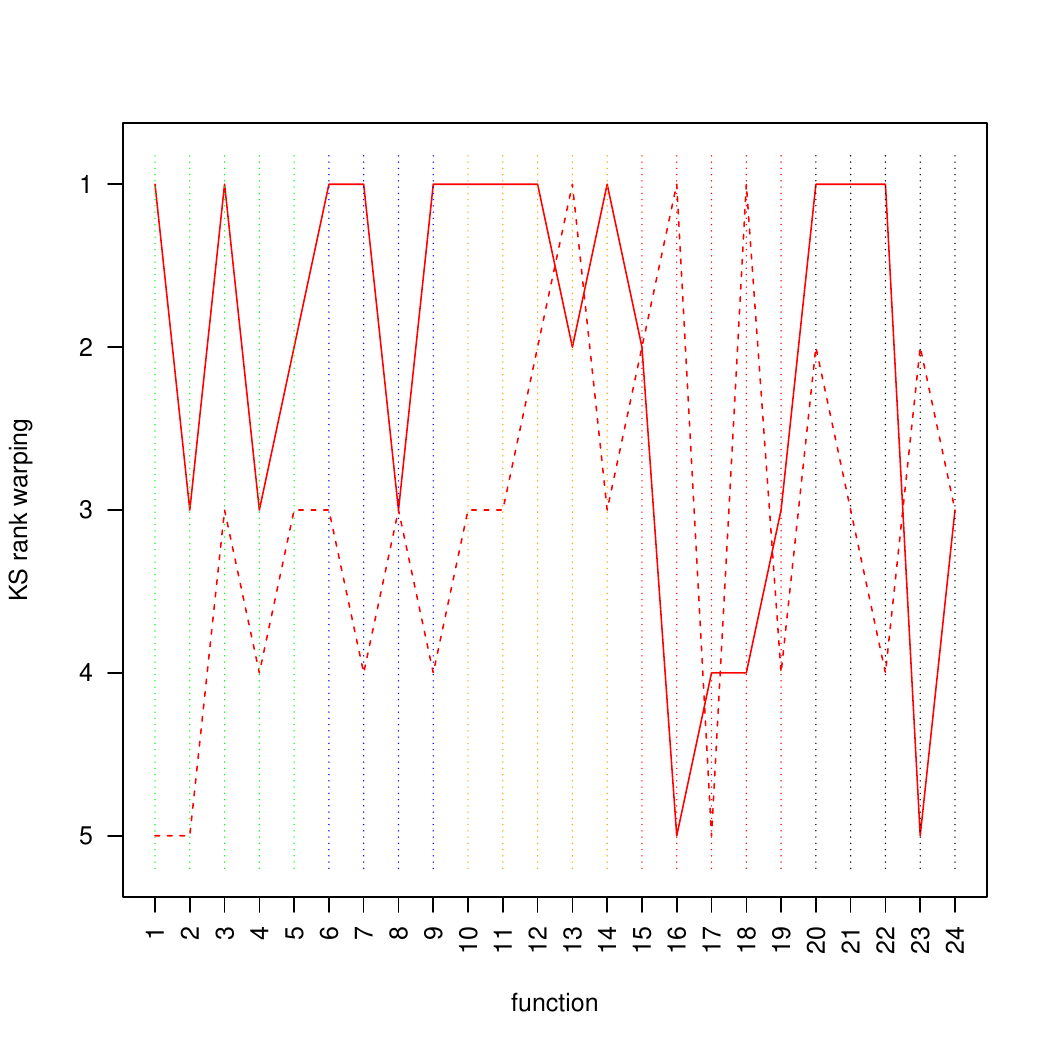}
\caption{
Comparison of the regression (continuous line) and optimization (dashed line) performances, as relative ranks, of the GP (and EGO) variant with warping.
The regression performance is quantified through the p-value of a normality test of the residuals.
\label{fig-pvalPoptrankWarp}
}
\end{figure}

No specific pattern was observed with the p-value of the normality test that could be correlated to either the Q2 or the optimization performance. It was only confirmed, cf. Fig.~\ref{fig-pvalPoptrankWarp}, that warping improves the normality of the residuals, which is precisely its goal. Out of the 24 functions, warping ranked first (best) in terms of residuals normality in 12 functions, far more than its performance in Q2 or $\Popt$.


\section{Complementary results}
\subsection{Initial budget}
\begin{figure}
\mbox{
\includegraphics[width=0.33\textwidth]{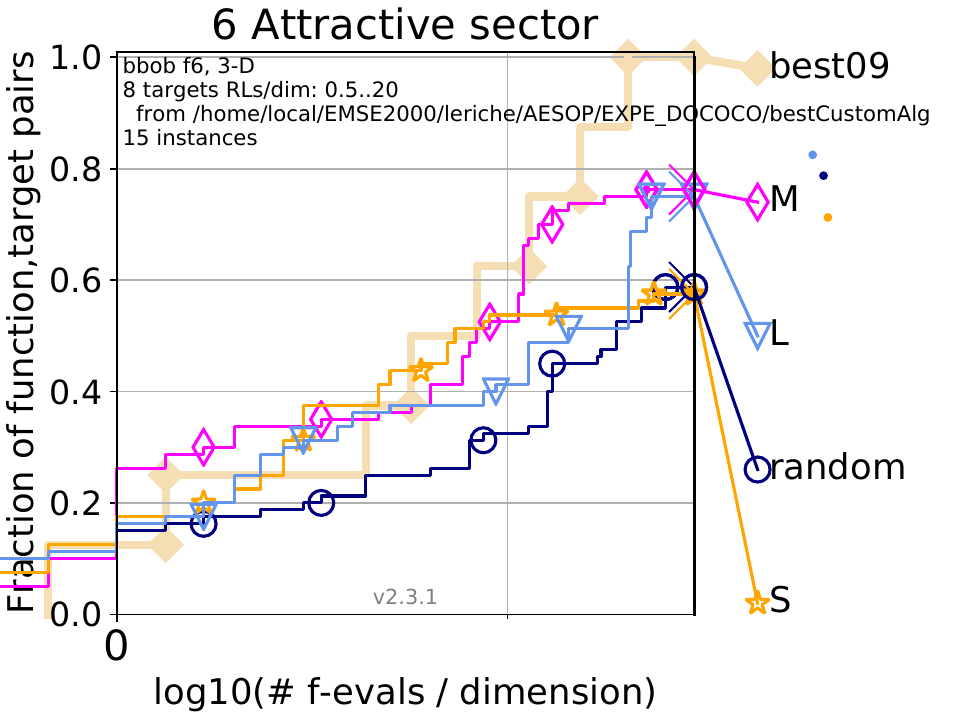}
\includegraphics[width=0.33\textwidth]{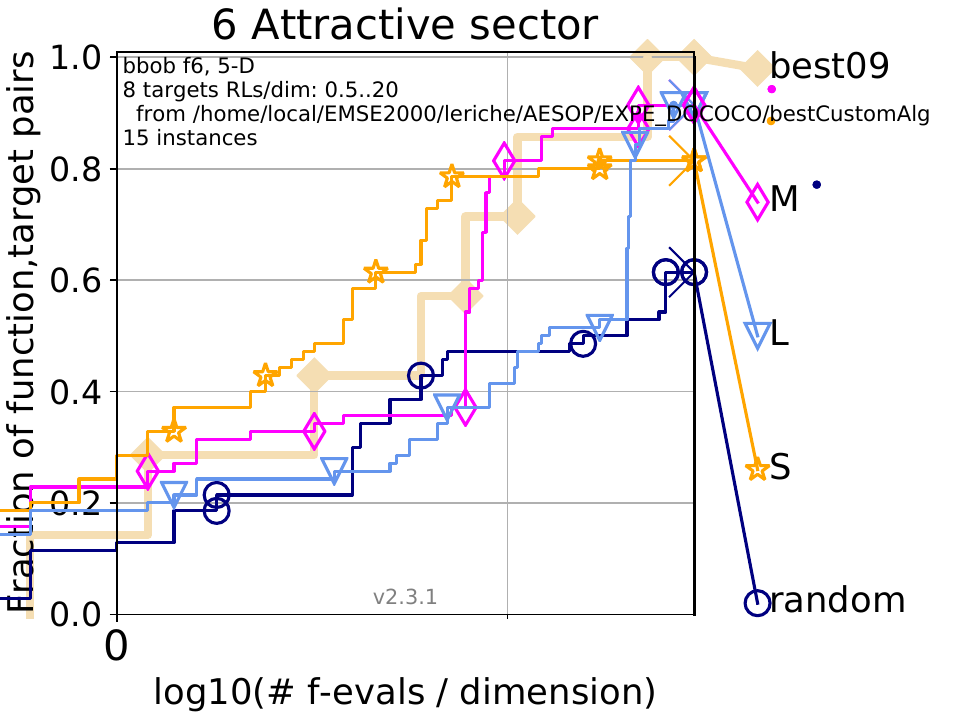}
\includegraphics[width=0.33\textwidth]{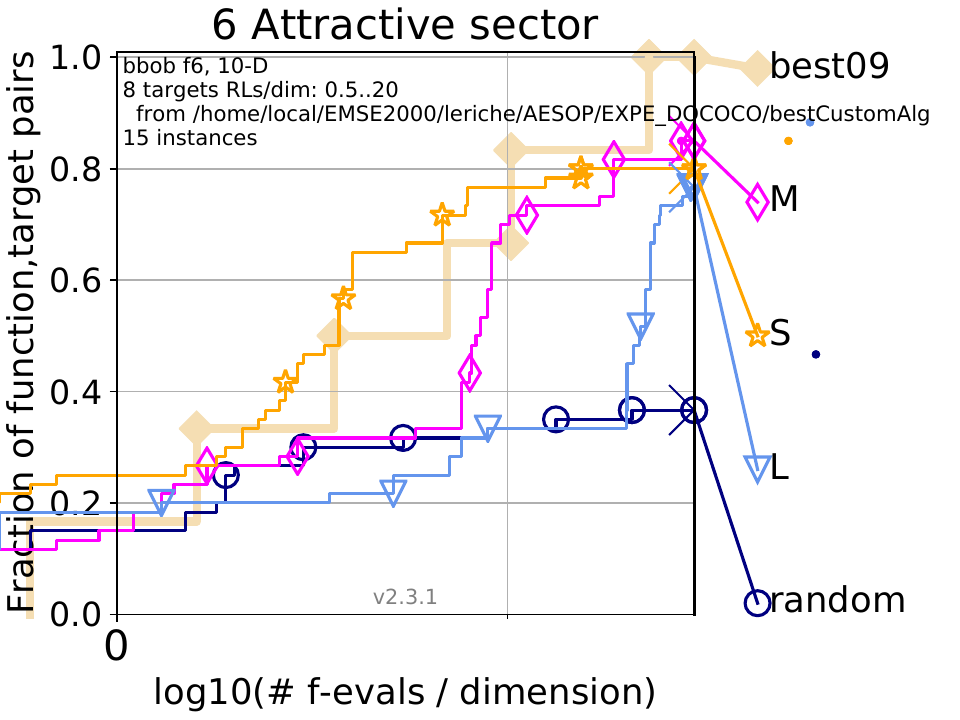}
}
\mbox{
\includegraphics[width=0.33\textwidth]{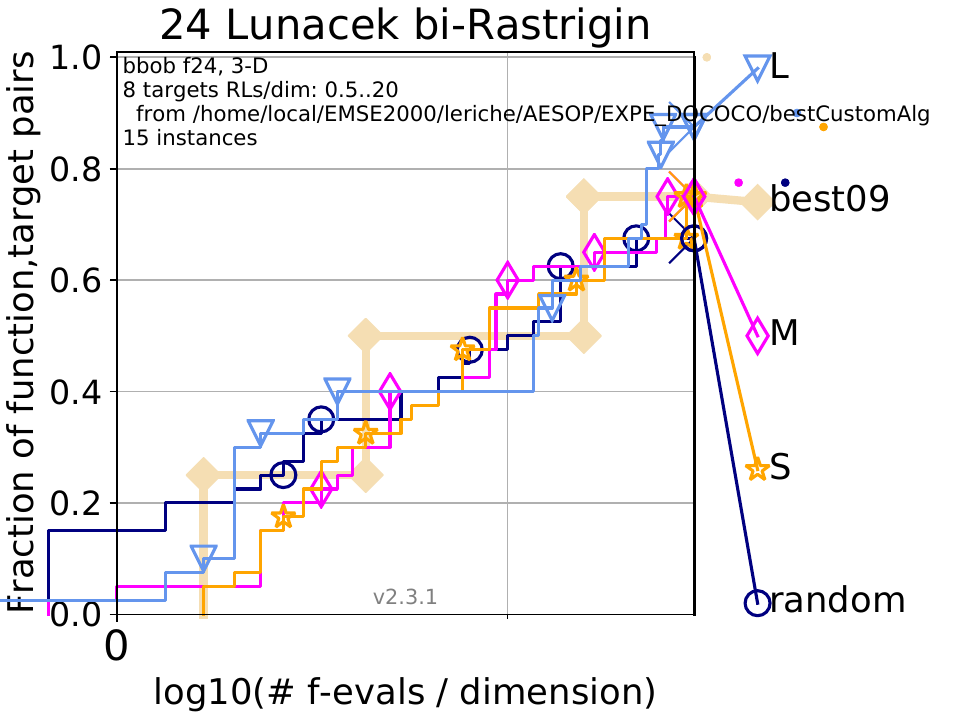}
\includegraphics[width=0.33\textwidth]{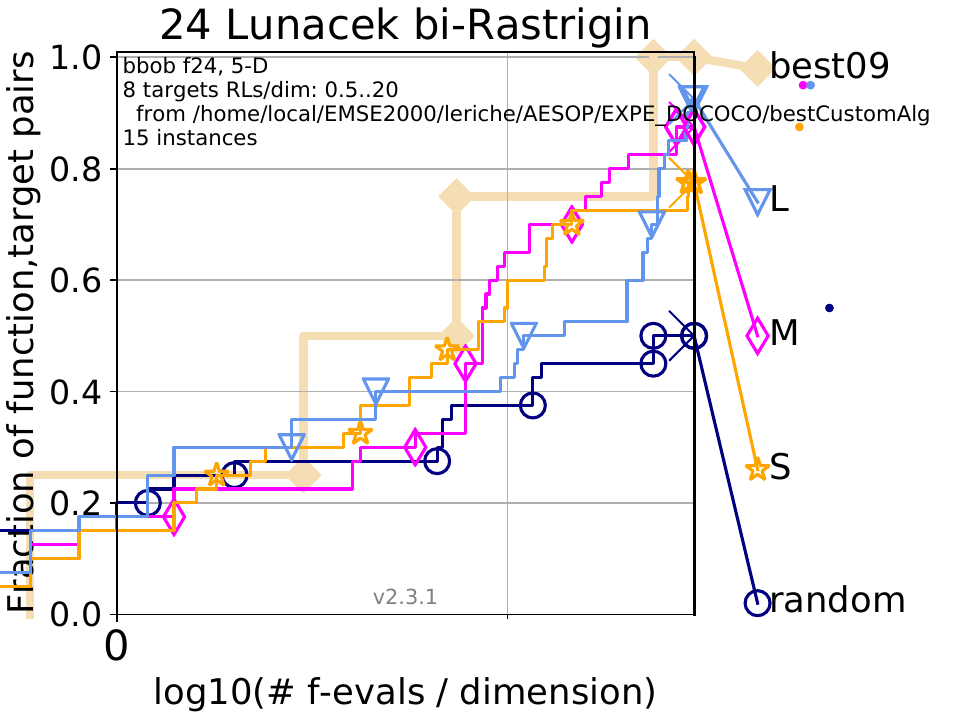}
\includegraphics[width=0.33\textwidth]{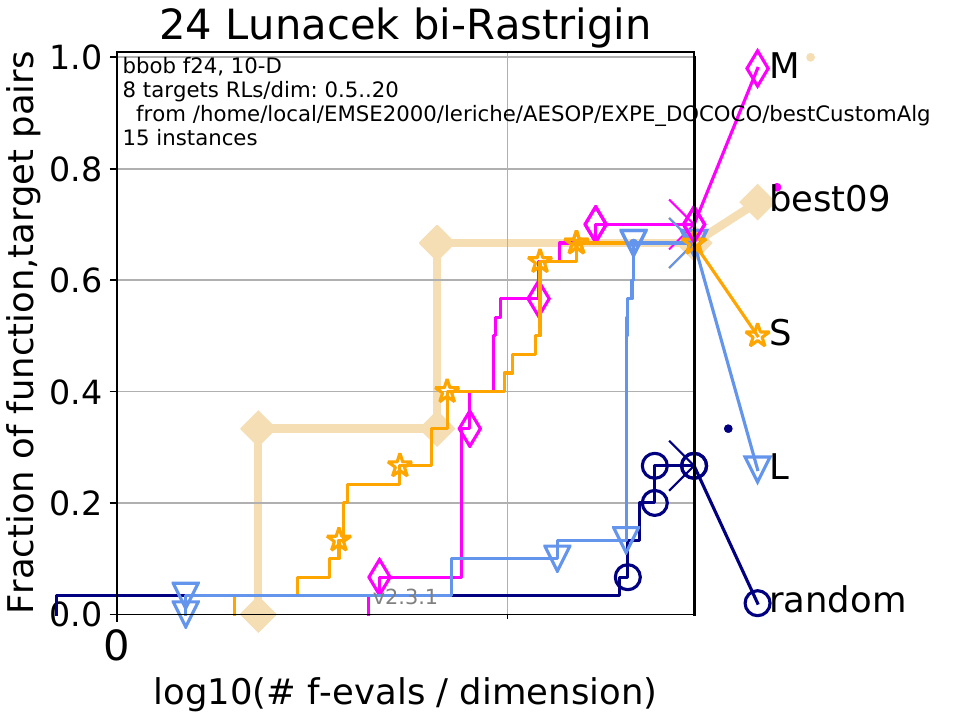}
}
\caption{
Empirical Run Time Distributions (ERTD in $d=3,5$ and 10) of the 2 main exceptions to the rule that a small initial budget is better. 
Top row: Attractive sector function (f6).
Bottom row: Lunacek bi-Rastrigin function (f24).
\label{fig-ERTD_budget_exceptions}
}
\end{figure}

We provide in Fig.~\ref{fig-ERTD_budget_exceptions} performance plots for the 2 only functions that need a large initial sampling 
before starting the optimization. 
Both functions have funnels that need sufficient global exploration to be located.

\subsection{Exponential kernel}
\label{sec-complement_expScalWarp}

\begin{figure}
\mbox{
\includegraphics[width=0.33\textwidth]{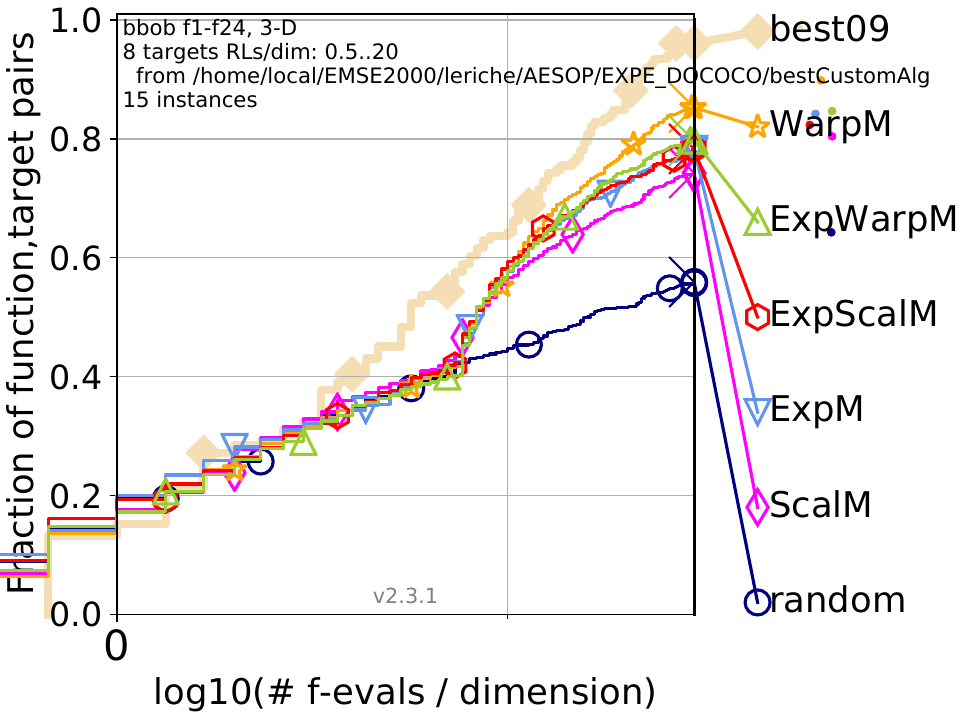}
\includegraphics[width=0.33\textwidth]{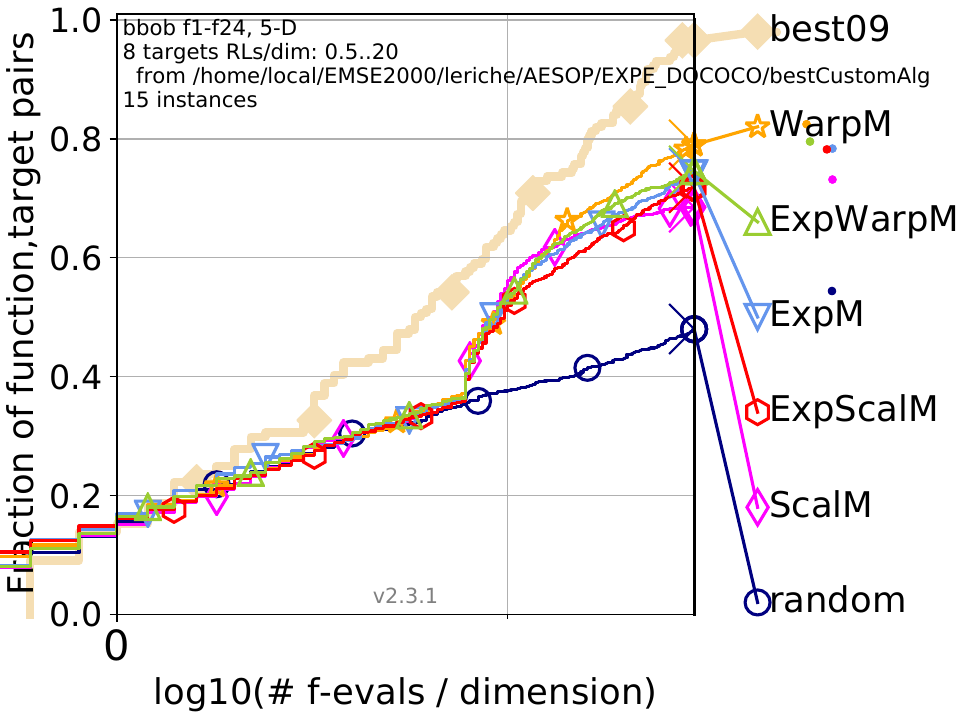}
\includegraphics[width=0.33\textwidth]{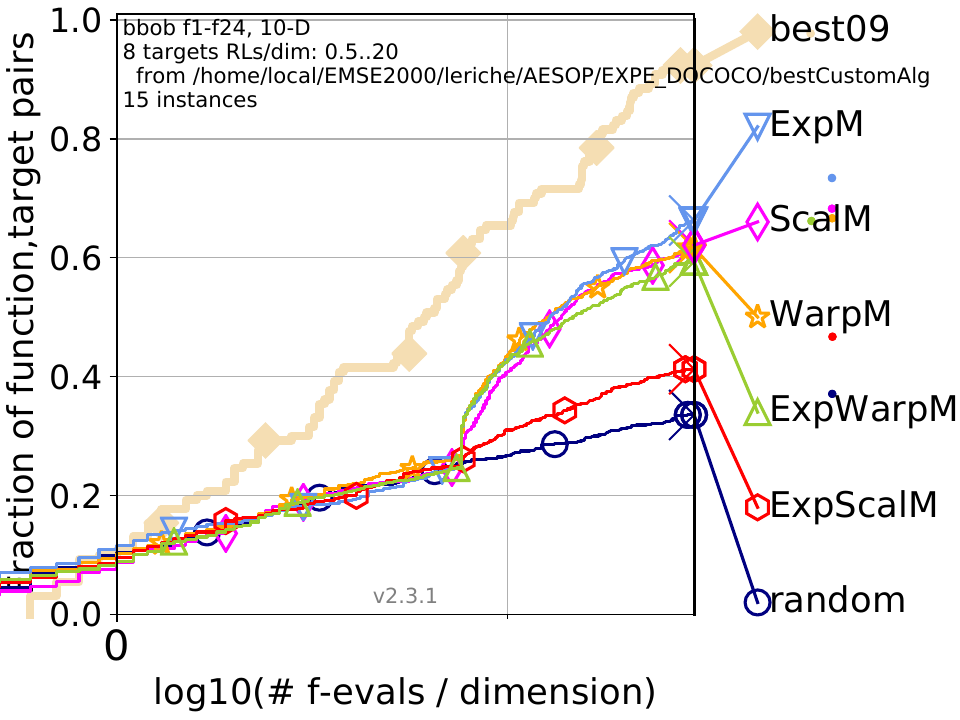}
}
\caption{
Empirical Run Time Distributions (ERTD) over all functions in $d=3,5$ and 10 dimensions testing for potential interactions between the exponential kernel, scaling and warping : algorithms with ``Exp'' in the name use an anisotropic exponential kernel, the others use an anisotropic Mat\'ern 5/2; Scal and Warp correspond to scaling and warping; See Table~\ref{tab:allruns} for further details. 
\label{fig-ERTD_exponScalWarp}
}
\end{figure}

The potential interactions between the exponential kernel and the (input) scaling or the (output) warping have also been investigated for a medium initial DoE, cf. Figure~\ref{fig-ERTD_exponScalWarp}. On the average of all functions, starting in dimension 5 up, no positive interaction is found as ExpM is always better than ExpWarpM which in turn is better than ExpScalM. 

\end{document}